\documentclass[preprint,12pt]{elsarticle}

\usepackage{amssymb}
\usepackage{amsmath}
\usepackage{bm}
\usepackage{color}
\usepackage{exscale}
\usepackage{relsize}
\usepackage{graphicx}
\usepackage{wrapfig}
\usepackage{multicol}
\usepackage{verbatim}

\usepackage{hyperref}
\usepackage{subcaption}
\usepackage{booktabs}
\usepackage{caption}
\usepackage{microtype}

\newtheorem{theorem}{Theorem}
\newtheorem{rmk}{Remark}

\newtheorem{lemma}{Lemma}
\newtheorem{corollary}{Corollary}
\DeclareSymbolFont{yhlargesymbols}{OMX}{yhex}{m}{n}
\DeclareMathAccent{\wideparen}{\mathord}{yhlargesymbols}{"F3}
\let\autocite\citep

\journal{Computers and Mathematics with Applications}

\begin{document}

\begin{frontmatter}



\title{A high-order multi-scale method and its convergence analysis for temperature-dependent nonlinear thermal radiation problems of composite structures}


\author[label1]{Hao Dong\corref{cor1}}\ead{donghao@mail.nwpu.edu.cn}
\cortext[cor1]{Corresponding author.}
\author[label1]{Yongfei Hu}
\author[label1]{Jiale Linghu}
\author[label5]{Yaochuang Han}

\address[label1]{School of Mathematics and Statistics, Xidian University, Xi'an 710071, PR China}
\address[label5]{School of Mathematical Sciences, Luoyang Normal University, Luoyang 471934, PR China}
\begin{abstract}
Accurate prediction of the nonlinear radiation thermal transfer in composite structures with temperature-dependent
properties is significant in high-temperature applications of the materials. This study establishes a high-accuracy
multi-scale computational model incorporating novel high-order correction terms for the high-fidelity simulation of nonlinear thermal radiation in composite structures, enabling local balance preserving of heat quantity. Moreover, an explicit convergence rate is also derived for the resulting high-order multi-scale solutions. Furthermore, an efficient multi-scale algorithm consisting of off-line and on-line computation stages is developed for high-accuracy simulation of nonlinear thermal radiation behavior in composite structures, and corresponding convergence analysis is also obtained. Two- and three-dimensional numerical examples are presented to validate the competitive advantages of the proposed multi-scale approach, not only exceptional numerical accuracy, but also reduced computational cost in both storage requirements and computational time.
\end{abstract}

\begin{keyword}
Temperature-dependent composite structures \sep Nonlinear thermal radiation \sep High-order multi-scale method \sep Efficient multi-scale model \sep Error estimation
\end{keyword}

\end{frontmatter}


\section{Introduction}
In the field of modern aerospace engineering, composite materials and manufactured structures are extensively utilized owing to their excellent thermal performances, such as high temperature resistance, thermal stability and low thermal conductivity. When composite structures are in service in high-temperature environments, the effects of temperature-dependence and thermal radiation must be taken into account due to their significant impact on the service performance of materials~\cite{R1,R2,R3}. For design-oriented evaluation and optimization of composite structures, it is essential to achieve accurate and efficient simulation of the nonlinear thermal radiation problems of composite structures with temperature-dependent properties in high-temperature environments.

Performing direct and high-fidelity numerical simulations for composite structures is prohibitive due to the requirement of an excessively fine mesh. To circumvent this challenging issue, approaches such as the homogenization method and various multi-scale methods have been proposed by scientists and engineers. Early studies about homogenization approaches for heterogeneous composite structures have been extensively pioneered by Cioranescu~\cite{R4}, Babu{\v s}ka~\cite{R5}, De~Giorgi~\cite{R6}, Lewinski~\cite{R7}, Bakhvalov~\cite{R8}, and Oleinik~\cite{R9}, among others, which advanced significantly through the averaging of differential equations with rapidly oscillating coefficients. However, homogenization methods are limited to predicting the macroscopic effective response of composite structures, but fail to resolve the microscopic oscillatory information due to spatial heterogeneity. To accurately simulate the multi-scale behavior of composite structures, various multi-scale solvers have been developed, including Asymptotic homogenization method (AHM)~\cite{R10}, Multi-scale finite element method (MsFEM)~\cite{R11}, Heterogeneous multi-scale method (HMM)~\cite{R12,R13}, Variational multi-scale method (VMS)~\cite{R14}, Multi-scale eigenelement method (MEM)~\cite{R15}, Localized orthogonal decomposition method (LOD)~\cite{R16}, Finite volume based asymptotic homogenization theory (FVBAHT)~\cite{R17}, and Multi-scale generalized finite element method (MS-GFEM)~\cite{R18}. It should be mentioned that existing multi-scale methods still exhibit insufficient numerical accuracy for composite structures with high-contrast material properties~\cite{R19,R20}. Driven by real-world engineering applications, Cui and his research team systematically established a family of high-order multi-scale approaches for high-accuracy and efficient simulation of the thermal, mechanical and multi-physics problems of composite structures, as shown in references \cite{R21,R22,R23,R24,R25,R26,R27,R28} for further details, which can greatly improve the computational efficiency and reduce the consumption of computational resources without loss of numerical accuracy.

For radiative heat transfer in heterogeneous media, Amosov established the semidiscrete problem and an asymptotic approximation model for nonstationary radiative-conductive heat transfer problem in periodic system~\cite{R29,R30}. Bakhvalov~\cite{R31} first proposed an averaging method for the heat transfer process inside periodic media with radiation and obtained the asymptotic expansion of the temperature field. Allaire et al. studied the homogenization method~\cite{R32,R33} and its second-order corrector~\cite{R34} for nonlinear radiative heat transfer problems in heterogeneous domains. Their research demonstrated that the second order corrector is essential to reconstruct the multi-scale solutions. For the radiation heat transfer in composite materials with a periodic microstructure, Huang et al. employed Rothe's method to establish the global existence and uniqueness of nonnegative weak solutions~\cite{R35}, and derived multi-scale asymptotic expansions with an explicit convergence rate $\varepsilon^{1/2}$. Furthermore, the novel multi-scale asymptotic expansions with an explicit convergence rate $\varepsilon$ are established for radiation heat transfer equations with non-homogeneous initial-boundary conditions, and an efficient numerical algorithm is developed for effectively solving these problems~\cite{R36}. Han et al.~\cite{R37} proposed a fast multipole accelerated boundary element method for radiative heat transfer in three-dimensional semitransparent media, whose computational efficiency is superior to that of conventional boundary element method. Nevertheless, existing studies on thermal radiation problems have not considered the nonlinear effect of temperature-dependent material properties. Hence, widespread engineering demands strongly necessitate continued research into this challenging issue, especially for nonlinear radiation thermal transfer in composite structures with temperature-dependent properties in high-temperature environment.

The rest of this study is organized as follows. Section~2 presents the detailed construction of the high-order multi-scale (HOMS) computational model for temperature-dependent nonlinear thermal radiation problems of composite structures, especially for novel high-order correction terms. Section~3 conducts both local and global error analyses for the proposed multi-scale solutions, and an explicit convergence rate is derived. Section~4 develops a two-stage numerical algorithm comprising off-line micro-scale computation, and on-line macro-scale and multi-scale computations. Moreover, corresponding convergence analysis is reported in Section~5. Section~6 presents extensive numerical experiments to verify the computational performance and theoretical results of the proposed approach. Finally, meaningful conclusions and potential directions are summarized in Section~7.

Throughout this study, the Einstein summation convention is applied to streamline repetitive indices.

\section{High-order multi-scale asymptotic analysis}
\subsection{The statement and setting of multi-scale nonlinear thermal radiation problem}
In this paper, the following nonlinear thermal radiation equation with rapidly oscillating coefficients is studied, which originates from the study of radiative heat transfer within composite media with temperature-dependent properties.
\begin{equation}
\left\{ \begin{aligned}
&\rho^{\varepsilon}({\bm{x}},T^{{\varepsilon}})c^{\varepsilon}({\bm{x}},T^{{\varepsilon}})\frac{\partial T^{{\varepsilon}}(\bm{x},t)}{\partial {t}} - \frac{\partial }{\partial {x_i}}\Big( {k_{ij}^{{\varepsilon}}({\bm{x}},T^{{\varepsilon}})\frac{\partial T^{{\varepsilon}}(\bm{x},t)}{\partial {x_j}}} \Big)\\
&\quad\quad\quad\quad\quad\quad+\beta^{{\varepsilon}}({\bm{x}},T^{{\varepsilon}})\sigma_{\mathrm{B}}\big[T^{{\varepsilon}}(\bm{x},t)\big]^4=h(\bm{x},t),\;\;\text{in}\;\;\Omega\times(0,t^*),\\
&T^{{\varepsilon}}(\bm{x},t) = \hat T(\bm{x},t),\;\;\text{on}\;\;\partial {\Omega_T}\times(0,t^*),\\
&{k_{ij}^{{\varepsilon}}({\bm{x}},T^{{\varepsilon}})\frac{\partial T^{{\varepsilon}}(\bm{x},t)}{\partial {x_j}}} {n_i} = \hat q(\bm{x},t),\;\;\text{on}\;\;\partial {\Omega_q}\times(0,t^*),\\
&T^{{\varepsilon}}(\bm{x},0) = \widetilde T(\bm{x}),\;\;\text{on}\;\;{\Omega},
\end{aligned} \right.
\label{eq:all}
\end{equation}
where ${\Omega}$ is a bounded convex domain in $\mathbb{R}^d(d=2,3)$ with the Lipschitz continuous boundary
$\partial\Omega=\partial {\Omega_T}\cup\partial {\Omega_q}$, and assume $\partial {\Omega_T}$ is of positive measure. Moreover, ${\Omega}$ is a composite solid medium, which physically is composed of different materials and has a periodically microscopic configuration with period parameter $\varepsilon$. In the mathematical model \eqref{eq:all}, the temperature field $T^{{\varepsilon}}$ is the unknown to be determined. $\rho^{\varepsilon}({\bm{x}},T^{{\varepsilon}})$, $c^{\varepsilon}({\bm{x}},T^{{\varepsilon}})$, $k_{ij}^{\varepsilon}({\bm{x}},T^{{\varepsilon}})$ and
$\beta^{{\varepsilon}}({\bm{x}},T^{{\varepsilon}})$ denote, respectively, the mass density, the specific heat, the thermal conductivity tensor and the absorption coefficient, all of which are temperature-dependent. $\sigma_{\mathrm{B}}$ is the Stefan--Boltzmann constant. $h(\bm{x},t)$, $\hat T(\bm{x},t)$, $\hat q(\bm{x},t)$ and $\widetilde T(\bm{x})$ denote the internal heat source in $\Omega$, the prescribed temperature on $\partial {\Omega_T}$, the prescribed heat flux
on $\partial {\Omega_q}$ and the initial temperature of $\Omega$.

In accordance with the asymptotic homogenization framework, we denote the micro-scale coordinates
$\bm{y}\!=\!{\bm{x}}/{\varepsilon}\!=\!(x_1/\varepsilon,\cdots,x_d/\varepsilon)\!=\!(y_1,\cdots,y_d)$ of periodic
unit cell (PUC) $\Omega_{\bm{y}}=(0,1)^d$. As a result, $\rho^{\varepsilon}({\bm{x}},T^{{\varepsilon}})$, $c^{\varepsilon}({\bm{x}},T^{{\varepsilon}})$, $k_{ij}^{\varepsilon}({\bm{x}},T^{{\varepsilon}})$ and $\beta^{{\varepsilon}}({\bm{x}},T^{{\varepsilon}})$ can be formulated in new expressions $\rho({\bm{y}},T^{{\varepsilon}})$, $c({\bm{y}},T^{{\varepsilon}})$, $k_{ij}({\bm{y}},T^{{\varepsilon}})$ and $\beta({\bm{y}},T^{{\varepsilon}})$, which entail that these material parameters are 1-periodic functions in micro-variable $\bm{y}$.

According to previous works \cite{R32,R33,R34,R35,R36,R38}, some assumptions are presented for the governing equation \eqref{eq:all} as below.
\begin{enumerate}
\item[(A$_1$)]
Function $k_{ij}^{\varepsilon}({\bm{x}},T^{{\varepsilon}})$ is symmetric and satisfies uniform ellipticity condition such that
\begin{displaymath}
\begin{aligned}
&k_{ij}^{\varepsilon}({\bm{x}},T^{{\varepsilon}})=k_{ji}^{\varepsilon}({\bm{x}},T^{{\varepsilon}}),\;\underline{k}|\bm{\xi}|^2\leq k_{ij}^{\varepsilon}({\bm{x}},T^{{\varepsilon}})\xi_i\xi_j  \leq\overline{k}|\bm{\xi}|^2,
\end{aligned}
\end{displaymath}
where $\underline{k}$ and $\overline{k}$ are two positive constants irrespective of $\varepsilon$ for arbitrary vector $\bm{\xi}=(\xi_1,\cdots,\xi_d)\in \mathbb{R}^d$.
\item[(A$_2$)]
$\rho^{\varepsilon}({\bm{x}},T^{{\varepsilon}})$, $c^{\varepsilon}({\bm{x}},T^{{\varepsilon}})$ and $\beta^{{\varepsilon}}({\bm{x}},T^{{\varepsilon}})\in L^\infty (\Omega)$; $\underline{\rho}\leq\rho^{\varepsilon}({\bm{x}},T^{{\varepsilon}})\leq \overline{\rho}$, $\underline{c}\leq c^{\varepsilon}({\bm{x}},T^{{\varepsilon}})\leq \overline{c}$ and $\underline{\beta}\leq \beta^{\varepsilon}({\bm{x}},T^{{\varepsilon}})\leq \overline{\beta}$, where $\underline{\rho}$, $\overline{\rho}$, $\underline{c}$, $\overline{c}$, $\underline{\beta}$ and $\overline{\beta}$ are positive constants irrespective of $\varepsilon$.
\item[(A$_3$)]
For any $({\bm{y}},T^{{\varepsilon}})\in \Omega_{\bm{y}}\times[T_{min},T_{max}-\delta]$, assume that for all $\delta$ with $0<\delta<T_{max}-T_{min}$
\begin{displaymath}
\begin{aligned}
|k_{ij}({\bm{y}},T^{{\varepsilon}}+\delta)-k_{ij}({\bm{y}},T^{{\varepsilon}})|\leq C_k\delta,\;\;|\beta({\bm{y}},T^{{\varepsilon}}+\delta)-\beta({\bm{y}},T^{{\varepsilon}})|\leq C_\beta\delta,
\end{aligned}
\end{displaymath}
where $C_k$ and $C_\beta$ are positive constants irrespective of $\varepsilon$ and $\delta$.
\item[(A$_4$)]
$h(\bm{x},t)\in H^{2,1}(\Omega\times(0,t^*)),\;\hat{T}(\bm{x},t)\in H^{\frac{7}{2},\frac{7}{4}}(\Omega\times(0,t^*)),\;\hat{q}(\bm{x},t)\in H^{\frac{5}{2},\frac{5}{4}}(\Omega\times(0,t^*)),\;\widetilde T(\bm{x})\in H^3(\Omega)$.
\end{enumerate}

\subsection{The establishment of high-order multi-scale computational model}\label{HOMS}
Considering the definitions of macroscopic and microscopic coordinates, the following rule of chain differentiation shall be derived for the multi-scale nonlinear system.
\begin{equation}
\frac{\partial \Phi^{\varepsilon}(\bm{x},t)}{\partial x_i}=\frac{\partial \Phi(\bm{x},\bm{y},t)}{\partial x_i}+\frac{1}{\varepsilon}\frac{\partial\Phi(\bm{x},\bm{y},t)}{\partial y_i},
\label{eq:chain}
\end{equation}
which will be widely used in the sequel.

Next, in a standard multi-scale asymptotic analysis framework, we define that the exact solution $T^{{\varepsilon}}(\bm{x},t)$ has the asymptotic expansion form as below:
\begin{equation}
\begin{aligned}
T^{{\varepsilon}}(\bm{x},t) = {T_0}(\bm{x},\bm{y},t) + \varepsilon {T_1}(\bm{x},\bm{y},t) + {\varepsilon^2}{T_2}(\bm{x},\bm{y},t) + \mathrm{O}({\varepsilon^3}).
\end{aligned}
\end{equation}

After that, employing Taylor's formula, the temperature-dependent multi-scale function
$\Psi^{\varepsilon}(\bm{x}, T^{\varepsilon})$ admits the following asymptotic expansion.
\begin{align*}
\Psi^{\varepsilon}(\bm{x}, T^{\varepsilon})
&= \Psi(\bm{y}, T^{\varepsilon}) \\
&= \Psi(\bm{y}, T_0 + \varepsilon T_1 + \varepsilon^2 T_2 + \mathrm{O}(\varepsilon^3)) \\
&= \Psi(\bm{y}, T_0)
   + \varepsilon T_1 \mathbf{D}^{(0,1)} \Psi(\bm{y}, T_0) \\
&\quad + \varepsilon^2 \left[ T_2 \mathbf{D}^{(0,1)} \Psi(\bm{y}, T_0)
   + \frac{1}{2} T_1^2 \mathbf{D}^{(0,2)} \Psi(\bm{y}, T_0) \right]
   + \mathrm{O}(\varepsilon^3) \\
&= \Psi^{(0)}(\bm{y}, T_0)
   + \varepsilon \Psi^{(1)}(\bm{x}, \bm{y}, T_0)
   + \varepsilon^2 \Psi^{(2)}(\bm{x}, \bm{y}, T_0)
   + \mathrm{O}(\varepsilon^3),
\end{align*}
where the multi-index notations are defined as $f_y(x,y) = \mathbf{D}^{(0,1)} f(x,y)$ and
$f_{yy}(x,y) = \mathbf{D}^{(0,2)} f(x,y)$. Hence, the nonlinear material property functions
$\rho^{\varepsilon}(\bm{x}, T^{\varepsilon})$, $c^{\varepsilon}(\bm{x}, T^{\varepsilon})$,
$k_{ij}^{\varepsilon}(\bm{x}, T^{\varepsilon})$ and $\beta^{\varepsilon}(\bm{x}, T^{\varepsilon})$ admit the following expansion expressions.
\begin{equation}
	\begin{aligned}
		{\rho ^{{\varepsilon}}}(\bm{x},{T^{{\varepsilon}}}) = {\rho ^{(0)}}(\bm{y},{T_0}) +
    {\varepsilon}{\rho ^{(1)}}(\bm{x},\bm{y},{T_0}) + \varepsilon^2{\rho ^{(2)}}
    (\bm{x},\bm{y},{T_0}) + {\mathrm{O}}(\varepsilon^3),
	\end{aligned}
  \label{eq:3}
\end{equation}
\begin{equation}
	\begin{aligned}		
		{c^{{\varepsilon}}}(\bm{x},{T^{{\varepsilon}}}) = {c^{(0)}}(\bm{y},{T_0}) +
    {\varepsilon}{c^{(1)}}(\bm{x},\bm{y},{T_0}) + \varepsilon^2{c^{(2)}}(\bm{x},\bm{y},{T_0})
    + {\mathrm{O}}(\varepsilon^3),
	\end{aligned}
\end{equation}
\begin{equation}
	\begin{aligned}
		k_{ij}^{{\varepsilon}}(\bm{x},{T^{{\varepsilon}}}) = k_{ij}^{(0)}(\bm{y},{T_0}) +
    {\varepsilon}k_{ij}^{(1)}(\bm{x},\bm{y},{T_0}) + \varepsilon^2k_{ij}^{(2)}
    (\bm{x},\bm{y},{T_0}) + {\mathrm{O}}(\varepsilon^3),
	\end{aligned}
\end{equation}
\begin{equation}
	\begin{aligned}		
		{\beta^{{\varepsilon}}}(\bm{x},{T^{{\varepsilon}}}) = {\beta^{(0)}}(\bm{y},{T_0}) + {\varepsilon}{\beta ^{(1)}}(\bm{x},\bm{y},{T_0}) + \varepsilon^2{\beta^{(2)}}(\bm{x},\bm{y},{T_0}) + {\mathrm{O}}(\varepsilon^3).
	\end{aligned}
  \label{eq:3-7}
\end{equation}

Then substituting \eqref{eq:3}-\eqref{eq:3-7} into multi-scale initial-boundary value system \eqref{eq:all} and utilizing
the chain rule \eqref{eq:chain}, a hierarchy of equations is obtained by collecting terms at each power of the small parameter $\varepsilon$ as follows.
\begin{equation}
O(\varepsilon^{-2}):
\begin{aligned}
\frac{\partial }{\partial y_i}\Big(k_{ij}^{(0)}(\bm{y},{T_0})\frac{\partial T_0}{\partial y_j} \Big) = 0.
\end{aligned}\label{eq:epsilon-2}
\end{equation}
\begin{equation}
O(\varepsilon^{ - 1}):
\begin{aligned}
&\frac{\partial }{\partial y_i}\Big( {{k_{ij}^{(1)}}(\bm{y},{T_0})\frac{\partial T_0}{\partial y_j}} \Big)+\frac{\partial }{\partial y_i}\Big( {k_{ij}^{(0)}(\bm{y},{T_0})\big( {\frac{\partial T_0}{\partial {x_j}} + \frac{\partial T_1}{\partial y_j}} \big)} \Big)\\
&+\frac{\partial }{\partial {x_i}}\Big( {k_{ij}^{(0)}(\bm{y},{T_0})\frac{\partial {T_0}}{\partial y_j}} \Big)= 0.
\end{aligned}\label{eq:epsilon-1}
\end{equation}
\begin{equation}
O(\varepsilon^0):
\begin{aligned}
&\rho^{(0)}(\bm{y},{T_0})c^{(0)}(\bm{y},{T_0})\frac{\partial T_0}{\partial {t}}-\frac{\partial }{\partial y_i}\Big( {{k_{ij}^{(2)}(\bm{y},{T_0})}\frac{\partial T_0}{\partial y_j}} \Big)\\
&- \frac{\partial }{\partial {x_i}}\Big( {k_{ij}^{(1)}(\bm{y},{T_0})\frac{\partial {T_0}}{\partial y_j}} \Big)-\frac{\partial }{\partial {y_i}}\Big( {k_{ij}^{(1)}(\bm{y},{T_0})\big( {\frac{\partial {T_0}}{\partial {x_j}} + \frac{\partial {T_1}}{\partial y_j}} \big)} \Big)\\
&\!-\!\frac{\partial }{\partial y_i}\Big( {k_{ij}^{(0)}(\bm{y},{T_0})\big( {\frac{\partial {T_1}}{\partial {x_j}}\! + \!\frac{\partial {T_2}}{\partial y_j}} \big)} \Big)\!\!-\! \frac{\partial }{\partial {x_i}}\Big( {k_{ij}^{(0)}(\bm{y},{T_0})\big( {\frac{\partial {T_0}}{\partial {x_j}}\! + \!\frac{\partial {T_1}}{\partial y_j}} \big)} \Big)\\
&+\beta^{{(0)}}(\bm{y},{T_0})\sigma_{\mathrm{B}}\big[T_0\big]^4=h.
\end{aligned}\label{eq:epsilon0}
\end{equation}

On account of $O(\varepsilon^{-2})$-order equation \eqref{eq:epsilon-2}, we can certainly conclude that $T_0(\bm{x},\bm{y},t)$ is independent of the microscopic variable $\bm{y}$, namely:
\begin{equation}
T_0(\bm{x},\bm{y},t) =T_0(\bm{x},t).
\label{eq:T0}
\end{equation}

And then, with the help of \eqref{eq:T0}, the term $\displaystyle\frac{\partial T_0}{\partial y_j}$ in $O(\varepsilon^{-1})$-order equation \eqref{eq:epsilon-1} equals zero. From this, we thus define the concrete expression for $T_1(\bm{x},\bm{y},t)$ as below:
\begin{equation}
T_1({\bm{x}},{{\bm{y}},t}) = M_{{\alpha _1}}(\bm{y},{T_0})\frac{\partial {T_0(\bm{x},t)}}{\partial {x_{{\alpha _1}}}},
\end{equation}
where $M_{{\alpha _1}}(\bm{y},{T_0})$ is defined as the first-order auxiliary cell functions defined in unit cell $\Omega_{\bm{y}}$, which satisfies the following first-order auxiliary cell problem.
\begin{equation}
\left\{ \begin{aligned}
&\frac{\partial }{\partial y_i}\Big( k_{ij}^{(0)}\frac{\partial {M_{\alpha _1}(\bm{y},{T_0})}}{{\partial y_j}} \Big) =-\frac{\partial k_{i\alpha_1}^{(0)}(\bm{y},{T_0})}{\partial y_i},\;\;\;{\bm{y}} \in {\Omega_{\bm{y}}},\\
&M_{{\alpha _1}}(\bm{y},{T_0})\;\mathrm{is}\;1-\mathrm{periodic}\;\mathrm{in}\;{\bm{y}}{\rm{, }}\;\;\;\int_{\Omega_{\bm{y}}}M_{{\alpha _1}}\mathrm{d}\Omega_{\bm{y}}=0.
\end{aligned} \right.
\label{eq:Malpha1}
\end{equation}
\begin{rmk}
On the basis of the Lax-Milgram theorem and the assumption (A$_1$), it can be rigorously proved that
first-order auxiliary cell problem \eqref{eq:Malpha1} has a unique solution for any fixed macroscopic temperature value $T_0$.
\end{rmk}

Subsequently, performing a volume integral of $O(\varepsilon^{0})$-order equation \eqref{eq:epsilon0} over microscopic unit cell $\Omega_{\bm{y}}$ and applying the Gauss theorem, we immediately derive the macroscopic homogenized equation associated with multi-scale nonlinear system \eqref{eq:all} as below:
\begin{equation}
\left\{ \begin{aligned}
&\bar S({T_0})\frac{\partial T_0(\bm{x},t)}{\partial {t}} - \frac{\partial }{{\partial {x_i}}}\Big( {\bar k_{ij}({T_0})\frac{{\partial {T_0}(\bm{x},t)}}{{\partial {x_j}}}} \Big)\\
&\quad\quad\quad\quad\quad\quad+\bar\beta({T_0})\sigma_{\mathrm{B}}\big[T_0(\bm{x},t)\big]^4=h(\bm{x},t),\;\;\text{in}\;\;\Omega\times(0,t^*),\\
&T_0(\bm{x},t) = \hat T(\bm{x},t),\;\;\text{on}\;\;\partial {\Omega_T}\times(0,t^*),\\
&{\bar k_{ij}(T_0)\frac{\partial T_0(\bm{x},t)}{\partial {x_j}}} {n_i} = \hat q(\bm{x},t),\;\;\text{on}\;\;\partial {\Omega_q}\times(0,t^*),\\
&T_0(\bm{x},0) = \widetilde T(\bm{x}),\;\;\text{on}\;\;{\Omega}.
\end{aligned} \right.
\label{eq:homogenized}
\end{equation}
Hereafter, the macroscopic homogenized material parameters in \eqref{eq:homogenized} are summarized as:
\begin{equation}
\begin{aligned}
&\bar S({T_0})=\frac{1}{|\Omega_{\bm{y}}|}{\int_{\Omega_{\bm{y}}}}\rho^{(0)}(\bm{y},{T_0})c^{(0)}(\bm{y},{T_0})\,\mathrm{d}\Omega_{\bm{y}},\\
&{\bar k_{ij}}({T_0}) = \frac{1}{|\Omega_{\bm{y}}|}{\int_{\Omega_{\bm{y}}}}\big({k_{ij}^{(0)}(\bm{y},{T_0})+ k_{ik}^{(0)}(\bm{y},{T_0}){\frac{\partial M_j(\bm{y},{T_0})}{\partial y_k}}}\big)\,\mathrm{d}\Omega_{\bm{y}},\\
&\bar\beta({T_0})=\frac{1}{|\Omega_{\bm{y}}|}{\int_{\Omega_{\bm{y}}}}\beta^{{(0)}}(\bm{y},{T_0})\,\mathrm{d}\Omega_{\bm{y}},
\end{aligned}
\label{eq:homogenized_params}
\end{equation}
where $|\Omega_{\bm{y}}|$ denotes the Lebesgue measure of microscopic unit cell $\Omega_{\bm{y}}$.
\begin{rmk}
All macroscopic homogenized material parameters of the investigated composite structures obviously vary with the macroscopic temperature solution $T_{0}$ because of the quasi-periodic properties of first-order cell functions with respect to $T_{0}$.
\end{rmk}
\begin{rmk}
Following the approach in \cite{R4}, $\bar{S}_0\leq \bar S(T_0)\le\bar{S}_1$, $\bar{\kappa}_0|\bm{\xi}|^2\leq \bar k_{ij}(T_0)\xi_i\xi_j\le\bar{\kappa}_1|\bm{\xi}|^2$, and $\bar{B}_0\leq \bar \beta(T_0)\le\bar{B}_1$, can be proved, where $\bar{S}_0$, $\bar{S}_1$, $\bar{\kappa}_0$, $\bar{\kappa}_1$, $\bar{B}_0$ and $\bar{B}_1$ are positive constants.
\end{rmk}

Additionally, substituting the term $h(\bm{x},t)$ in $O(\varepsilon^{0})$-order equation \eqref{eq:epsilon0} with its macroscopic homogenized counterpart \eqref{eq:homogenized}, we obtain that $T_2(\bm{x},\bm{y},t)$ satisfies the following equality.
\begin{equation}
\begin{aligned}
&\frac{\partial}{\partial y_i}\Big( {k_{ij}^{(0)}(\bm{y},{T_0}){\frac{\partial T_2}{\partial y_j}}} \Big)\!=\!\Big[ { {\bar k}_{\alpha_1\alpha_2}({T_0}) \!- \! k_{\alpha_1\alpha_2}^{(0)}(\bm{y},{T_0})}\!-\!\frac{\partial}{\partial y_i}\big( {k_{i\alpha_1}^{(0)}(\bm{y},{T_0}){M_{\alpha_2}}} \big)\\
&{-{k_{\alpha_1j}^{(0)}(\bm{y},{T_0})\frac{\partial M_{\alpha_2}}{\partial y_j}}} \Big]\frac{\partial^2 T_{0}}{\partial x_{\alpha_1}\partial x_{\alpha_2}}+\Big[\rho^{(0)}(\bm{y},{T_0})c^{(0)}(\bm{y},{T_0})-\bar S({T_0})\Big]\frac{\partial T_0}{\partial {t}}\\
&+\Big[\beta^{{(0)}}(\bm{y},{T_0})\sigma_{\mathrm{B}}-\bar\beta({T_0})\sigma_{\mathrm{B}}\Big]\big[T_0\big]^4+ \Big[ { \frac{\partial{\bar k}_{i\alpha_1}(T_{0})}{\partial x_i} -  {\frac{\partial{k}_{i\alpha_1}^{(0)}(\bm{y},T_{0})}{\partial x_i}}}\\
&{- \frac{\partial}{\partial y_i}\big( {k_{ij}^{(0)}(\bm{y},T_{0})\frac{\partial M_{\alpha_1}}{\partial x_{j}}} \big)-\frac{\partial}{\partial x_i}\big( {k_{ij}^{(0)}(\bm{y},T_{0})\frac{\partial M_{\alpha_1}}{\partial y_{j}}} \big)} \Big]\frac{\partial T_{0}}{\partial x_{\alpha_1}}\\
&- \frac{\partial}{\partial y_i}\Big({{M_{\alpha_1}}\mathbf{D}^{(0,1)}{k_{i\alpha_2}^{(0)}(\bm{y},T_{0})}}+ {{M_{\alpha_1}}\mathbf{D}^{(0,1)}{k_{ij}^{(0)}(\bm{y},T_{0})}\frac{\partial M_{\alpha_2}}{\partial y_j}}\Big)\frac{\partial T_{0}}{\partial x_{\alpha_1}}\frac{\partial T_{0}}{\partial x_{\alpha_2}}.
\end{aligned}\label{eq:T2}
\end{equation}
According to the pivotal equality \eqref{eq:T2}, the detailed expression for $T_2(\bm{x},\bm{y},t)$ can be established as follows.
\begin{equation}
\begin{aligned}
{T_2}({\bm{x}},{{\bm{y}}},t)& ={M_{{\alpha _1}{\alpha _2}}}(\bm{y},T_{0})\frac{{{\partial ^2}{T_0(\bm{x},t)}}}{{\partial {x_{{\alpha _1}}}\partial {x_{{\alpha _2}}}}}+G(\bm{y},T_{0})\frac{\partial T_0(\bm{x},t)}{\partial {t}}\\
&+N(\bm{y},T_{0})\big[T_0(\bm{x},t)\big]^4+{N_{{\alpha _1}}}(\bm{y},T_{0})\frac{{{\partial }{T_0(\bm{x},t)}}}{{\partial {x_{{\alpha _1}}}}}\\
&+{N_{{\alpha _1\alpha _2}}}(\bm{y},T_{0})\frac{\partial T_{0}(\bm{x},t)}{\partial x_{\alpha_1}}\frac{\partial T_{0}(\bm{x},t)}{\partial x_{\alpha_2}},
\end{aligned}\label{eq:T2expression}
\end{equation}
where $M_{\alpha_1\alpha_2}(\bm{y},T_{0})$, $G(\bm{y},T_{0})$, $N(\bm{y},T_{0})$, $N_{\alpha_1}(\bm{y},T_{0})$ and $N_{\alpha_1\alpha_2}(\bm{y},T_{0})$ are denoted as the second-order auxiliary cell functions defined on PUC $\Omega_{\bm{y}}$.

By substituting \eqref{eq:T2expression} into \eqref{eq:T2}, a sequence of equations \eqref{eq:Malpha1alpha2}-\eqref{eq:Nalpha1alpha2}, which are imposed with the homogeneous Dirichlet boundary condition on outer boundary of unit cell $\Omega_{\bm{y}}$, are defined in turn for solving second-order auxiliary cell functions as follows:
\begin{equation}
\left\{
\begin{aligned}
&\frac{\partial}{\partial y_i}\big[ k_{ij}^{(0)}\frac{\partial M_{\alpha_1\alpha_2}}{\partial y_j} \big]\!\!
=\!\bar k_{\alpha_1\alpha_2}\!-\!k_{\alpha_1\alpha_2}^{(0)}\!-\!\frac{\partial}{\partial y_i}\bigl(k_{i\alpha_1}^{(0)}M_{\alpha_2}\bigr)
\!-\!k_{\alpha_1j}^{(0)}\frac{\partial M_{\alpha_2}}{\partial y_j},\bm{y}\!\in\!\Omega_{\bm{y}},\\
&M_{{\alpha _1\alpha_2}}(\bm{y},{T_0})\;\mathrm{is}\;1-\mathrm{periodic}\;\mathrm{in}\;{\bm{y}}{\rm{, }}\;\;\;\int_{\Omega_{\bm{y}}}M_{{\alpha _1\alpha_2}}\mathrm{d}\Omega_{\bm{y}}=0.
\end{aligned}
\right.
\label{eq:Malpha1alpha2}
\end{equation}
\begin{equation}
\left\{
\begin{aligned}
&\frac{\partial}{\partial y_i}\big[ { k_{ij}^{(0)}{\frac{\partial G}{\partial y_j}}} \big] = \rho^{(0)} {c}^{(0)}-\bar S,\;\;\;\bm{y}\in \Omega_{\bm{y}},\\
&G(\bm{y},{T_0})\;\mathrm{is}\;1-\mathrm{periodic}\;\mathrm{in}\;{\bm{y}}{\rm{, }}\;\;\;\int_{\Omega_{\bm{y}}}G\mathrm{d}\Omega_{\bm{y}}=0.
\end{aligned} \right.
\label{eq:G}
\end{equation}
\begin{equation}
\left\{
\begin{aligned}
&\frac{\partial}{\partial y_i}\big[ { k_{ij}^{(0)}{\frac{\partial N}{\partial y_j}}} \big] =\beta^{{(0)}}\sigma_{\mathrm{B}}-\bar\beta\sigma_{\mathrm{B}},\;\;\;\bm{y}\in \Omega_{\bm{y}},\\
&N(\bm{y},{T_0})\;\mathrm{is}\;1-\mathrm{periodic}\;\mathrm{in}\;{\bm{y}}{\rm{, }}\;\;\;\int_{\Omega_{\bm{y}}}N\mathrm{d}\Omega_{\bm{y}}=0.
\end{aligned} \right.
\label{eq:N}
\end{equation}
\begin{equation}
\!\!\!\!\!\!\!\left\{
\begin{aligned}
&\frac{\partial}{\partial y_i}\!\Bigl[ k_{ij}^{(0)}\frac{\partial N_{\alpha_1}}{\partial y_j} \!\Bigr]\!
\!=\!\frac{\partial\bar k_{i\alpha_1}}{\partial x_i}
\!-\!\frac{\partial k_{i\alpha_1}^{(0)}}{\partial x_i}
\!-\!\frac{\partial}{\partial y_i}\!\Bigl(\! k_{ij}^{(0)}\!\frac{\partial M_{\alpha_1}}{\partial x_j} \!\Bigr)
\!-\!\frac{\partial}{\partial x_i}\!\Bigl(\! k_{ij}^{(0)}\!\frac{\partial M_{\alpha_1}}{\partial y_j} \!\Bigr),\! \bm{y}\in \Omega_{\bm{y}},\\
&N_{\alpha_1}(\bm{y},{T_0})\;\mathrm{is}\;1-\mathrm{periodic}\;\mathrm{in}\;{\bm{y}}{\rm{, }}\;\;\;\int_{\Omega_{\bm{y}}} N_{\alpha_1}\mathrm{d}\Omega_{\bm{y}}=0.
\end{aligned}
\right.
\label{eq:Nalpha1}
\end{equation}
\begin{equation}
\!\!\!\!\!\left\{
\begin{aligned}
&\frac{\partial}{\partial y_i}\!\Bigl[
k_{ij}^{(0)}
\frac{\partial N_{\alpha_1\alpha_2}}{\partial y_j}
\Bigr]
\!=\!
-\frac{\partial}{\partial y_i}\!\Bigl(\!
M_{\alpha_1}\mathbf{D}^{(0,1)}k_{i\alpha_2}^{(0)}
\!+\!
M_{\alpha_1}\mathbf{D}^{(0,1)}k_{ij}^{(0)}
\frac{\partial M_{\alpha_2}}{\partial y_j}
\!\Bigr),\!
\bm{y}\in\Omega_{\bm{y}},\\
&N_{\alpha_1\alpha_2}(\bm{y},{T_0})\;\mathrm{is}\;1-\mathrm{periodic}\;\mathrm{in}\;{\bm{y}}{\rm{, }}\;\;\;\int_{\Omega_{\bm{y}}} N_{\alpha_1\alpha_2}\mathrm{d}\Omega_{\bm{y}}=0.
\end{aligned}
\right.
\label{eq:Nalpha1alpha2}
\end{equation}
\begin{rmk}
Under assumption (A$_1$), the Lax-Milgram theorem ensures the existence and uniqueness of solutions to the auxiliary cell problems \eqref{eq:Malpha1alpha2}-\eqref{eq:Nalpha1alpha2} for any fixed macroscopic temperature value $T_0$.
\end{rmk}

Summarizing, the following high-order multi-scale solution can be concluded for the multi-scale nonlinear problem \eqref{eq:all}.
\begin{equation}
\begin{aligned}
T^{\varepsilon}(\bm{x},t) &\approx T_{0}(\bm{x},t) + \varepsilon M_{\alpha_1}(\bm{y},T_0) \frac{\partial T_0(\bm{x},t)}{\partial x_{\alpha_1}} \\
&\quad + \varepsilon^2 \bigg[ M_{\alpha_1\alpha_2}(\bm{y},T_0) \frac{\partial^2 T_0(\bm{x},t)}{\partial x_{\alpha_1} \partial x_{\alpha_2}}
+ G(\bm{y},T_0) \frac{\partial T_0(\bm{x},t)}{\partial t} \\
&\quad + N(\bm{y},T_0) \big[ T_0(\bm{x},t) \big]^4
+ N_{\alpha_1}(\bm{y},T_0) \frac{\partial T_0(\bm{x},t)}{\partial x_{\alpha_1}} \\
&\quad + N_{\alpha_1\alpha_2}(\bm{y},T_0) \frac{\partial T_0(\bm{x},t)}{\partial x_{\alpha_1}} \frac{\partial T_0(\bm{x},t)}{\partial x_{\alpha_2}}\bigg].
\end{aligned}
\label{eq:HOMS}
\end{equation}

\section{The error analyses of multi-scale asymptotic solutions}
This section first demonstrates the local balance preserving of heat quantity of high-order multi-scale solution by the local error analysis in the pointwise sense. Furthermore, the global error estimate of high-order multi-scale solution is derived in the integral sense, which theoretically guarantees the computational reliability of the proposed high-order multi-scale approach.

\subsection{Error analysis of multi-scale solutions in the pointwise sense}
Firstly, define the low-order multi-scale (LOMS) solution and high-order multi-scale (HOMS) solution for multi-scale nonlinear thermal radiation equation \eqref{eq:all} as below
\begin{equation}
\begin{aligned}
{T^{(1\varepsilon)} }(\bm{x},t)={T_0} + \varepsilon {T_1},\;
{T^{(2\varepsilon)} }(\bm{x},t)={T^{(1\varepsilon)}}+\varepsilon ^2 {T_2}.
\end{aligned}
\end{equation}

Next, the residual functions for LOMS and HOMS approximate solutions are defined as below.
\begin{equation}
\begin{aligned}
{T_{\Delta}^{(1\varepsilon)} }(\bm{x},t)={T^{\varepsilon} }(\bm{x},t)-{T^{(1\varepsilon)} }(\bm{x},t),\;
{T_{\Delta}^{(2\varepsilon)} }(\bm{x},t)={T^{\varepsilon} }(\bm{x},t)-{T^{(2\varepsilon)} }(\bm{x},t).
\end{aligned}
\label{eq:residual_second_order}
\end{equation}

To analyze the local heat quantity balance preserving by multi-scale asymptotic solutions, substituting the above residual function $T^{(1\varepsilon)}_{\Delta}$ into \eqref{eq:all}, the residual equation for LOMS solution is derived as follows.
\begin{equation}
\left\{
\begin{aligned}
&\rho^{\varepsilon}c^{\varepsilon}\frac{{\partial {T_{\Delta}^{(1\varepsilon)}}}}{{\partial t}}- \frac{\partial }{{\partial {x_i}}}\Big( {{k_{ij}^{\varepsilon}}\frac{{\partial {T_{\Delta}^{(1\varepsilon)} }}}{{\partial {x_j}}}}\Big)\\
&\;\;=F_0(\bm{x},\bm{y},t)+\varepsilon F_1(\bm{x},\bm{y},t),\;\;\text{in}\;\;\Omega\times(0,t^*),\\
&{T_{\Delta}^{(1\varepsilon)}(\bm{x},t) } = -\varepsilon{M_{\alpha_1}}\frac{\partial T_{0}}{\partial x_{\alpha_1}}:=\varepsilon \widehat{\psi}_1(\bm{x},t) ,\;\;\text{on}\;\;\partial {\Omega_T}\times(0,t^*),\\
&{k_{ij}^{{\varepsilon}}\frac{\partial T_{\Delta}^{(1\varepsilon)}(\bm{x},t) }{\partial {x_j}}}{n_i}:= \bar \zeta_{1i}(\bm{x},t){n_i},\;\;\text{on}\;\;\partial {\Omega_q}\times(0,t^*),\\
&T_{\Delta}^{(1\varepsilon)}({\bm{x}},0)=-\varepsilon{M_{\alpha_1}}\frac{\partial \widetilde T}{\partial x_{\alpha_1}}:=\varepsilon \widetilde \omega_1(\bm{x}),\;\;\text{in}\;\;\Omega,
\end{aligned} \right.
\label{eq:residual_LOMS}
\end{equation}
where the specific expressions of functions ${F}_{0}(\bm{x},\bm{y},t)$, ${F}_{1}(\bm{x},\bm{y},t)$ are uncomplicated to obtain and be exhibited in Appendix A of the present study because of their lengthy expressions.

Then putting the residual function $T^{(2\varepsilon)}_{\Delta}$ into \eqref{eq:all}, we derive the following residual equation for the HOMS solution.
\begin{equation}
\left\{
\begin{aligned}
  &\rho^{\varepsilon}c^{\varepsilon}\frac{{\partial {T_{\Delta}^{(2\varepsilon)}}}}{{\partial t}}- \frac{\partial }{{\partial {x_i}}}\Big( {{k_{ij}^{\varepsilon}}\frac{{\partial {T_{\Delta}^{(2\varepsilon)}}}}{{\partial {x_j}}}}\Big)
    =\varepsilon F_2(\bm{x},\bm{y},t),\;\;\text{in}\;\;\Omega\times(0,t^*),\\
  &{T_{\Delta}^{(2\varepsilon)}(\bm{x},t) }
    =-\varepsilon{M_{\alpha_1}}\frac{\partial T_{0}}{\partial x_{\alpha_1}}
    - \varepsilon^2\Bigg[ G\frac{{\partial {T_{0}}}}{{\partial t}}+{M_{\alpha_1\alpha_2}}\frac{\partial^2 T_{0}}{\partial x_{\alpha_1}\partial x_{\alpha_2}}+N\big[T_0\big]^4 \\
  &\;\; +N_{\alpha_1}\frac{\partial T_{0}}{\partial x_{\alpha_1}}+N_{\alpha_1\alpha_2}\frac{\partial T_{0}}{\partial x_{\alpha_1}}\frac{\partial T_{0}}{\partial x_{\alpha_2}}\Bigg]
    :=\varepsilon \widehat \psi_2(\bm{x},t),\;\;\text{on}\;\;\partial {\Omega_T}\times(0,t^*),\\
  &{k_{ij}^{{\varepsilon}}\frac{\partial T_{\Delta}^{(2\varepsilon)}(\bm{x},t) }{\partial {x_j}}}{n_i}
    := \bar \zeta_{2i}(\bm{x},t){n_i},\;\;\text{on}\;\;\partial {\Omega_q}\times(0,t^*),\\
  &T_{\Delta}^{(2\varepsilon)}({\bm{x}},0)
    =-\varepsilon{M_{\alpha_1}}\frac{\partial \widetilde T}{\partial x_{\alpha_1}}
    - \varepsilon^2\Bigg[ G\frac{{\partial {T_{0}}}}{{\partial t}}\bigg|_{t=0}+{M_{\alpha_1\alpha_2}}\frac{\partial^2 \widetilde T}{\partial x_{\alpha_1}\partial x_{\alpha_2}}\\
  &\qquad +N\big[\widetilde T\big]^4 \!+N_{\alpha_1}\frac{\partial \widetilde T}{\partial x_{\alpha_1}}+N_{\alpha_1\alpha_2}\frac{\partial \widetilde T}{\partial x_{\alpha_1}}\frac{\partial \widetilde T}{\partial x_{\alpha_2}}\Bigg]
    :=\varepsilon \widetilde\omega_2(\bm{x}),\;\;\text{in}\;\;\Omega,
\end{aligned}
\right.
\label{eq:HOMS-res}
\end{equation}
where the detailed expressions of functions ${F}_{2}(\bm{x},\bm{y},t)$ are also uncomplicated to achieve and be displayed in Appendix A of the present study owing to their lengthy expressions.

From the preceding error results in the pointwise sense, it is apparent that LOMS solution fails to maintain local balance of heat quantity since the $\varepsilon$-independent terms ${F}_{0}(\bm{x},\bm{y},t)$ in \eqref{eq:residual_LOMS} cannot converge toward zero as the microstructural parameter $\varepsilon$ approaches zero. By virtue of the high-order correction terms, the HOMS solution can guarantee the local heat quantity balance of the thermal radiation equation in the original governing equation \eqref{eq:all} due to their $O(\varepsilon)$-order pointwise errors. This serves as the principal motivation for this study to establish the HOMS solutions exhibiting high-accuracy computation performance for composite structures.

\subsection{Error analysis of high-order multi-scale solution in the integral sense}
Next, we proceed with the global error analysis of HOMS solution in the integral sense. Firstly, the following assumptions are presented.
\begin{description}
    \item[($\mathrm{B}_{1}$)] Assume that $\Omega$ is a bounded, convex, and smooth domain with the boundary $\partial\Omega \in C^4$.
    \item[($\mathrm{B}_{2}$)] Assume that $\dfrac{\partial \rho^\varepsilon(\bm{x}, T^\varepsilon)}{\partial t}$, $\dfrac{\partial c^\varepsilon(\bm{x}, T^\varepsilon)}{\partial t}$ and $\dfrac{\partial k_{ij}^\varepsilon(\bm{x}, T^\varepsilon)}{\partial t} \in L^\infty(\Omega \times (0,t^*))$.
\end{description}
\begin{theorem}\label{thm:error_estimate}
    Let \( T^{\varepsilon} \) be the solution of multi-scale problem \eqref{eq:all}, and \( T_0 \) be the solution
    of the macroscopic homogenized problem~\eqref{eq:homogenized}. Suppose the assumptions \((A_1)\)-\((A_4)\) and \((B_1)\)-\((B_2)\) hold, and that \( T_0(\boldsymbol{x},t) \in L^\infty\left(0,t^*;H^{ 4}(\Omega)\right) \) and
    \( \displaystyle\frac{\partial T_0}{\partial t} \in L^\infty\left(0,t^{*};H^{ 2}(\Omega)\right) \) are satisfied.
    Then the following error estimate holds:
\begin{equation}
\begin{aligned}
\bigl\|T_{\Delta}^{ (2\varepsilon)}(\bm{x},t)\bigr\|_{ L^\infty\left(0,t^*;L^2(\Omega) \right)}
+ \bigl\|T_{ \Delta}^{ (2\varepsilon)}(\bm{x},t)\bigr\|_{ L^2 \left(0,t^*;H^1_0(\Omega) \right)}
\leq C(t^*) \varepsilon^{1/2},
\end{aligned}
\label{eq:error_estimate}
\end{equation}
where \( C(t^{*}) \) is a constant independent of the small periodic parameter \( \varepsilon \).
\end{theorem}
$\mathbf{Proof:}$ Firstly, let us introduce the cut-off function $m_{T,\varepsilon}(\bm{x})\in C^{\infty}(\bar{\Omega})$ in \cite{R4} for temperature field defined as follows
\begin{equation}
\left\{ {\begin{aligned}
&m_{T,\varepsilon}(\bm{x})=1,\;\mathrm{if}\;\mathrm{dist}(\bm{x},\partial\Omega_{T})\leq \varepsilon,\\
&m_{T,\varepsilon}(\bm{x})=0,\;\mathrm{if}\;\mathrm{dist}(\bm{x},\partial\Omega_{T})\geq 2\varepsilon,\\
&\|\bigtriangledown m_{T,\varepsilon}(\bm{x})\|_{L^\infty(\Omega)}\leq C\varepsilon^{-1},\;\mathrm{if}\; \varepsilon<\mathrm{dist}(\bm{x},\partial\Omega_{T})<2\varepsilon.
\end{aligned}} \right.
\end{equation}
Then, define novel residual function for temperature field as below
\begin{equation}
\wideparen T^{(2\varepsilon)}_{\Delta}(\bm{x},t)=T^{(2\varepsilon)}_{\Delta}(\bm{x},t)-\varepsilon m_{T,\varepsilon}(\bm{x})\widehat\psi_2(\bm{x},t).
\end{equation}
Next, multiplying the residual equation in \eqref{eq:HOMS-res} by $\wideparen T^{(2\varepsilon)}_{\Delta}$ and integrating over $\Omega$, the following equality is formulated by further employing Green's formula.
\begin{equation}
\begin{aligned}
& \int_{\Omega} \rho^{\varepsilon}(\bm{x},T^{\varepsilon}) c^{\varepsilon}(\bm{x},T^{\varepsilon}) \frac{\partial \wideparen T^{(2\varepsilon)}_{\Delta}}{\partial t}\wideparen T^{(2\varepsilon)}_{\Delta} \, \mathrm{d}\Omega
+ \int_{\Omega} k_{ij}^{\varepsilon}(\bm{x},T^{\varepsilon}) \frac{\partial \wideparen T^{(2\varepsilon)}_{\Delta}}{\partial x_j} \frac{\partial \wideparen T^{(2\varepsilon)}_{\Delta}}{\partial x_i} \, \mathrm{d}\Omega \\
&= \int_{\Omega} \varepsilon F_2(\bm{x},\bm{y},t) \wideparen T^{(2\varepsilon)}_{\Delta} \, \mathrm{d}\Omega+ \int_{\partial {\Omega_q}} \bar \zeta_{2i}(\bm{x},t){n_i}\wideparen T^{(2\varepsilon)}_{\Delta} \mathrm{d}s\\
&-\int_{\Omega}\rho^{\varepsilon}({\bm{x}},{T^\varepsilon })c^{\varepsilon}({\bm{x}},{T^\varepsilon })\frac{{\partial \big[{\varepsilon m_{T,\varepsilon}(\bm{x})\widehat\psi_2(\bm{x},t)}}\big]}{{\partial t}}\wideparen T_\Delta^{(2\varepsilon)}\mathrm{d}\Omega\\
&-\int_{\Omega} {k_{ij}^{\varepsilon}({\bm{x}},{T^\varepsilon })\frac{{\partial \big[\varepsilon m_{T,\varepsilon}(\bm{x})\widehat\psi_2(\bm{x},t)}\big]}{\partial {x_j}}}\frac{{\partial {\wideparen T_\Delta^{(2\varepsilon)} }}}{\partial {x_i}}\mathrm{d}\Omega.
\end{aligned}
\label{eq:simplified_form}
\end{equation}
From \eqref{eq:simplified_form}, it is easily verified that the following equations hold:
\begin{equation}
\begin{aligned}
& \frac{1}{2} \frac{\partial}{\partial t} \int_{\Omega} \rho^{\varepsilon} c^{\varepsilon} \left(\wideparen T^{(2\varepsilon)}_{\Delta}\right)^2 \, \mathrm{d}\Omega
+ \int_{\Omega} k_{ij}^{\varepsilon} \frac{\partial \wideparen T^{(2\varepsilon)}_{\Delta}}{\partial x_j} \frac{\partial \wideparen T^{(2\varepsilon)}_{\Delta}}{\partial x_i} \, \mathrm{d}\Omega \\
&= \int_{\Omega} \varepsilon F_2 \wideparen T^{(2\varepsilon)}_{\Delta} \, \mathrm{d}\Omega+ \int_{\partial {\Omega_q}} \bar \zeta_{2i}{n_i}\wideparen T^{(2\varepsilon)}_{\Delta} \mathrm{d}s\\
&+ \frac{1}{2}\int_{\Omega} \frac{\partial \rho^{\varepsilon}}{\partial t} c^{\varepsilon} \left(\wideparen T^{(2\varepsilon)}_{\Delta}\right)^2 \, \mathrm{d}\Omega
+ \frac{1}{2}\int_{\Omega} \rho^{\varepsilon} \frac{\partial c^{\varepsilon}}{\partial t} \left(\wideparen T^{(2\varepsilon)}_{\Delta}\right)^2 \, \mathrm{d}\Omega\\
&-\int_{\Omega}\rho^{\varepsilon}c^{\varepsilon}\frac{{\partial \big[{\varepsilon m_{T,\varepsilon}\widehat\psi_2}}\big]}{{\partial t}}\wideparen T_\Delta^{(2\varepsilon)}\mathrm{d}\Omega-\int_{\Omega} {k_{ij}^{\varepsilon}\frac{{\partial \big[\varepsilon m_{T,\varepsilon}\widehat\psi_2}\big]}{\partial {x_j}}}\frac{{\partial {\wideparen T_\Delta^{(2\varepsilon)} }}}{\partial {x_i}}\mathrm{d}\Omega.
\end{aligned}
\label{eq:energy_1}
\end{equation}
Then, integrating both sides over the time interval \([0,t]\) (\(0<t\leq t^*\)), and simplifying the calculation, we obtain:
\begin{equation}
\begin{aligned}
& \int_{\Omega} \rho^{\varepsilon} c^{\varepsilon} \left(\wideparen T^{(2\varepsilon)}_{\Delta}\right)^2 \, \mathrm{d}\Omega
+ \int_{0}^{t} \int_{\Omega} 2k_{ij}^{\varepsilon} \frac{\partial \wideparen T^{(2\varepsilon)}_{\Delta}}{\partial x_j} \frac{\partial \wideparen T^{(2\varepsilon)}_{\Delta}}{\partial x_i} \, \mathrm{d}\Omega\mathrm{d}\tau \\
&=  \int_{0}^{t} \int_{\Omega} 2\varepsilon F_2 \wideparen T^{(2\varepsilon)}_{\Delta} \, \mathrm{d}\Omega\mathrm{d}\tau
+ \int_{\Omega} \rho^{\varepsilon} c^{\varepsilon} \left(\varepsilon \widetilde\omega_2\right)^2 \, \mathrm{d}\Omega+ \int_{0}^{t}\int_{\partial {\Omega_q}}2 \bar \zeta_{2i}{n_i}\wideparen T^{(2\varepsilon)}_{\Delta} \mathrm{d}s\mathrm{d}t\\
& + \int_{0}^{t} \int_{\Omega} \frac{\partial \rho^{\varepsilon}}{\partial \tau} c^{\varepsilon} \left(\wideparen T^{(2\varepsilon)}_{\Delta}\right)^2 \, \mathrm{d}\Omega\mathrm{d}\tau
+ \int_{0}^{t} \int_{\Omega} \rho^{\varepsilon} \frac{\partial c^{\varepsilon}}{\partial \tau} \left(\wideparen T^{(2\varepsilon)}_{\Delta}\right)^2 \, \mathrm{d}\Omega\mathrm{d}\tau\\
&- \int_{0}^{t} \int_{\Omega}2\rho^{\varepsilon}c^{\varepsilon}\frac{{\partial \big[{\varepsilon m_{T,\varepsilon}\widehat\psi_2}}\big]}{{\partial \tau}}\wideparen T_\Delta^{(2\varepsilon)}\mathrm{d}\Omega\mathrm{d}\tau- \int_{0}^{t} \int_{\Omega}2 {k_{ij}^{\varepsilon}\frac{{\partial \big[\varepsilon m_{T,\varepsilon}\widehat\psi_2}\big]}{\partial {x_j}}}\frac{{\partial {\wideparen T_\Delta^{(2\varepsilon)} }}}{\partial {x_i}}\mathrm{d}\Omega\mathrm{d}\tau.
\end{aligned}
\label{eq:with_initial}
\end{equation}

In the following, we estimate both sides of \eqref{eq:with_initial} separately. First, for the left-hand side, using Assumption $(A_1)$-$(A_3)$ and Poincar\'{e}-Friedrichs inequality, we obtain the following inequality:
\begin{equation}
\begin{aligned}
& \int_{\Omega} \rho^{\varepsilon} c^{\varepsilon} \left(\wideparen T^{(2\varepsilon)}_{\Delta}(\bm{x},t)\right)^2 \, \mathrm{d}\Omega
+ \int_{0}^{t} \int_{\Omega} 2k_{ij}^{\varepsilon} \frac{\partial \wideparen T^{(2\varepsilon)}_{\Delta}(\bm{x},\tau)}{\partial x_j} \frac{\partial \wideparen T^{(2\varepsilon)}_{\Delta}(\bm{x},\tau)}{\partial x_i} \, \mathrm{d}\Omega \, \mathrm{d}\tau \\
& \ge \underline{\rho}\underline{c}\left\| \wideparen T^{(2\varepsilon)}_{\Delta}(\bm{x},t) \right\|_{L^2(\Omega)}^2
+ C_1 \int_{0}^{t} \left\| \wideparen T^{(2\varepsilon)}_{\Delta}(\bm{x},\tau) \right\|_{H^1_0(\Omega)}^2 \, \mathrm{d}\tau \\
& \ge \underline{C}\left( \left\| \wideparen T^{(2\varepsilon)}_{\Delta}(\bm{x},t) \right\|_{L^2(\Omega)}^2
+ \int_{0}^{t} \left\| \wideparen T^{(2\varepsilon)}_{\Delta}(\bm{x},\tau) \right\|_{H^1_0(\Omega)}^2 \, \mathrm{d}\tau \right),
\end{aligned}
\label{eq:left_hand_inequality}
\end{equation}
where $\underline{C}= {\min}\left(\underline{\rho}\underline{c}, C_{1} \right)$.

Then, using Assumption $(A_1)$-$(A_3)$, Schwarz's inequality, Young's inequality $\displaystyle ab\leq\frac{1}{2}(\lambda a^2+\frac{1}{\lambda}b^2)$ for all parameter $\lambda\in\mathbb{R}^{+}$, Lemma 2.2 in Chapter 2 of reference \cite{R9}, and selecting the same parameter $\lambda$ for any arbitrary $\tau\in[0,t]$, the estimation of the right-hand side of \eqref{eq:with_initial} yields the following inequality:
\begin{equation}
\begin{aligned}
& \int_{0}^{t} \int_{\Omega} 2\varepsilon F_2 \wideparen T^{(2\varepsilon)}_{\Delta} \, \mathrm{d}\Omega \, \mathrm{d}\tau
+ \int_{\Omega} \rho^{\varepsilon} c^{\varepsilon} \left(\varepsilon \widetilde \omega_2\right)^2 \, \mathrm{d}\Omega+ \int_{0}^{t}\int_{\partial {\Omega_q}}2 \bar \zeta_{2i}{n_i}\wideparen T^{(2\varepsilon)}_{\Delta} \mathrm{d}s\mathrm{d}t\\
& + \int_{0}^{t} \int_{\Omega} \frac{\partial \rho^{\varepsilon}}{\partial \tau} c^{\varepsilon} \left(\wideparen T^{(2\varepsilon)}_{\Delta}\right)^2 \, \mathrm{d}\Omega\mathrm{d}\tau
+ \int_{0}^{t} \int_{\Omega} \rho^{\varepsilon} \frac{\partial c^{\varepsilon}}{\partial \tau} \left(\wideparen T^{(2\varepsilon)}_{\Delta}\right)^2 \, \mathrm{d}\Omega\mathrm{d}\tau\\
&- \int_{0}^{t} \int_{\Omega}2\rho^{\varepsilon}c^{\varepsilon}\frac{{\partial \big[{\varepsilon m_{T,\varepsilon}\widehat\psi_2}}\big]}{{\partial \tau}}\wideparen T_\Delta^{(2\varepsilon)}\mathrm{d}\Omega\mathrm{d}\tau- \int_{0}^{t} \int_{\Omega}2 {k_{ij}^{\varepsilon}\frac{{\partial \big[\varepsilon m_{T,\varepsilon}\widehat\psi_2}\big]}{\partial {x_j}}}\frac{{\partial {\wideparen T_\Delta^{(2\varepsilon)} }}}{\partial {x_i}}\mathrm{d}\Omega\mathrm{d}\tau\\
&\leq 2\int_{0}^{t} \int_{\Omega} \left( \frac{(\varepsilon F_2)^2}{2}
+ \frac{\big(\wideparen T^{(2\varepsilon)}_{\Delta}\big)^2}{2}\right) \, \mathrm{d}\Omega\mathrm{d}\tau
+ C_2\varepsilon^2\\
&+ C_3\int_{0}^{t}\varepsilon^{1/2}\left\|\wideparen T^{(2\varepsilon)}_{\Delta}\right\|_{H^1_0(\Omega)}^2 \mathrm{d}\tau
+ C_4\int_{0}^{t}\left\|\wideparen T^{(2\varepsilon)}_{\Delta}\right\|_{L^2(\Omega)}^2 \mathrm{d}\tau\\
&+C_5\int_{0}^{t} \int_{\Omega}\varepsilon\wideparen T_\Delta^{(2\varepsilon)}\mathrm{d}\Omega\mathrm{d}\tau+C_6 \int_{0}^{t} \int_{\Omega}\varepsilon\frac{{\partial {\wideparen T_\Delta^{(2\varepsilon)} }}}{\partial {x_i}}\mathrm{d}\Omega\mathrm{d}\tau\\
&\leq 2\int_{0}^{t} \int_{\Omega} \left( \frac{(\varepsilon F_2)^2}{2}
+ \frac{\big(\wideparen T^{(2\varepsilon)}_{\Delta}\big)^2}{2}\right) \, \mathrm{d}\Omega \, \mathrm{d}\tau
+ C_2\varepsilon^2\\
&+ C_3\int_0^t\Big[\frac{1}{2\lambda}\varepsilon+\frac{\lambda}{2}\left\|\wideparen T^{(2\varepsilon)}_{\Delta}\right\|_{H_0^1(\Omega)}^2\Big]\mathrm{d}\tau
+ C_4\int_{0}^{t}\left\|\wideparen T^{(2\varepsilon)}_{\Delta}\right\|_{L^2(\Omega)}^2 \mathrm{d}\tau\\
&+ C_5\int_0^t\Big[\frac{1}{2}\varepsilon^2+\frac{1}{2}\left\|\wideparen T^{(2\varepsilon)}_{\Delta}\right\|_{L^2(\Omega)}^2\Big]\mathrm{d}\tau+ C_6\int_0^t\Big[\frac{1}{2\lambda}\varepsilon^2+\frac{\lambda}{2}\left\|\wideparen T^{(2\varepsilon)}_{\Delta}\right\|_{H_0^1(\Omega)}^2\Big]\mathrm{d}\tau\\
&\leq \overline{C}\varepsilon + \overline{C}\int_{0}^{t} \left\|\wideparen T^{(2\varepsilon)}_{\Delta} \right\|_{L^2(\Omega)}^2\mathrm{d}\tau
+ \frac{\overline{C}\lambda}{2}\int_0^t\left\|\wideparen T^{(2\varepsilon)}_{\Delta}\right\|_{H_0^1(\Omega)}^2\mathrm{d}\tau.
\end{aligned}
\label{eq:right_hand_inequality}
\end{equation}

By selecting a sufficiently small $\lambda$, and combining \eqref{eq:left_hand_inequality} and \eqref{eq:right_hand_inequality}, we directly obtain the following inequality:
\begin{equation}
\begin{aligned}
&\underline{C}\bigg( \big\|\wideparen T_{\Delta}^{(2\varepsilon)}(\bm{x},t) \big\|_{L^2(\Omega)}^2
+ \int_0^t \big\|\wideparen T_{\Delta}^{(2\varepsilon)}(\bm{x},\tau) \big\|_{H_0^1(\Omega)}^2 \, \mathrm{d}\tau \bigg) \\
&\leq \overline{C}\varepsilon+ \overline{C}\int_0^t \big\|\wideparen T_{\Delta}^{(2\varepsilon)}(\bm{x},\tau) \big\|_{L^2(\Omega)}^2 \, \mathrm{d}\tau
+ \overline{C}\int_0^t \int_0^{\tau} \big\|\wideparen T_{\Delta}^{(2\varepsilon)}(\bm{x},s) \big\|_{H_0^1(\Omega)}^2 \mathrm{d}s \, \mathrm{d}\tau.
\end{aligned}
\label{eq:final_inequality}
\end{equation}
Redefining $C=\overline{C}/\underline{C}$, and $\displaystyle\Theta(t) = \left\|\wideparen T^{(2\varepsilon)}_{\Delta} \right\|_{L^2(\Omega)}^2+ \int_{0}^{t} \left\|\wideparen T^{(2\varepsilon)}_{\Delta} \right\|_{H^1_0(\Omega)}^2 \, \mathrm{d}\tau$, the above inequality \eqref{eq:final_inequality} can be rewritten as:
\begin{equation}
\Theta(t) \leq C\varepsilon + C\int_{0}^{t} \Theta(\tau) \, \mathrm{d}\tau.
\label{eq:simplified_inequality}
\end{equation}
Subsequently, applying Gronwall's inequality to solve the above integral inequality, we obtain:
\begin{equation}
\left\|\wideparen T^{(2\varepsilon)}_{\Delta}(\bm{x},t) \right\|_{L^2(\Omega)}^2
+ \int_{0}^{t} \left\|\wideparen T^{(2\varepsilon)}_{\Delta}(\bm{x},\tau) \right\|_{H^1_0(\Omega)}^2 \, \mathrm{d}\tau
\leq C(t^*) \varepsilon.
\label{eq:after_gronwall}
\end{equation}
Furthermore, in virtue of the triangle inequality and by adopting the error estimate $\left\|m_{T,\varepsilon}(\bm{x})\widehat\psi_2(\bm{x},t)\right\|_{H^1(\Omega)} =\left\|m_{T,\varepsilon}(\bm{x})\widehat\psi_2(\bm{x},t)\right\|_{H^1(K_{\varepsilon})}
\leq C\varepsilon^{-{1}/{2}}$ in chapter 7 of reference \cite{R4} where $K_{\varepsilon}=\left\{{\bm{x}|\mathrm{dist}(\bm{x},\partial\Omega_{T})\leq 2\varepsilon}\right\}\cap\Omega$, one can readily obtain
\begin{equation}
\begin{aligned}
&\big \|T_{\Delta}^{(2\varepsilon)}\big\|_{L^2(\Omega)}^2
+\mathlarger{\int}_0^t\big \|T_{\Delta}^{(2\varepsilon)}\big\|_{H^1(\Omega)}^2\mathrm{d}\tau\\
&\leq C(t^*) \varepsilon+\big \|\varepsilon m_{T,\varepsilon}(\bm{x})\widehat\psi_2(\bm{x},t)\big\|_{L^2(\Omega)}^2
+\mathlarger{\int}_0^t\big \|\varepsilon m_{T,\varepsilon}(\bm{x})\widehat\psi_2(\bm{x},t)\big\|_{H^1(\Omega)}^2d\tau\\
&\leq C(t^*) \varepsilon+(\varepsilon C\varepsilon^{-1/2})^2+(\varepsilon C\varepsilon^{-1/2})^2\leq C(t^*) \varepsilon.
\end{aligned}
\label{eq:after_new}
\end{equation}
Then, applying the mean value inequality to \eqref{eq:after_new}, we obtain the following inequality:
\begin{equation}
\left\| T^{(2\varepsilon)}_{\Delta} \right\|_{L^2(\Omega)}
+ \left\| T^{(2\varepsilon)}_{\Delta} \right\|_{L^2(0,t;H^1_0(\Omega))}
\leq C(t^*) \varepsilon^{1/2}.
\label{eq:final_norm_estimate}
\end{equation}
Finally, utilizing the arbitrariness of the time variable \(t\), the following inequality holds:
\begin{equation}
\left\| T^{(2\varepsilon)}_{\Delta} \right\|_{L^\infty(0,t^*;L^2(\Omega))}
+ \left\| T^{(2\varepsilon)}_{\Delta} \right\|_{L^2(0,t^*;H^1_0(\Omega))}
\leq C(t^*) \varepsilon^{1/2}.
\label{eq:norm_estimate}\end{equation}
In summary, Theorem~\ref{thm:error_estimate} is proved.

Since the HOMS solution of multi-scale nonlinear problem (1) does not conform to the boundary conditions when $\Omega$ is a general domain, obtaining the optimal convergence order is impeded by the resulting boundary error. To gain the optimal
error estimation, certain hypotheses are further presented as follows:
\begin{description}
    \item[(I)] Assume that $\Omega$ can be exactly covered by a union of periodic unit cells. Let $\mathcal{T}_\varepsilon$ be
      the index set defined by $\mathcal{T}_\varepsilon = \bigl\{\mathbf{z} = (z_1,\cdots,z_d)
      \in \mathbb{Z}^d : \varepsilon(\mathbf{z} + \overline{\Omega}_{\bm{y}})
      \subset \overline{\Omega}\bigr\},$
      and $\overline{\Omega} = \bigcup_{\mathbf{z}\in\mathcal{T}_\varepsilon}
      \varepsilon(\mathbf{z} + \overline{\Omega}_{\bm{y}})$. Furthermore, let
      $E_{\mathbf{z}} = \varepsilon(\mathbf{z} + \Omega_{\bm{y}})$ be the
      translational unit cell and $\partial E_{\mathbf{z}}$ be its boundary.
    \item[(II)] ${\rho^\varepsilon(\bm{x}, T^\varepsilon)}$, $c^\varepsilon(\bm{x}, T^\varepsilon)$ and
      $ k_{ii}^\varepsilon(\bm{x}, T^\varepsilon)$ are symmetric with respect to
      all mid-planes of the reference unit cell, while $ k_{ij}^\varepsilon(\bm{x}, T^\varepsilon)$ is skew-symmetric with respect to all mid-planes of the reference unit cell $\Omega_{\bm{y}}$ \cite{R26,R35,R36,R39,R40,R41}.
    \item[(III)] Apply the homogeneous Dirichlet boundary condition to replace the periodic
boundary condition for whole auxiliary cell functions \cite{R39,R40,R41}.
    \item[(IV)] The multi-scale nonlinear problem (1) is imposed with pure Dirichlet boundary condition.
\end{description}
\begin{lemma}\label{lemma:cell_functions_continuity}
  If the multi-scale radiative heat transfer problem~\eqref{eq:all} satisfies assumptions $(A_1)$-$(A_3)$ and $(II)$-$(III)$, then the normal derivatives of all cell functions $\sigma_{TY}(M_{\alpha_1})$,
$\sigma_{TY}(M_{\alpha_1\alpha_2})$, $\sigma_{TY}(G)$, $\sigma_{TY}(N)$, $\sigma_{TY}(N_{\alpha_1})$, $\sigma_{TY}(N_{\alpha_1\alpha_2})$ are continuous on the boundary $\partial \Omega_{\bm{y}}$ of the microscopic unit cell domain $\Omega_{\bm{y}}$ by employing the proof approach in \cite{R39,R40,R41}, where the differential operator $\displaystyle\sigma_{TY} (\Phi)= n_i k_{ij}(\bm{y},{T^\varepsilon})\frac{\partial \Phi}{\partial y_j}$.
\end{lemma}
\begin{corollary}
Assume that $\Omega$ is the integral periodic region. Let $T^\varepsilon(\bm{x},t)$ be the weak solution of multi-scale nonlinear equation (1), $T_0(\bm{x},t)$ be the weak solution of corresponding homogenized equation (14), ${T^{(2\varepsilon)} }(\bm{x},t)$ be the HOMS solution given by formula (23). Under the above hypotheses (A$_1$)-(A$_4$), (B$_1$)-(B$_2$) and (I)-(IV), the following global error estimation is obtained.
\begin{equation}
\begin{aligned}
\bigl\|T_{\Delta}^{ (2\varepsilon)}(\bm{x},t)\bigr\|_{ L^\infty\left(0,t^*;L^2(\Omega) \right)}
+ \bigl\|T_{ \Delta}^{ (2\varepsilon)}(\bm{x},t)\bigr\|_{ L^2 \left(0,t^*;H^1_0(\Omega) \right)}
\leq C(t^*) \varepsilon,
\end{aligned}
\label{eq:error_estimate_integral_periodic}
\end{equation}
\end{corollary}
where \( C(t^{*}) \) is a constant independent of the small periodic parameter \( \varepsilon \).
$\mathbf{Proof:}$ Recalling the above proof again, based on assumptions (I)-(IV), the error order $\varepsilon^{{1}/{2}}$ generating from boundary $\partial \Omega$ will disappear in the proof. However, it should be noted that the new auxiliary cell functions equipped with Dirichlet boundary condition lack enough regularity on the outer boundary of microscopic unit cell $\Omega_{\bm{y}}$ in general case.

Firstly, on the basis of the multi-scale chain rule~\eqref{eq:chain}, we can derive the following equality for HOMS solution defined in \eqref{eq:HOMS}.
\begin{equation}
\begin{aligned}
&\sigma_T\left(T^{(2\varepsilon)}\right)
= n_i k_{ij}^{\varepsilon}(\bm{x},T^{\varepsilon}) \frac{\partial T^{(2\varepsilon)}}{\partial x_j} \\
&= n_ik_{ij}(\bm{y},{T^\varepsilon})\left(\frac{\partial}{\partial x_j} \!+\! \frac{1}{\varepsilon} \frac{\partial}{\partial y_j}\right)\biggl[T_0 + \varepsilon M_{\alpha_1} \frac{\partial T_0}{\partial x_{\alpha_1}}
+ \varepsilon^2 \Big(M_{\alpha_1\alpha_2} \frac{\partial^2 T_0}{\partial x_{\alpha_1} \partial x_{\alpha_2}} \\
&\;\;\; + G \frac{\partial T_0}{\partial t}+ N [T_0]^4+ N_{\alpha_1} \frac{\partial T_0}{\partial x_{\alpha_1}}+ N_{\alpha_1\alpha_2} \frac{\partial T_0}{\partial x_{\alpha_1}} \frac{\partial T_0}{\partial x_{\alpha_2}}\Big)\biggr]\\
&= n_i k_{ij}(\bm{y},{T^\varepsilon})\frac{\partial}{\partial x_j}\biggl[T_0 + \varepsilon M_{\alpha_1} \frac{\partial T_0}{\partial x_{\alpha_1}}+ \varepsilon^2 \Big(M_{\alpha_1\alpha_2} \frac{\partial^2 T_0}{\partial x_{\alpha_1} \partial x_{\alpha_2}} \\
&\;\;\;+ G \frac{\partial T_0}{\partial t}+ N [T_0]^4+ N_{\alpha_1} \frac{\partial T_0}{\partial x_{\alpha_1}}
+ N_{\alpha_1\alpha_2} \frac{\partial T_0}{\partial x_{\alpha_1}} \frac{\partial T_0}{\partial x_{\alpha_2}}
\Big)\biggr] \\
&\;\;\;+ n_i k_{ij}(\bm{y},{T^\varepsilon})\biggl[\frac{\partial M_{\alpha_1}}{\partial y_j} \frac{\partial T_0}{\partial x_{\alpha_1}}+ \varepsilon \Bigl(\frac{\partial M_{\alpha_1\alpha_2}}{\partial y_j} \frac{\partial^2 T_0}{\partial x_{\alpha_1} \partial x_{\alpha_2}} \\
&\;\;\; + \frac{\partial G}{\partial y_j} \frac{\partial T_0}{\partial t}+ \frac{\partial N}{\partial y_j} [T_0]^4
+ \frac{\partial N_{\alpha_1}}{\partial y_j} \frac{\partial T_0}{\partial x_{\alpha_1}}
+ \frac{\partial N_{\alpha_1\alpha_2}}{\partial y_j} \frac{\partial T_0}{\partial x_{\alpha_1}} \frac{\partial T_0}{\partial x_{\alpha_2}}\Bigr)\biggr]\\
&= n_i k_{ij}(\bm{y},{T^\varepsilon})\frac{\partial}{\partial x_j}\biggl[
T_0 + \varepsilon M_{\alpha_1} \frac{\partial T_0}{\partial x_{\alpha_1}}
+ \varepsilon^2 \Big(M_{\alpha_1\alpha_2} \frac{\partial^2 T_0}{\partial x_{\alpha_1} \partial x_{\alpha_2}} \\
&\;\;\; + G \frac{\partial T_0}{\partial t}+ N [T_0]^4+ N_{\alpha_1} \frac{\partial T_0}{\partial x_{\alpha_1}}
+ N_{\alpha_1\alpha_2} \frac{\partial T_0}{\partial x_{\alpha_1}} \frac{\partial T_0}{\partial x_{\alpha_2}}
\Big)\biggr]\\
&\;\;\; + \sigma_{TY}(M_{\alpha_1}) \frac{\partial T_0}{\partial x_{\alpha_1}}
+ \varepsilon \sigma_{TY}(M_{\alpha_1\alpha_2}) \frac{\partial^2 T_0}{\partial x_{\alpha_1} \partial x_{\alpha_2}}+ \varepsilon \sigma_{TY}(G) \frac{\partial T_0}{\partial t} \\
&\;\;\; + \varepsilon \sigma_{TY}(N) [T_0]^4 + \varepsilon \sigma_{TY}(N_{\alpha_1}) \frac{\partial T_0}{\partial x_{\alpha_1}}
+ \varepsilon \sigma_{TY}(N_{\alpha_1\alpha_2}) \frac{\partial T_0}{\partial x_{\alpha_1}} \frac{\partial T_0}{\partial x_{\alpha_2}}.
\end{aligned}
\label{eq:integral_form}
\end{equation}
Then, combining \eqref{eq:integral_form} and Lemma 1, we calculate to obtain:
\begin{equation}
\begin{aligned}
&\quad\sum_{\mathbf{z}\in \mathcal{T}_\varepsilon} \int_{\partial E_\mathbf{z}} \sigma_T\left(T^{(2\varepsilon)}_{\Delta}\right) \cdot T^{(2\varepsilon)}_{\Delta} \, \mathrm{d}\Gamma_{\bm{y}} \\
&= \sum_{\mathbf{z}\in \mathcal{T}_\varepsilon} \int_{\partial E_\mathbf{z}} \sigma_T\left(T^\varepsilon - T^{(2\varepsilon)}\right) \cdot T^{(2\varepsilon)}_{\Delta} \, \mathrm{d}\Gamma_{\bm{y}}\\
&= -\sum_{\mathbf{z}\in \mathcal{T}_\varepsilon} \int_{\partial E_\mathbf{z}} \sigma_T\left(T^{(2\varepsilon)}\right) \cdot T^{(2\varepsilon)}_{\Delta} \, \mathrm{d}\Gamma_{\bm{y}}= 0.
\end{aligned}
\label{eq:I_results}
\end{equation}
Finally, following along the lines of the proof of Theorem 1, the proof of Corollary 1 is completed.

\section{Two-stage multi-scale numerical algorithm}
\begin{figure}[!htb]
\centering
\includegraphics[width=\linewidth]{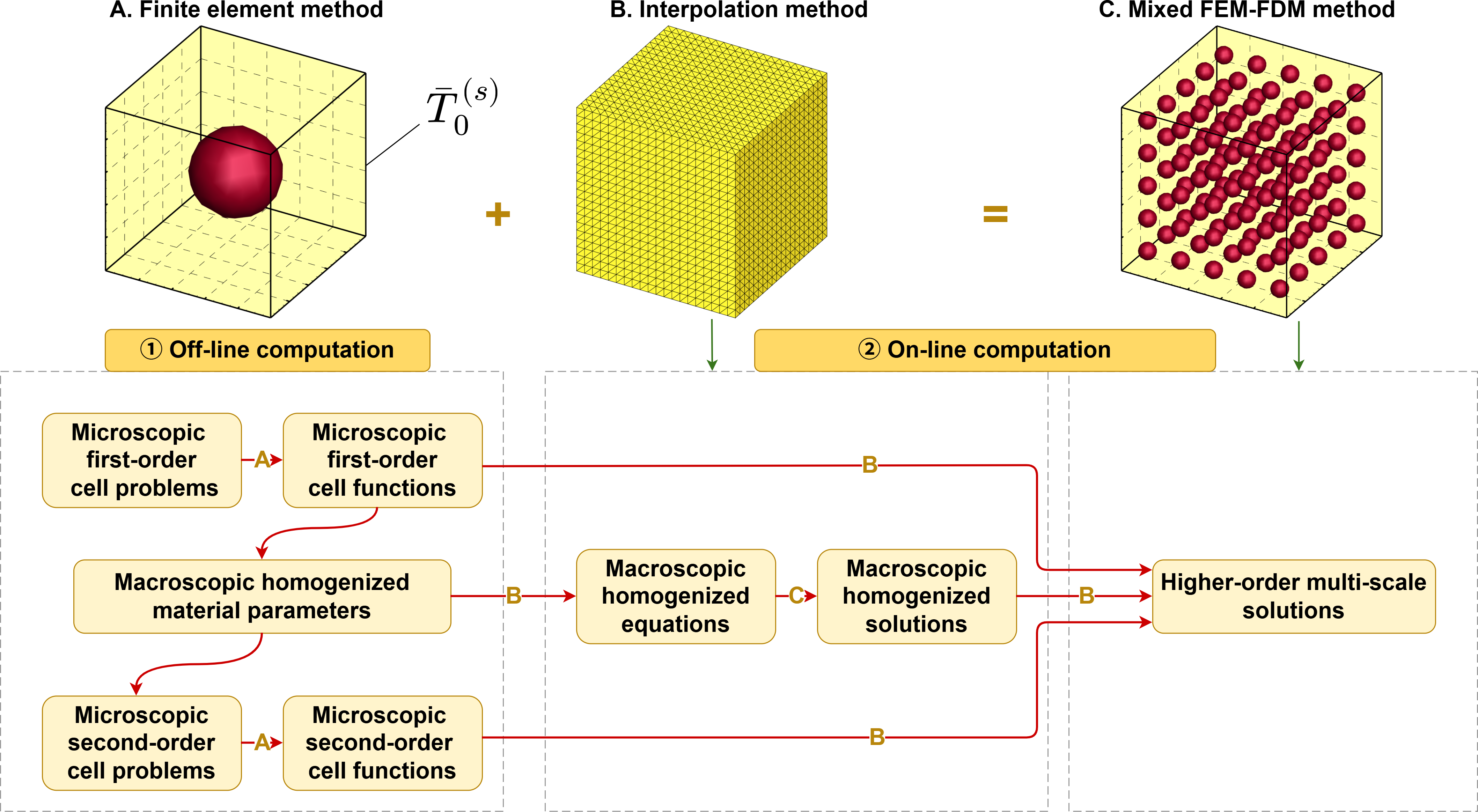}
\caption{The schematic diagram of two-stage multi-scale numerical algorithm.}
\label{fig:algorithm}
\end{figure}

The proposed multi-scale computational framework consists of microscopic cell problems, macroscopic
homogenized problem and high-order multi-scale solutions, which comprise a closed solving system.
Noting that all microscopic cell functions defined in Section~\ref{HOMS} are dependent on the macroscopic
temperature $T_0$, the continuous property of microscopic cell functions can be proved similarly as in
\cite{R23,R27}. According to this continuous property, microscopic cell functions only need to be evaluated
corresponding to a few representative macroscopic temperatures rather than all temperature points, and
then the interpolation technique can be exploited to obtain the auxiliary cell functions during the
simulation process \cite{R22,R23,R27,R28}. In the following, we present a two-stage numerical algorithm
comprising off-line and on-line stages for efficiently simulating the time-dependent nonlinear thermal
radiation problem \eqref{eq:all} of composite structures, as elaborated in Fig.~\ref{fig:algorithm}.

\subsection{Off-line stage: computation for microscopic cell problems}
\label{sec:offline}
\begin{enumerate}
\item[(1)]
Determine the geometric configuration of the periodic unit cell $\Omega_{\bm{y}}=(0,1)^d$ and generate a tetrahedral (or triangular) finite element mesh $\mathcal{T}_{h_1}=\{K\}$ for the unit cell $\Omega_{\bm{y}}$, where $h_1=\max_K\{h_K\}$. Then define the linear conforming finite element space $V_{h_1}(\Omega_{\bm{y}})=\{\nu\in C^0(\overline{\Omega}_{\bm{y}}):\nu|_{\partial \Omega_{\bm{y}}}=0,\ \nu|_K\in P_1(K)\}\subset H_0^1(\Omega_{\bm{y}})$ for the auxiliary cell problems.
\item[(2)]
Define the computational temperature range $[T_0^{\min},T_0^{\max}]$ and select a certain number of representative macroscopic temperatures $\{\bar T_0^{(s)}\}_{s=0}^{N_T}$ within this range. Note that the classical periodic boundary condition of auxiliary cell problems is replaced by the homogeneous Dirichlet boundary condition for practical numerical implementation \cite{R39,R40,R41}. Next, employ the finite element method to solve the first-order auxiliary cell functions \eqref{eq:Malpha1} on $V_{h_1}(\Omega_{\bm{y}})$ corresponding to each distinct representative macroscopic temperature $\bar T_0^{(s)}$. The specific finite element scheme for the first-order cell problem \eqref{eq:Malpha1} is established as follows:
\begin{equation}
\begin{aligned}
& \int_{\Omega_{\bm{y}}} k_{ij}^{(0)}(\bm{y},\bar T_0^{(s)}) \frac{\partial M_{\alpha_1}(\bm{y},\bar T_0^{(s)})}{\partial y_j} \frac{\partial \nu^{h_1}}{\partial y_i} \, \mathrm{d}\Omega_{\bm{y}} \\
& \qquad = - \int_{\Omega_{\bm{y}}} k_{i\alpha_1}^{(0)}(\bm{y},\bar T_0^{(s)}) \frac{\partial \nu^{h_1}}{\partial y_i} \, \mathrm{d}\Omega_{\bm{y}}, \quad \forall \nu^{h_1} \in V_{h_1}(\Omega_{\bm{y}}).
\label{eq:FEM_Malpha1_offline}
\end{aligned}
\end{equation}
\item[(3)]
The macroscopic homogenized material parameters $\bar S(T_0)$, $\bar k_{ij}(T_0)$ and $\bar\beta(T_0)$ are
evaluated by formula \eqref{eq:homogenized_params} corresponding to each representative macroscopic temperature $\bar T_0^{(s)}$.
\item[(4)]
Using the same mesh as for the first-order cell functions, the second-order auxiliary cell functions \eqref{eq:Malpha1alpha2}, \eqref{eq:G}, \eqref{eq:N}, \eqref{eq:Nalpha1} and \eqref{eq:Nalpha1alpha2}, which correspond to each distinct representative macroscopic temperature $\bar T_0^{(s)}$, are computed on $V_{h_1}(\Omega_{\bm{y}})$ by employing the finite element method respectively.
\end{enumerate}
\subsection{On-line stage: computation for macroscopic homogenized problem}
\label{sec:online}
\begin{enumerate}
\item[(1)]
Let $\mathcal{T}_{h_0}=\{e\}$ be a finite element mesh of the macroscopic domain $\Omega$, where $h_0=\max_e\{h_e\}$. Then define the linear conforming finite element space $V_{h_0}(\Omega)=\{\nu\in C^0(\bar{\Omega}):\nu|_{\partial\Omega_T}=0,\ \nu|_e\in P_1(e)\}\subset H^1(\Omega)$. The homogenized material parameters at each node $\bm{x}$ of $V_{h_0}(\Omega)$ can be obtained by the interpolation approach using the precomputed values at representative temperatures.
\item[(2)]
Solve the macroscopic homogenized problem \eqref{eq:homogenized} on a coarse mesh with a larger time step.
Using the equidistant time step $\Delta t = t^* / N$ to discretize the time domain $(0,t^*)$ as
$0=t_0<t_1<\cdots<t_N=t^*$ with $t_n=n\Delta t$ $(n=0,\cdots,N)$. For a given temperature $T_0$, let
$\{\bar T_0^{(s)}\}_{s=0}^{N_T}$ be the off-line representative temperatures with precomputed parameters
$\bar S(\bar T_0^{(s)})$, $\bar k_{ij}(\bar T_0^{(s)})$, $\bar\beta(\bar T_0^{(s)})$.
For $\bar T_0^{(s)} \leq T_0 \leq \bar T_0^{(s+1)}$, the parameters at
arbitrary $T_0$ are obtained by the linear interpolation as below.
\begin{equation}
\label{eq:interpolation}
\bar S(T_0) = \frac{\bar T_0^{(s+1)} - T_0}{\bar T_0^{(s+1)} - \bar T_0^{(s)}} \bar S(\bar T_0^{(s)}) +
\frac{T_0 - \bar T_0^{(s)}}{\bar T_0^{(s+1)} - \bar T_0^{(s)}} \bar S(\bar T_0^{(s+1)}).
\end{equation}
The same linear interpolation approach is applied to compute $\bar k_{ij}(T_0)$ and $\bar\beta(T_0)$. Furthermore, the fully discrete scheme then reads: find $T_0^{n+1} \in V_{h_0}(\Omega)$ such that
\begin{equation}
\begin{aligned}
& \int_{\Omega} \bar S(T_0^{n+1}) \frac{T_0^{n+1}-T_0^{n}}{\Delta t} \upsilon \, \mathrm{d}\Omega
 + \int_{\Omega} \bar k_{ij}(T_0^{n+1}) \frac{\partial T_0^{n+1}}{\partial x_j} \frac{\partial \upsilon}{\partial x_i} \, \mathrm{d}\Omega \\
& + \int_{\Omega} \bar\beta(T_0^{n+1}) \sigma_{\mathrm{B}} (T_0^{n+1})^4 \upsilon \, \mathrm{d}\Omega  = \int_{\Omega} h(\bm{x},t_{n+1}) \upsilon \, \mathrm{d}\Omega\\
&+ \int_{\partial\Omega_q} \hat q(\bm{x},t_{n+1}) \upsilon \, \mathrm{d}s, \quad \forall \upsilon \in V_{h_0}(\Omega).
\label{eq:fullyDiscreteScheme}
\end{aligned}
\end{equation}
\item[(3)]
The fully discrete formulation at each time step leads to a nonlinear algebraic system arising from the temperature-dependent material parameters and nonlinear source terms. To efficiently solve this system, a fixed-point iteration scheme is adopted. At the ($k+1$)-th iterative step, all temperature-dependent material parameters are updated by linear interpolation based on the previous iterate $T_0^{n+1,(k)}$. In addition, the nonlinear radiative term is linearized by lagging the cubic factor: $\left(T_0^{n+1}\right)^4\approx\big(T_0^{n+1,(k)}\big)^3 \, T_0^{n+1,(k+1)}$. The detailed computational scheme is presented as below.
\begin{equation}
\begin{aligned}
& \int_{\Omega} \bar S(T_0^{n+1,(k)}) \frac{T_0^{n+1,(k+1)}-T_0^{n}}{\Delta t} \upsilon \, \mathrm{d}\Omega \\
& + \int_{\Omega} \bar k_{ij}(T_0^{n+1,(k)}) \frac{\partial T_0^{n+1,(k+1)}}{\partial x_j} \frac{\partial \upsilon}{\partial x_i} \, \mathrm{d}\Omega \\
& + \int_{\Omega} \bar\beta(T_0^{n+1,(k)}) \sigma_{\mathrm{B}} (T_0^{n+1,(k)})^{3} T_0^{n+1,(k+1)} \upsilon \, \mathrm{d}\Omega \\
& = \int_{\Omega} h(\bm{x},t_{n+1}) \upsilon \, \mathrm{d}\Omega
 + \int_{\partial\Omega_q} \hat q(\bm{x},t_{n+1}) \upsilon \, \mathrm{d}s,
\quad \forall \upsilon \in V_{h_0}(\Omega).
\label{eq:fixedPointIteration}
\end{aligned}
\end{equation}
To stabilize the fixed-point iteration, a relaxation strategy is employed as follows:
\begin{equation}
T_0^{n+1,(k+1)}=\alpha\,T_0^{n+1,(k+1)} + (1-\alpha)\,T_0^{n+1,(k)}.
\end{equation}
The iteration continues until
$\| T_0^{n+1,(k+1)} - T_0^{n+1,(k)} \|_{L^\infty(\Omega)}
\leq E_{\mathrm{tol}}$, where $E_{\mathrm{tol}}$ denotes the prescribed iteration threshold.
\end{enumerate}

\subsection{On-line stage: reconstruction for high-order multi-scale solutions}
\begin{enumerate}
\item[(1)]
For any point $(\bm{x},t) \in \Omega \times (0,t^*)$, the interpolation technique is employed to retrieve the precomputed auxiliary cell functions and homogenized solutions.
\item[(2)]
The average technique on related elements \cite{R39,R40,R41} is used to evaluate the spatial derivatives $\partial T_0/\partial x_{\alpha_1}$ and $\partial^2 T_0/\partial x_{\alpha_1}\partial x_{\alpha_2}$, and the difference scheme is utilized to evaluate the time derivative $\partial T_0/\partial t$ at each temporal step.
\item[(3)]
Finally, the high-order multi-scale temperature field $T^{(2\varepsilon)}(\bm{x},t)$ is computed by the
formula \eqref{eq:HOMS}, high-order interpolation and post-processing techniques
can be utilized to obtain even more accurate multi-scale solutions \cite{R41}.
\end{enumerate}

\begin{rmk}
In this work, representative macroscopic temperatures are selected as equidistant values across the temperature range. The quantity of these representative points is determined by balancing computational accuracy against efficiency.
\end{rmk}

\section{The error estimate for two-stage multi-scale algorithm}
\begin{lemma}
Let $M_{\alpha_1}^{h_1}$ denote the corresponding finite element solution for the microscopic cell function $M_{\alpha_1}$. If all microscopic cell functions belong to $H^2(\Omega_{\bm{y}})$ for any fixed $T_0$, then the following inequalities hold.
\begin{equation}
\begin{aligned}
&\left\|{M}_{\alpha_1}^{h_1}(\bm{y},T_0)-{M}_{\alpha_1}(\bm{y},T_0)\right\|_{H^m(\Omega_{\bm{y}})}\leq Ch_1^{2-m}\left\| {M}_{\alpha_1}(\bm{y},T_0)\right\|_{H^2(\Omega_{\bm{y}})},
\end{aligned}
\end{equation}
where $m=0,1$ and $C$ is the finite element estimate constant independent of $h_1$ and dependent on ${\Omega_{\bm{y}}}$.
\end{lemma}
$\mathbf{Proof:}$ The inequality (50) is derived by employing classical finite element theory.
\begin{lemma}
Denote $\bar k_{ij}^{h_1}(T_{0})$ be the finite element approximation of the corresponding homogenized thermal parameters, the following results hold.
\begin{equation}
\begin{aligned}
&\left|\bar k_{ij}^{h_1}(T_{0})-\bar k_{ij}(T_{0})\right|\leq Ch_1^2\left\| {M}_i(\bm{y},T_0)\right\|_{H^2(\Omega_{\bm{y}})}\left\|{M}_j(\bm{y},T_0)\right\|_{H^2(\Omega_{\bm{y}})},\\
&\widetilde\kappa_0|\bm{\xi}|^2\leq \bar k_{ij}^{h_1}(T_{0})\xi_i\xi_j\leq\widetilde\kappa_1|\bm{\xi}|^2,
\end{aligned}
\end{equation}
where $C$ is a constant independent of $h_1$.
\end{lemma}
$\mathbf{Proof:}$ By employing the definition of macroscopic homogenized thermal parameters given in (15), along with assumption $A_1$, Remark 2 and Lemma 2, it can be concluded that
\begin{equation}
\begin{aligned}
&\left|\bar k_{ij}^{h_1}(T_{0})-\bar k_{ij}(T_{0})\right|\\
&=\left|\frac{1}{|\Omega_{\bm{y}}|}{\int_{\Omega_{\bm{y}}}}\big({k_{ij}^{(0)} + k_{i\alpha_1}^{(0)}{\frac{\partial M_j^{h_1}}{\partial y_{\alpha_1}}}}\big)\mathrm{d}\Omega_{\bm{y}}-\frac{1}{|\Omega_{\bm{y}}|}{\int_{\Omega_{\bm{y}}}}\big({k_{ij}^{(0)} + k_{i\alpha_1}^{(0)}{\frac{\partial M_j}{\partial y_{\alpha_1}}}}\big)\mathrm{d}\Omega_{\bm{y}}\right|\\
&=\left|\frac{1}{|\Omega_{\bm{y}}|}{\int_{\Omega_{\bm{y}}}}{k_{i\alpha_1}^{(0)}{\frac{\partial \big(M_j^{h_1}-M_j\big)}{\partial y_{\alpha_1}}}}\mathrm{d}\Omega_{\bm{y}}\right|\\
&=\frac{1}{|\Omega_{\bm{y}}|}\left|-{\int_{\Omega_{\bm{y}}}}\frac{\partial M_i}{\partial y_{\alpha_1}}k_{\alpha_1\alpha_2}^{(0)}{\frac{\partial }{\partial y_{\alpha_2}}\big(M_j^{h_1}-M_j\big)}\mathrm{d}\Omega_{\bm{y}}\right|\\
&=\frac{1}{|\Omega_{\bm{y}}|}\left|{\int_{\Omega_{\bm{y}}}}\frac{\partial M_i^{h_1}}{\partial y_{\alpha_1}}k_{\alpha_1\alpha_2}^{(0)}{\frac{\partial }{\partial y_{\alpha_2}}\!\big(M_j^{h_1}\!-\!M_j\big)}\mathrm{d}\Omega_{\bm{y}}\!-\!{\int_{\Omega_{\bm{y}}}}\frac{\partial M_i}{\partial y_{\alpha_1}}k_{\alpha_1\alpha_2}^{(0)}{\frac{\partial }{\partial y_{\alpha_2}}\!\big(M_j^{h_1}\!-\!M_j\big)}\mathrm{d}\Omega_{\bm{y}}\right|\\
&=\frac{1}{|\Omega_{\bm{y}}|}\left|{\int_{\Omega_{\bm{y}}}}\frac{\partial}{\partial y_{\alpha_1}}\big(M_i^{h_1}-M_i\big)k_{\alpha_1\alpha_2}^{(0)}{\frac{\partial }{\partial y_{\alpha_2}}\big(M_j^{h_1}-M_j\big)}\mathrm{d}\Omega_{\bm{y}}\right|\\
&\leq C\left\|M_i^{h_1}-M_i\right\|_{H^1(\Omega_{\bm{y}})}\left\|M_j^{h_1}-M_j\right\|_{H^1(\Omega_{\bm{y}})}\leq Ch_1^2\left\|M_i\right\|_{H^2(\Omega_{\bm{y}})}\left\|M_j\right\|_{H^2(\Omega_{\bm{y}})}.
\end{aligned}
\end{equation}
Furthermore, choosing a sufficiently small $h_1>0$ satisfies
\begin{equation}
Ch_1^2\left\| M_i(\bm{y},T_0)\right\|_{H^2(\Omega_{\bm{y}})}\left\| M_j(\bm{y},T_0)\right\|_{H^2(\Omega_{\bm{y}})}\leq\bar\kappa_0/2.
\end{equation}
Hence, we can verify that the lower bound in (51) holds
\begin{equation}
\bar k_{ij}^{h_1}(T_{0})\xi_i\xi_j=\bar k_{ij}(T_{0})\xi_i\xi_j+\big[\bar k_{ij}^{h_1}(T_{0})-\bar k_{ij}(T_{0})\big]\xi_i\xi_j\geq(\bar\kappa_0-\bar\kappa_0/2)\xi_i\xi_i=\widetilde\kappa_0|\bm{\xi}|^2,
\end{equation}
where $\widetilde\kappa_0=\bar\kappa_0/2$ is a constant independent of $h_1$. Finally, the upper bound in (51) is easily derived when setting $\widetilde\kappa_1=\bar\kappa_1+\bar\kappa_0/2$.

As shown in Lemmas 2 and 3, the values of macroscopic homogenized thermal parameters $\bar k_{ij}(T_{0})$ depend on the finite element computations of the auxiliary cell functions $M_j(\bm{y},T_0)$. Therefore, in practice, we need to numerically solve the modified homogenized equation as below
\begin{equation}
\left\{ \begin{aligned}
&\bar S({T_0^{h_1}})\frac{\partial T_0^{h_1}(\bm{x},t)}{\partial {t}} - \frac{\partial }{{\partial {x_i}}}\Big( {\bar k_{ij}^{h_1}({T_0^{h_1}})\frac{{\partial {T_0^{h_1}}(\bm{x},t)}}{{\partial {x_j}}}} \Big)\\
&\quad\quad\quad\quad\quad\quad+\bar\beta({T_0^{h_1}})\sigma_{\mathrm{B}}\big[T_0^{h_1}(\bm{x},t)\big]^4=h(\bm{x},t),\;\;\text{in}\;\;\Omega\times(0,t^*),\\
&T_0^{h_1}(\bm{x},t) = \hat T(\bm{x},t),\;\;\text{on}\;\;\partial {\Omega_T}\times(0,t^*),\\
&{\bar k_{ij}(T_0^{h_1})\frac{\partial T_0^{h_1}(\bm{x},t)}{\partial {x_j}}} {n_i} = \hat q(\bm{x},t),\;\;\text{on}\;\;\partial {\Omega_q}\times(0,t^*),\\
&T_0^{h_1}(\bm{x},0) = \widetilde T(\bm{x}),\;\;\text{on}\;\;{\Omega}.
\end{aligned} \right.
\label{eq:homogenized1}
\end{equation}
\begin{lemma}
Denote $T_{0}^{h_1}$ be the exact solution of the modified homogenized equation (55), the following estimate holds when assuming ${T_0}\in W^{1,\infty}(\Omega)$
\begin{equation}
{\begin{aligned}
\big\|{T_0^{h_1}}-{T_0}\big\|_{L^\infty(0,t^*;L^2(\Omega))}\leq Ch_1^2,
\end{aligned}}
\end{equation}
where $C$ is a constant independent of $h_1$.
\end{lemma}
$\mathbf{Proof:}$ Through subtracting the macroscopic homogenized equation (14) from corresponding modified homogenized equation (55), one can directly check that
\begin{equation}
\begin{aligned}
&\bar S(T_{0}^{h_1})\frac{{\partial \big(T_0^{h_1}-T_0\big)}}{{\partial t}}- \frac{\partial }{{\partial {x_i}}}\Big[ {{\bar k_{ij}^{h_1}}(T_{0}^{h_1})\frac{{\partial \big(T_0^{h_1}-T_0\big)}}{{\partial {x_j}}}}\Big]=\big[\bar S(T_{0})-\bar S(T_{0}^{h_1})\big]\frac{{\partial T_0}}{{\partial t}}\\
&-\frac{\partial }{{\partial {x_i}}}\Big[ {{\big(\bar k_{ij}(T_{0})-\bar k_{ij}^{h_1}(T_{0})\big)}\frac{{\partial T_0}}{{\partial {x_j}}}}\Big]- \frac{\partial }{{\partial {x_i}}}\Big[ {{\big(\bar k_{ij}^{h_1}(T_{0})-\bar k_{ij}^{h_1}(T_{0}^{h_1})\big)}\frac{{\partial T_0}}{{\partial {x_j}}}}\Big]\\
&+ \big[\bar\beta({T_0})-\bar\beta({T_0^{h_1}})\big]\sigma_{\mathrm{B}}\big(T_0\big)^4+\bar\beta({T_0^{h_1}})\sigma_{\mathrm{B}}\Big[\big(T_0\big)^2+\big(T_0^{h_1}\big)^2\Big]\big[T_0+T_0^{h_1}\big]\big[T_0-T_0^{h_1}\big].
\end{aligned}
\end{equation}
Moreover, multiplying on both sides of equality (57) by $T_0^{h_1}-T_0$ and integrating on $\Omega$, it follows that
\begin{equation}
\begin{aligned}
&\frac{1}{2}\frac{\partial}{\partial t}\Big[\int_{\Omega}\bar S(T_{0}^{h_1})\big(T_0^{h_1}-T_0\big)^2\mathrm{d}\Omega\Big]+\int_{\Omega}{{\bar k_{ij}^{h_1}}(T_{0}^{h_1})\frac{{\partial \big(T_0^{h_1}-T_0\big)}}{{\partial {x_j}}}}\frac{\partial \big(T_0^{h_1}-T_0\big) }{{\partial {x_i}}}\mathrm{d}\Omega\\
&=\frac{1}{2}\int_{\Omega}\frac{\partial \bar S(T_{0}^{h_1})}{\partial t}\big(T_0^{h_1}-T_0\big)^2\mathrm{d}\Omega+\int_{\Omega}\big[\bar S(T_{0})-\bar S(T_{0}^{h_1})\big]\frac{{\partial T_0}}{{\partial t}}\big(T_0^{h_1}-T_0\big)\mathrm{d}\Omega\\
&+\int_{\Omega}{{\big(\bar k_{ij}(T_{0})-\bar k_{ij}^{h_1}(T_{0})\big)}\frac{{\partial T_0}}{{\partial {x_j}}}}\frac{\partial \big(T_0^{h_1}-T_0\big)}{{\partial {x_i}}}\mathrm{d}\Omega\\
&+\int_{\Omega}{{\big(\bar k_{ij}^{h_1}(T_{0})-\bar k_{ij}^{h_1}(T_{0}^{h_1})\big)}\frac{{\partial T_0}}{{\partial {x_j}}}}\frac{\partial \big(T_0^{h_1}-T_0\big)}{{\partial {x_i}}}\mathrm{d}\Omega\\
&+\int_{\Omega}\big[\bar\beta({T_0})-\bar\beta({T_0^{h_1}})\big]\sigma_{\mathrm{B}}\big(T_0\big)^4\big(T_0^{h_1}-T_0\big)\mathrm{d}\Omega\\
&+\int_{\Omega}\bar\beta({T_0^{h_1}})\sigma_{\mathrm{B}}\Big[\big(T_0\big)^2+\big(T_0^{h_1}\big)^2\Big]\big[T_0+T_0^{h_1}\big]\big[T_0-T_0^{h_1}\big]\big(T_0^{h_1}-T_0\big)\mathrm{d}\Omega.
\end{aligned}
\end{equation}
Subsequently, integrating both sides of (58) from $0$ to $t$ $(0<t\leq t^*)$, recalling the inequality (51), and employing Cauchy-Schwarz inequality and Young inequality, we can obtain the following inequality from equality (58) if $|\bar S(T_{0})-\bar S(T_{0}^{h_1})|\leq C|T_{0}-T_{0}^{h_1}|$, $|\bar k_{ij}^{h_1}(T_{0})-\bar k_{ij}^{h_1}(T_{0}^{h_1})|\leq C|T_{0}-T_{0}^{h_1}|$ and $|\bar \beta(T_{0})-\bar \beta(T_{0}^{h_1})|\leq C|T_{0}-T_{0}^{h_1}|$.
\begin{equation}
\begin{aligned}
&C\big\|T_0^{h_1}-T_0\big\|_{L^2(\Omega)}^2+\mathlarger{\int}_0^tC\big\|T_0^{h_1}-T_0\big\|_{H^1(\Omega)}^2\mathrm{d}\tau\\
&\leq \mathlarger{\int}_0^t C\big\|T_0^{h_1}-T_0\big\|_{L^2(\Omega)}^2\mathrm{d}\tau+\mathlarger{\int}_0^tC\big\|T_0^{h_1}-T_0\big\|_{L^2(\Omega)}^2\mathrm{d}\tau+\mathlarger{\int}_0^t\frac{C}{2\lambda}h_1^4\mathrm{d}\tau\\
&+\mathlarger{\int}_0^t\frac{\lambda}{2}\big\|T_0^{h_1}-T_0\big\|_{H^1(\Omega)}^2\mathrm{d}\tau+\mathlarger{\int}_0^t\frac{C}{2\lambda}\big\|T_0^{h_1}-T_0\big\|_{L^2(\Omega)}^2\mathrm{d}\tau\\
&+\mathlarger{\int}_0^t\frac{\lambda}{2}\big\|T_0^{h_1}-T_0\big\|_{H^1(\Omega)}^2\mathrm{d}\tau+\mathlarger{\int}_0^t C\big\|T_0^{h_1}-T_0\big\|_{L^2(\Omega)}^2\mathrm{d}\tau+\mathlarger{\int}_0^tC\big\|T_0^{h_1}-T_0\big\|_{L^2(\Omega)}^2\mathrm{d}\tau\\
&\leq \mathlarger{\int}_0^tCh_1^4\mathrm{d}\tau+\mathlarger{\int}_0^tC\big\|T_0^{h_1}-T_0\big\|_{L^2(\Omega)}^2\mathrm{d}\tau+\mathlarger{\int}_0^t{\lambda}\big\|T_0^{h_1}-T_0\big\|_{H^1(\Omega)}^2\mathrm{d}\tau
\end{aligned}
\end{equation}
When choosing a sufficiently small $\lambda$ and setting $\Upsilon(t)=\big\|T_0^{h_1}-T_0\big\|_{L^2(\Omega)}^2$, then we can derive $\Upsilon(t)\leq Ch_1^4+C\mathlarger{\int}_0^t\Upsilon(\tau)\mathrm{d}\tau$ from (59). Consequently, by taking advantage of Gronwall inequality and the arbitrariness of time variable $t$, there holds the inequality (56).
\begin{lemma}
Denote $T_{0}^{h_1,h_0,(k+1)}$ be the finite element solution of the modified homogenized equation (55) by semi-discrete iterative method proposed in reference \cite{R36}, the following estimate holds
\begin{equation}
{\begin{aligned}
\big\|T_{0}^{h_1,h_0,(k+1)}-T_{0}^{h_1}\big\|_{H^1(\Omega)}\leq CE_{\mathrm{tol}}+Ch_0,
\end{aligned}}
\end{equation}
where $C$ is a constant irrespective of $h_0$ and $E_{\mathrm{tol}}$.
\end{lemma}
$\mathbf{Proof:}$ As shown in the modified homogenized equation (55), it satisfies the conditions of error estimates in reference \cite{R36}. Hence, employing the same proof technique as reference \cite{R36}, the estimate (60) can be derived.
\begin{theorem}
Let $T_{0}^{h_1,h_0,(k+1)}$ be the semi-discrete iterative finite element solution of the modified homogenized equation (55), and $T_{0}$ be the exact solution of the macroscopic homogenized equations (14), then the following estimate holds
\begin{equation}
\begin{aligned}
\big\|T_0^{h_1,h_0,(k+1)}-{T_0}\big\|_{L^2(\Omega)}\leq CE_{\mathrm{tol}}+Ch_0+Ch_1^2,
\end{aligned}
\end{equation}
where $C$ is a positive constant independent of $h_1$, $h_0$ and $E_{\mathrm{tol}}$.
\end{theorem}
$\mathbf{Proof:}$ Firstly, employing the triangle inequality, there exists an inequality such that
\begin{equation}
\begin{aligned}
&\big\|T_0^{h_1,h_0,(k+1)}-{T_0}\big\|_{L^2(\Omega)}\\
&\leq \big\|T_0^{h_1,h_0,(k+1)}-{T_0^{h_1}}\big\|_{L^2(\Omega)}+\big\|T_0^{h_1}-{T_0}\big\|_{L^2(\Omega)}.
\end{aligned}
\end{equation}
Moreover, with the help of Lemmas 4 and 5, suffice it to obtain the error estimate (61).

To sum up, we remark that the above-mentioned theoretical analysis rigorously ensures the convergence of the proposed two-stage numerical algorithm in microscopic and macroscopic computation.
\section{Numerical examples and discussions}
In this section, two- and three-dimensional numerical examples are presented to validate the computational performance of the proposed HOMS method and corresponding two-stage algorithm. All numerical experiments are implemented based on FreeFEM++ software and conducted on a desktop workstation equipped with a 12th Gen Intel Core i5-12600KF processor
(3.70\,GHz) and 32\,GB internal memory.

Since it is impossible to obtain the exact solutions for the multi-scale radiative heat transfer problem \eqref{eq:all}, we replace the exact solution $T^{\varepsilon}(\bm{x},t)$ by direct FEM solution on a sufficiently fine mesh as the reference for evaluating the proposed method. For convenience, the following error notations are introduced for macroscopic homogenized solution $T_0$.
\begin{equation}
\mathrm{TL^2err0}(t)=\frac{\|T^{\varepsilon}-T_0\|_{L^2(\Omega)}}{\|T^{\varepsilon}\|_{L^2(\Omega)}},
\mathrm{TH^1err0}(t)=\frac{|T^{\varepsilon}-T_0|_{H^1(\Omega)}}{|T^{\varepsilon}|_{H^1(\Omega)}}.
\end{equation}
Errors $\mathrm{TL^2err1}(t)$ and $\mathrm{TH^1err1}(t)$, and $\mathrm{TL^2err2}(t)$ and $\mathrm{TH^1err2}(t)$ for $T^{(1\varepsilon)}$ and $T^{(2\varepsilon)}$ can be defined similarly.

\subsection{Example 1: Linear radiative heat transfer problem of 2D composite structure}\label{example1}
\begin{figure}[!htb]
\centering
\begin{subfigure}[b]{0.33\textwidth}
  \centering
  \includegraphics[width=\linewidth]{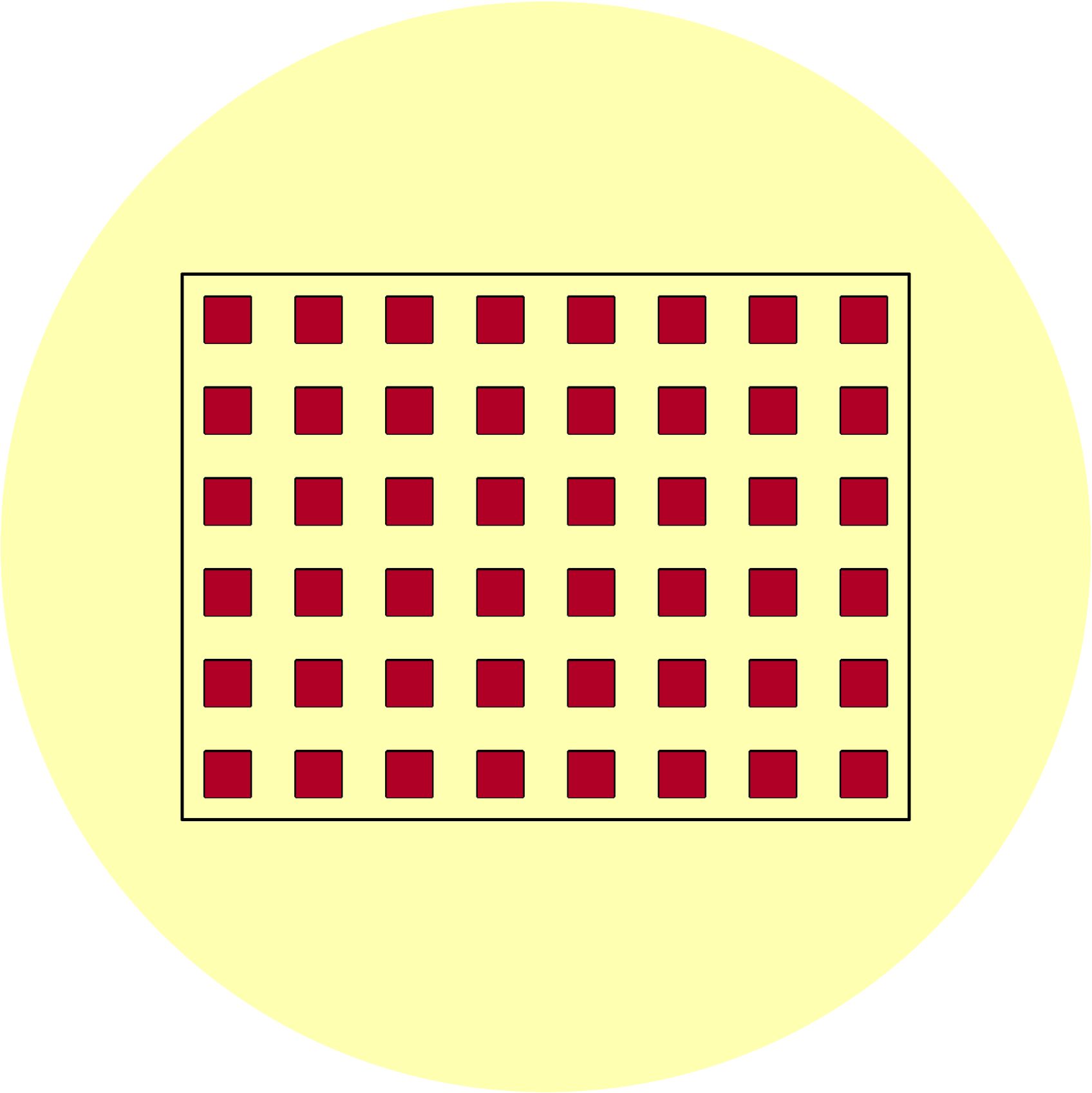}
  \caption{}
  \label{fig:macro}
\end{subfigure}
\hfill
\begin{subfigure}[b]{0.3\textwidth}
  \centering
  \includegraphics[width=\linewidth]{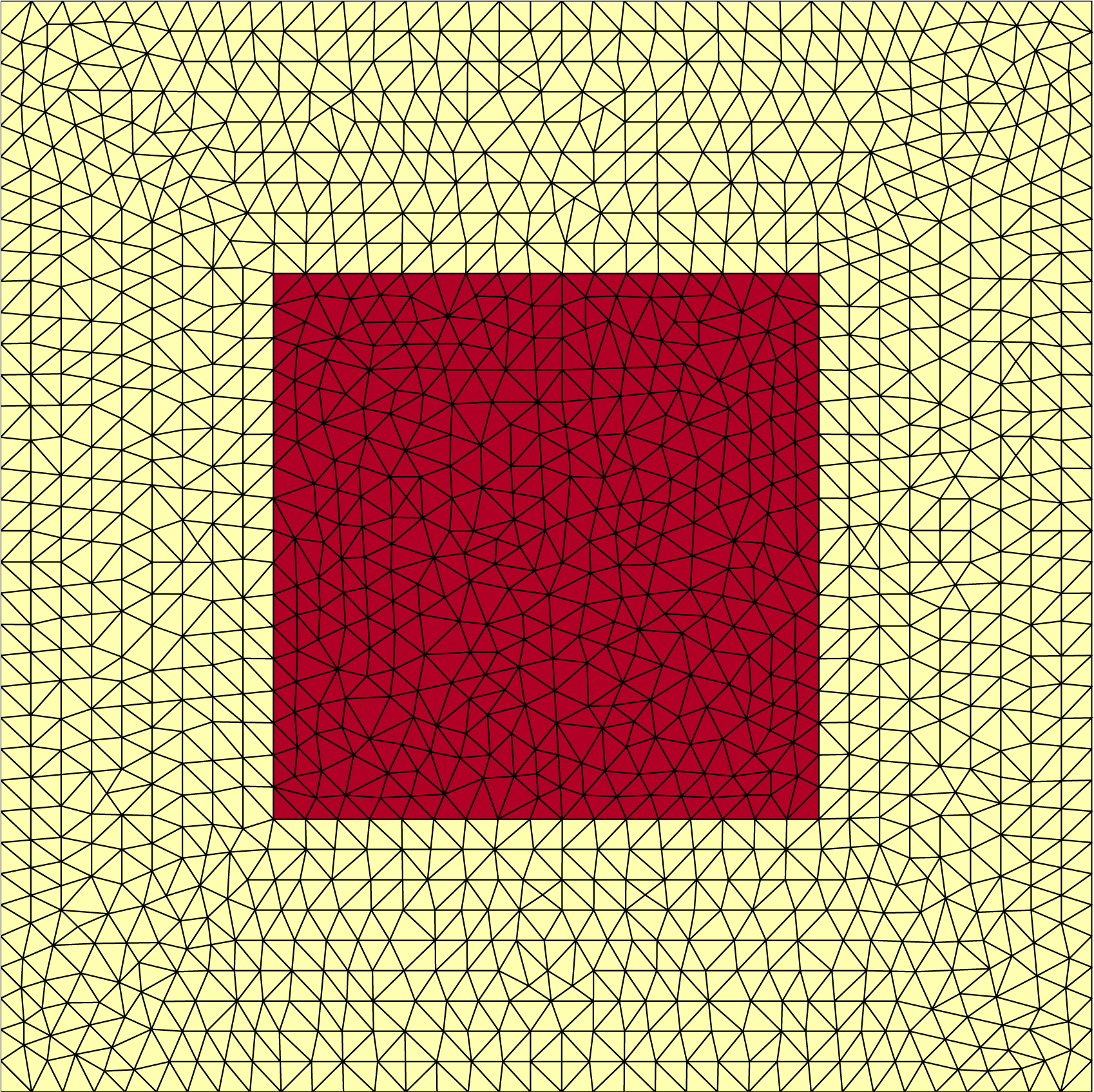}
  \caption{}
  \label{fig:cell}
\end{subfigure}
\hfill
\begin{subfigure}[b]{0.33\textwidth}
  \centering
  \includegraphics[width=\linewidth]{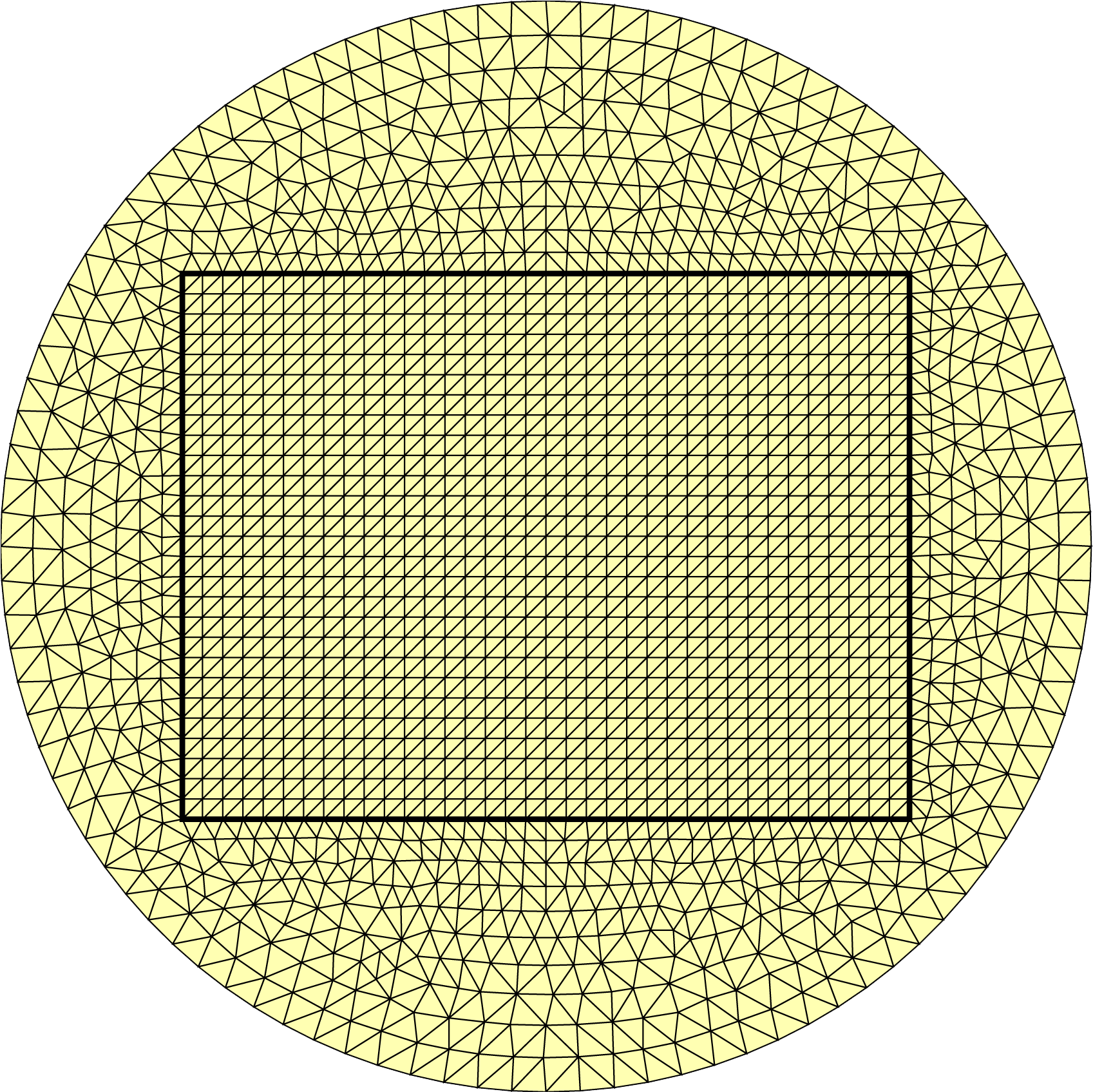}
  \caption{}
  \label{fig:homo}
\end{subfigure}
\caption{(a) The multi-scale structure $\Omega$; (b) the microscopic unit cell $\Omega_{\bm{y}}$; (c) the macroscopic homogenized structure.}
\label{e1structure}
\end{figure}

This example studies the linear radiative heat transfer problem of 2D composite structure. The multi-scale structure $\Omega$ is a circular area centered at $(0.75\mathrm{cm},0.75\mathrm{cm})$ with a radius of 0.75$\mathrm{cm}$, as shown in \autoref{e1structure}. Furthermore, corresponding macroscopic homogenization structure and microscopic unit cell $\Omega_{\bm{y}}$ are shown in \autoref{e1structure}, where $\Omega_{\bm{y}}=(y_1,y_2)=[0,1]\times[0,1]$ and $\varepsilon=1/8$.

This composite consists of two component materials. The yellow material is matrix phase, whereas the red material is inclusion phase. Moreover, the material parameters of both are listed in \autoref{tabe1canshu}.
\begin{table}[htbp]
    \centering
    \caption{Material parameters of 2D composite structure.}    \label{tabe1canshu}
\begin{tabular}{ccc}
    \toprule
    \textbf{Property} & \textbf{Matrix} & $\textbf{Inclusion}$ \\
    \midrule
     $\rho^\varepsilon {\cdot} c^\varepsilon$ ($\mathrm{J}\boldsymbol{\cdot}\mathrm{cm}^{-3}\boldsymbol{\cdot}\mathrm{K}^{-1}$) & $2.41$  & $3.60$  \\
     $k_{ij}^\varepsilon$ ($\mathrm{W}\boldsymbol{\cdot}\mathrm{cm}^{-1}\boldsymbol{\cdot}\mathrm{K}^{-1}$) & $3.20$  & $0.19$   \\
     $\sigma_{\mathrm{B}} {\cdot} \beta^\varepsilon$ ($\mathrm{W}\boldsymbol{\cdot}\mathrm{cm}^{-3}\boldsymbol{\cdot}\mathrm{K}^{-4}$) & $4.96 \times 10^{-12}$  & $1.70 \times 10^{-12}$  \\
    \bottomrule
\end{tabular}
\end{table}

The heat source function and initial-boundary conditions in multi-scale problem \eqref{eq:all} of this example are defined as follows:
\[
h(\bm{x},t) = 1.0\times 10^4\,\mathrm{J}/(\mathrm{cm}^2 \boldsymbol{\cdot} \mathrm{s}),\
\hat T(\bm{x},t) = 300.0\,\mathrm{K},\
\widetilde T(\bm{x}) = 300.0\,\mathrm{K}.
\]

The multi-scale problem \eqref{eq:all}, auxiliary cell problems \eqref{eq:Malpha1}, and \eqref{eq:Malpha1alpha2}--\eqref{eq:Nalpha1alpha2}, and the homogenized problem \eqref{eq:homogenized} are computed with a temporal step $\Delta t = 0.01 \,\mathrm{s}$ over the time interval $t \in [0,1]\,\mathrm{s} $ (the same time step and time interval are adopted for all subsequent examples). The computational cost is listed in \autoref{e1computation}, which depicts the advantage of the proposed approach by substantially saving computer storage and simulation time.
\begin{table}[!htb]
\caption{Comparison of the cost on computation resources.}
\label{e1computation}
\centering
\begin{tabular}{lccc}
\toprule
 & Cell eqs. & Homogenized eqs. & Multi-scale eqs. \\
\midrule
FEM elements & 2618 & 14280 & 55630 \\
FEM nodes & 1380 & 7241 & 27966 \\
\midrule
 &  \multicolumn{2}{c}{HOMS method} & FEM \\
\midrule
Computational time & \multicolumn{2}{c}{768.58s}  & 2560.26s \\
\bottomrule
\end{tabular}
\end{table}
\begin{rmk}
In \autoref{e1computation}, "FEM" denotes the direct finite element method on a fine mesh
for \eqref{eq:all}, while "HOMS" denotes the high-order multi-scale method.
\end{rmk}

\begin{figure}[!htb]
\centering
\begin{minipage}[c]{0.22\textwidth}
  \centering
  \includegraphics[width=\linewidth]{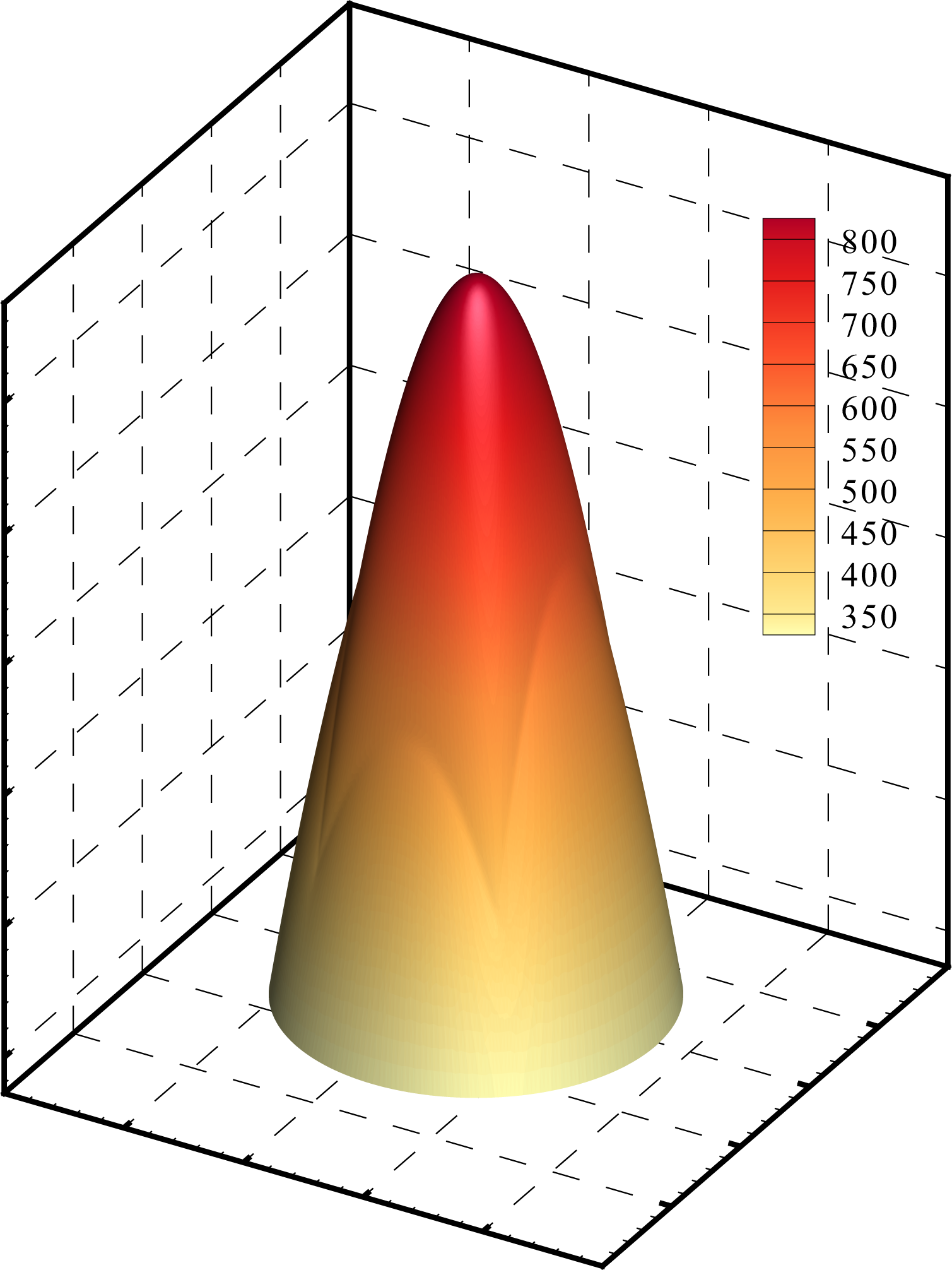}\\(a)
\end{minipage}\hspace{0.01\textwidth}
\begin{minipage}[c]{0.22\textwidth}
  \centering
  \includegraphics[width=\linewidth]{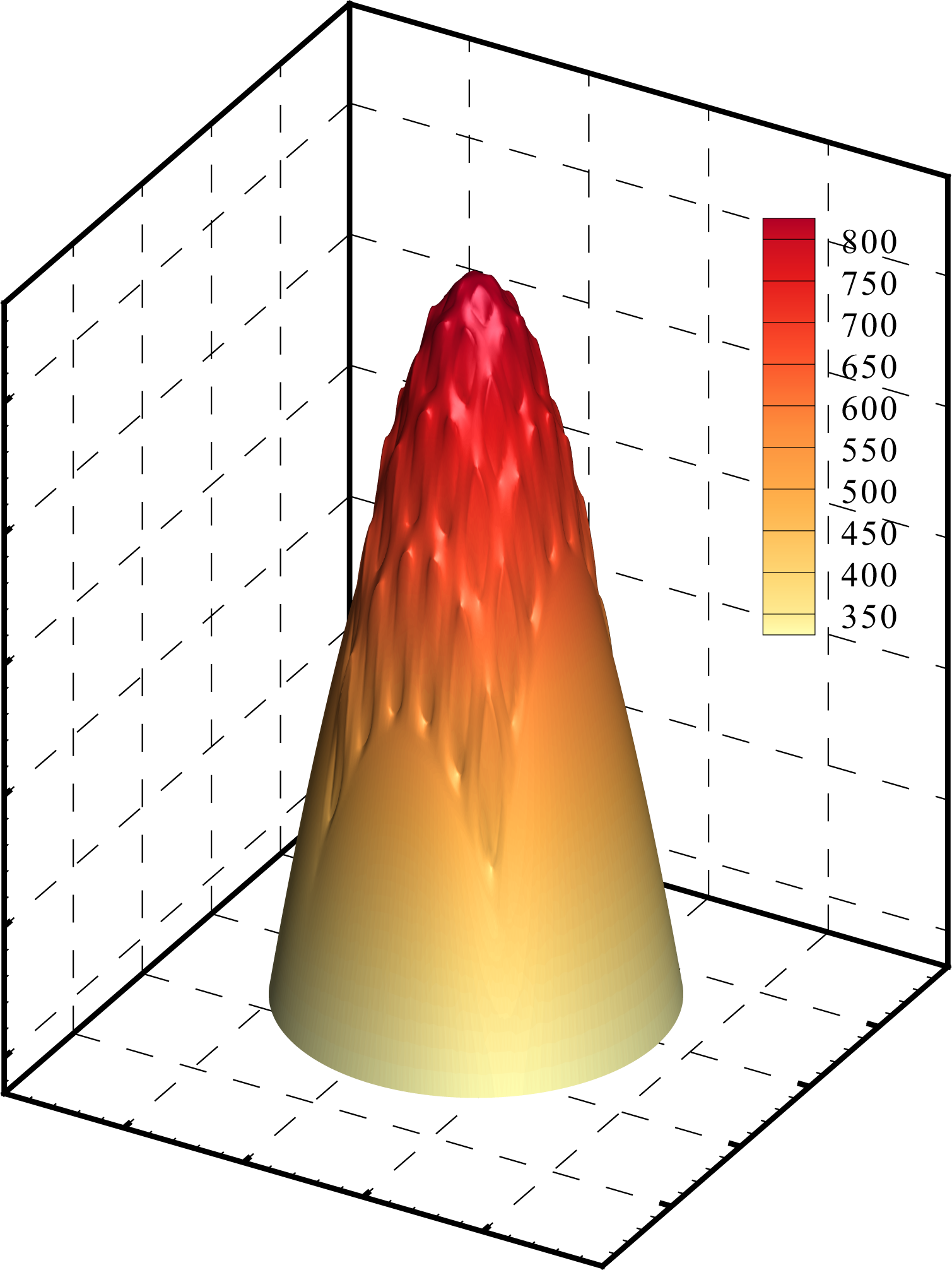}\\(b)
\end{minipage}\hspace{0.01\textwidth}
\begin{minipage}[c]{0.22\textwidth}
  \centering
  \includegraphics[width=\linewidth]{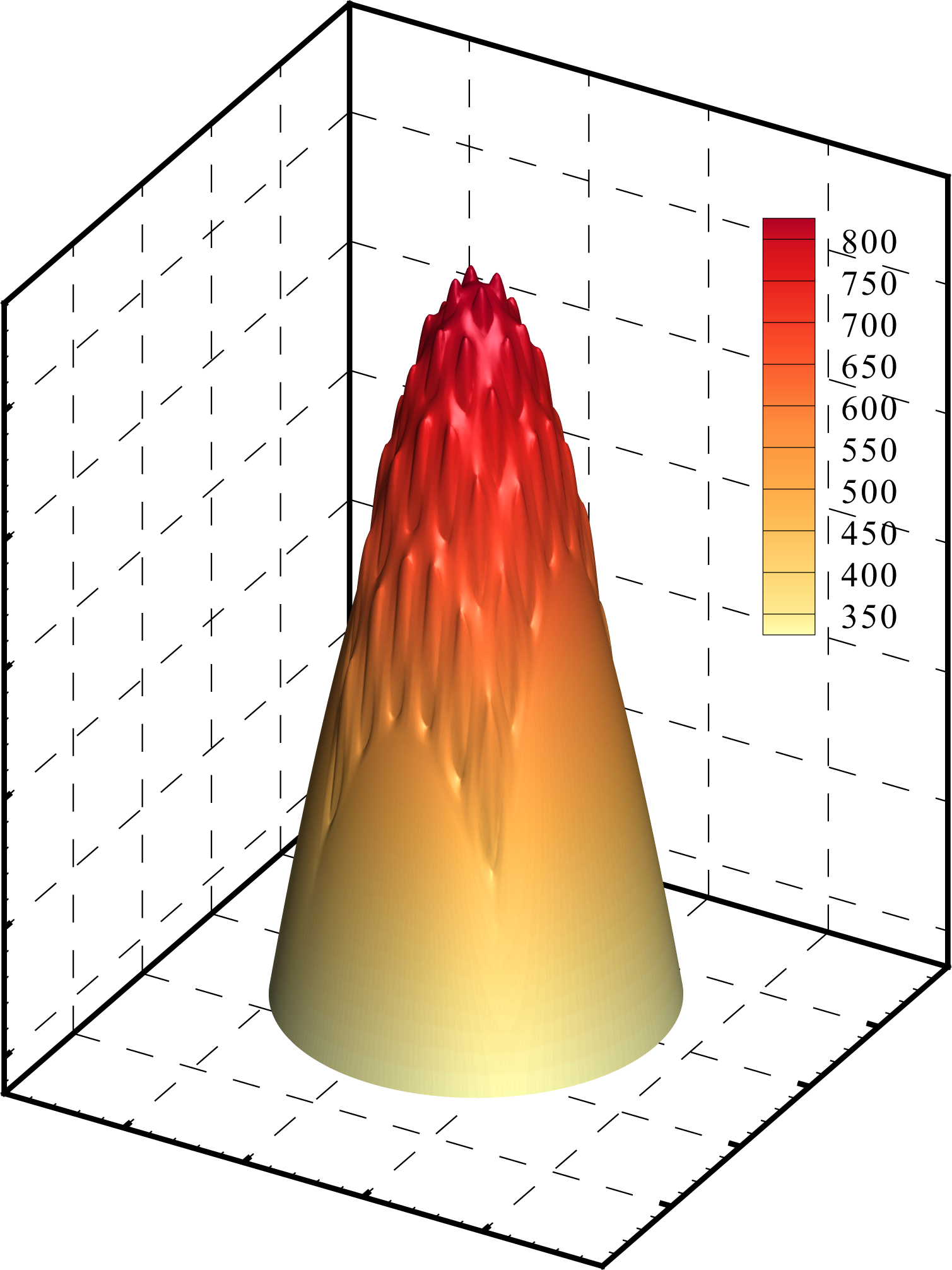}\\(c)
\end{minipage}\hspace{0.01\textwidth}
\begin{minipage}[c]{0.22\textwidth}
  \centering
  \includegraphics[width=\linewidth]{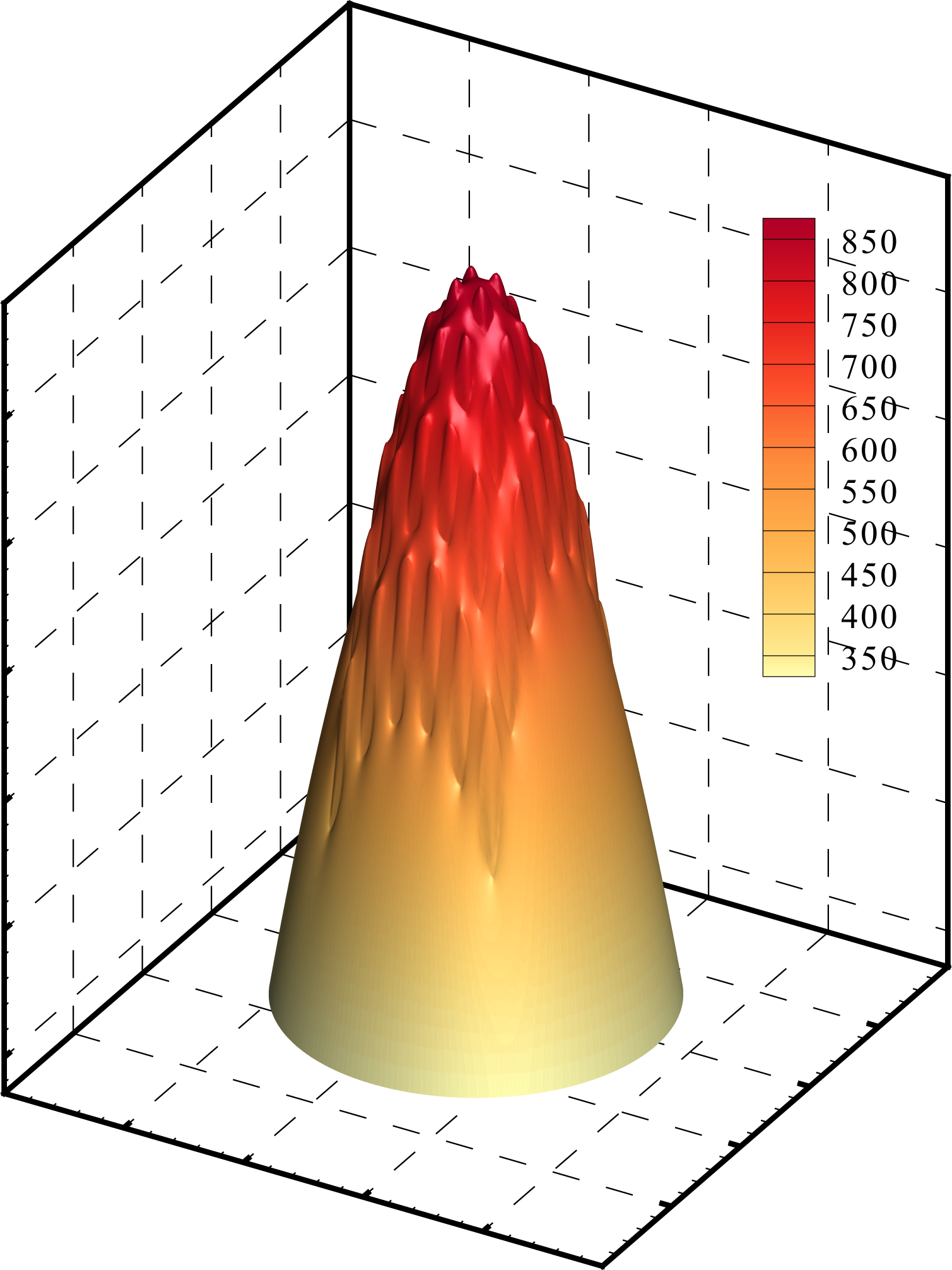}\\(d)
\end{minipage}
\caption{The temperature field at $t=1.0\,\mathrm{s}$: (a) $T_{0}$; (b) $T^{(1{\varepsilon})}$; (c) $T^{(2{\varepsilon})}$; (d) $T^{{\varepsilon}}$.}\label{e1temperature}
\end{figure}

After multi-scale simulation, the resulting temperature fields $T_{0}$, $T^{(1\bm{\varepsilon})}$,
$T^{(2\bm{\varepsilon})}$ and $T^{\bm{\varepsilon}}$ are shown in \autoref{e1temperature}.

\begin{figure}[!htb]
\centering
\begin{subfigure}[c]{0.45\textwidth}
  \centering
  \includegraphics[width=\linewidth]{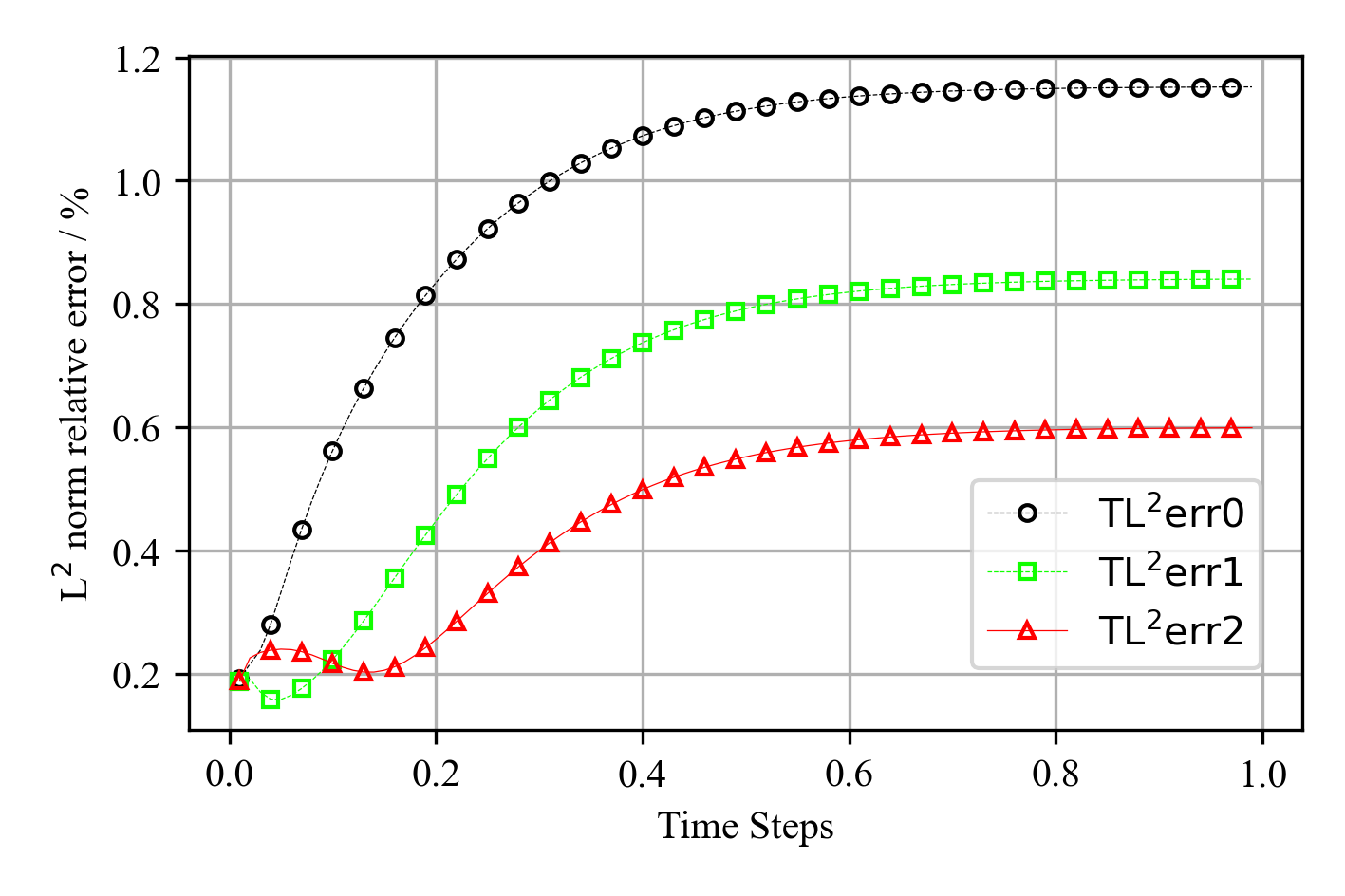}
  \caption{$L^2$ error}
  \label{fig:e1l2-error}
\end{subfigure}
\hfill
\begin{subfigure}[c]{0.45\textwidth}
  \centering
  \includegraphics[width=\linewidth]{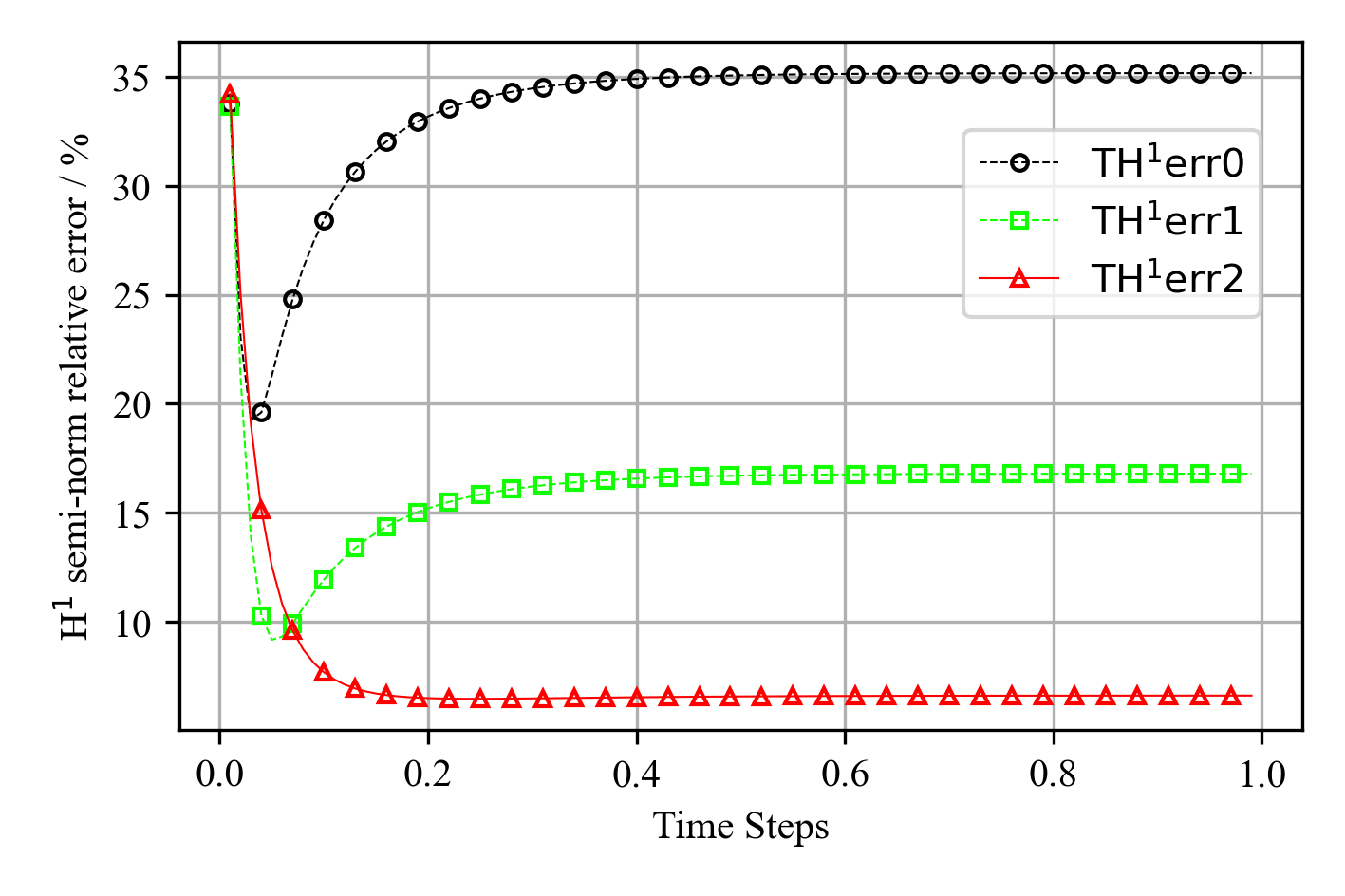}
  \caption{$H^1$ semi-norm error}
  \label{fig:e1h1-error}
\end{subfigure}
\caption{The evolution of relative errors in the $L^2$ norm and $H^1$ semi-norm.}
\label{fig:e1error-comparison}
\end{figure}

Furthermore, the detailed relative errors of the temperature fields are presented in \autoref{fig:e1error-comparison}. The HOMS method achieves an $L^2$ norm relative error of $0.6\%$ and an $H^1$ semi-norm relative error of $3\%$, compared to $0.9\%$ ($L^2$ norm) and $17\%$ ($H^1$ semi-norm) for the first-order method, and $1.1\%$ ($L^2$ norm) and $35\%$ ($H^1$ semi-norm) for the homogenization method. This demonstrates that the HOMS method substantially improves the gradient accuracy over both lower-order approaches. Moreover, the HOMS method can largely save computer memory and CPU time, which is quite important and meaningful in engineering computation.

\autoref{fig:e1error-comparison} also demonstrates that only the HOMS solutions are almost identical to the precise FEM
solutions and can accurately capture the steep oscillating information at the micro-scale arising from the heterogeneities
in the 2D composite, especially for the thermal field. By contrast, according to the error results in
\autoref{fig:e1error-comparison}, the homogenized and LOMS solutions are far from sufficient to provide a high-accuracy
solution for the multi-scale problem \eqref{eq:all}. This observation indicates that the HOMS solutions are vital for accurately capturing the microscopic oscillating information of the 2D composite. Consequently, it can be concluded that our HOMS numerical algorithm is stable and efficient for simulating the multi-scale problems of the 2D composite.

\subsection{Example 2: Nonlinear radiative heat transfer problem of 2D composite structure with temperature-dependent properties }\label{example2}
\begin{figure}[!htb]
\centering
\begin{minipage}[c]{0.3\textwidth}
  \centering
  \includegraphics[width=\linewidth]{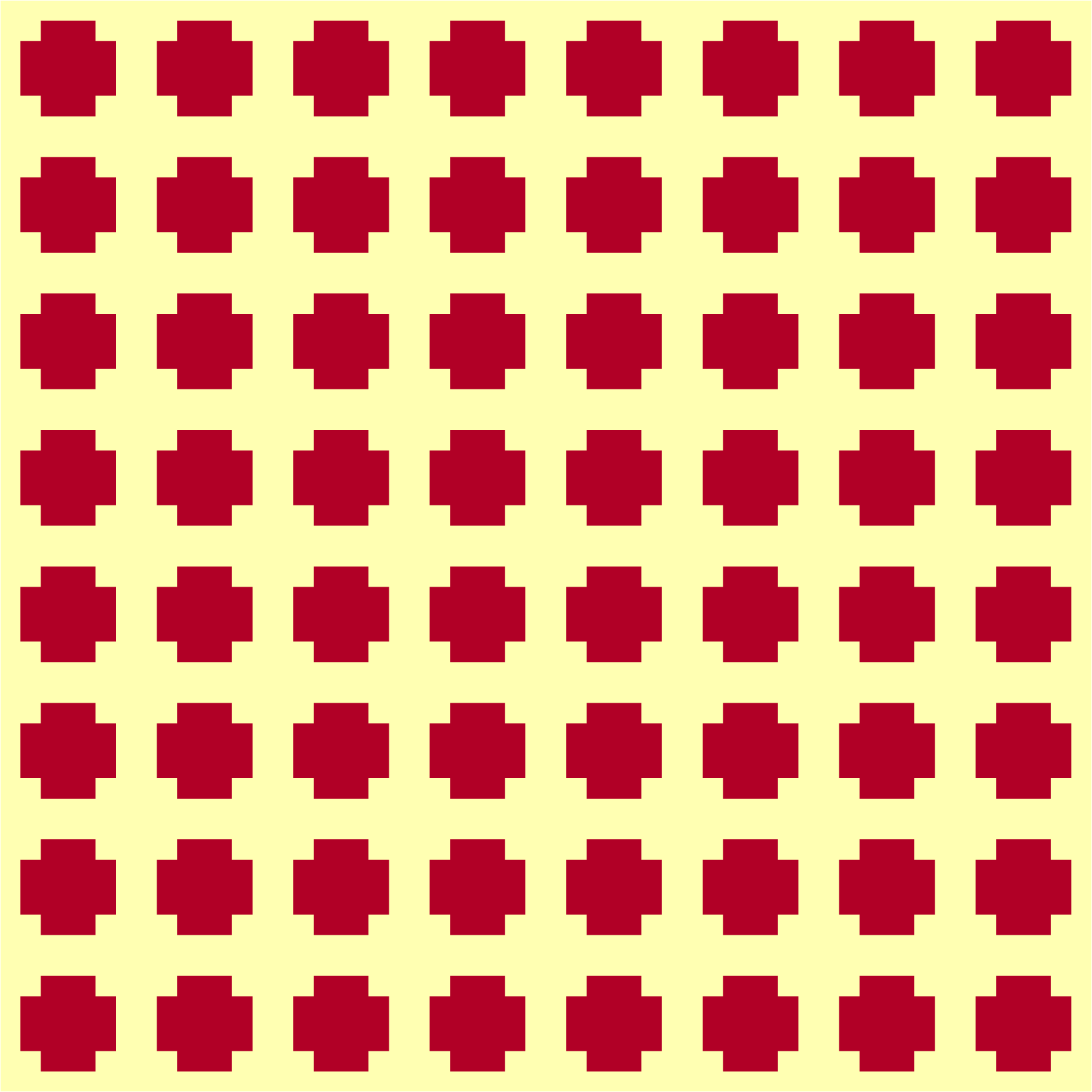}\\
  (a)
\end{minipage}
\hfill
\begin{minipage}[c]{0.3\textwidth}
  \centering
  \includegraphics[width=\linewidth]{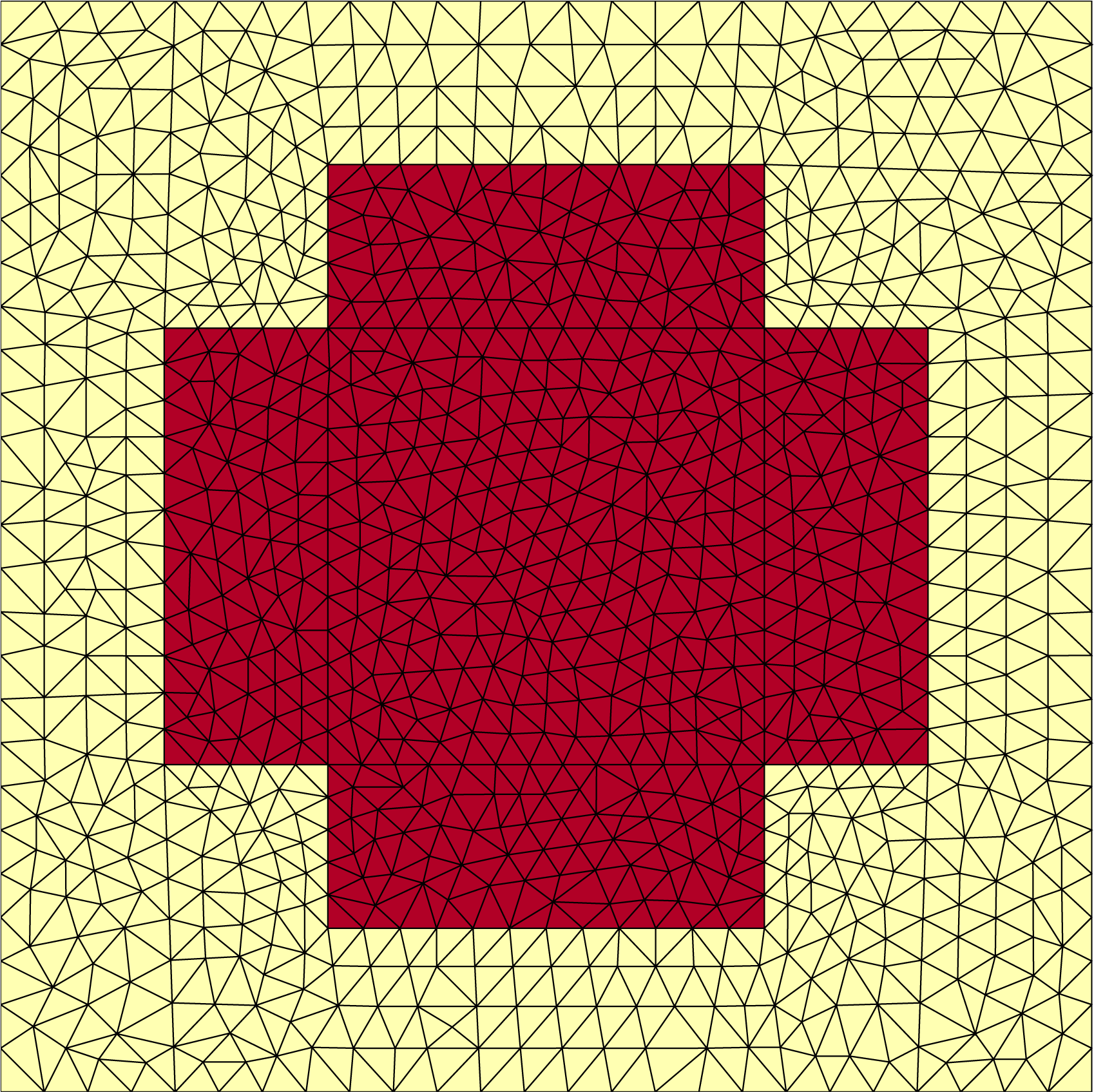}\\
  (b)
\end{minipage}
\hfill
\begin{minipage}[c]{0.3\textwidth}
  \centering
  \includegraphics[width=\linewidth]{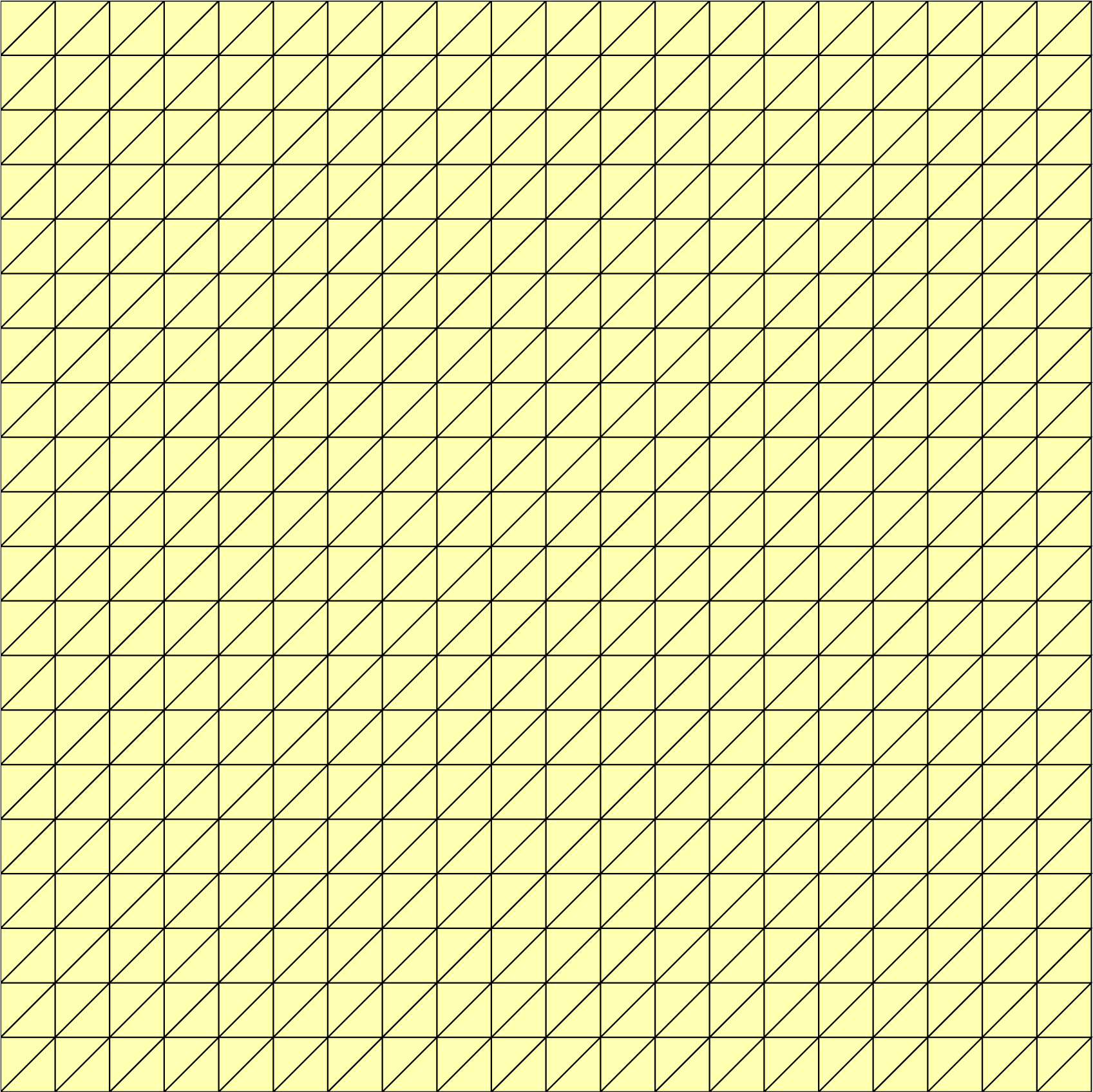}\\
  (c)
\end{minipage}
\caption{(a) The multi-scale structure $\Omega$; (b) the microscopic unit cell $\Omega_{\bm{y}}$; (c) the macroscopic homogenized structure.}
\label{2dstructure}
\end{figure}

This example studies the nonlinear radiative heat transfer problem of 2D composite structure with temperature-dependent properties. The multi-scale structure $\Omega$, macroscopic homogenization structure and microscopic unit cell $\Omega_{\bm{y}}$ are shown in \autoref{2dstructure}, where $\Omega=(x_1,x_2)=[0,1]\times[0,1] \mathrm{cm}^2$, and $\Omega_{\bm{y}}=(y_1,y_2)=[0,1]\times[0,1]$.

This composite structure consists of two component materials, where the characteristic periodic parameter is set to $\varepsilon=1/8$. The material parameters of the composite are defined in \autoref{tab2dcanshu}.
\begin{table}[htbp]
    \centering
    \caption{Material parameters of 2D composite structure}    \label{tab2dcanshu}
\begin{tabular}{ccc}
    \toprule
    \textbf{Property} & \textbf{Matrix} & \textbf{Inclusion} \\
    \midrule
     $\rho^\varepsilon \boldsymbol{\cdot} c^\varepsilon$ ($\mathrm{J}\boldsymbol{\cdot}\mathrm{cm}^{-3}\boldsymbol{\cdot}\mathrm{K}^{-1}$) & $4.5$  & $1.5$  \\
     $k_{ij}^\varepsilon$ ($\mathrm{W}\boldsymbol{\cdot}\mathrm{cm}^{-1}\boldsymbol{\cdot}\mathrm{K}^{-1}$) & $100.0+4.0\times10^{-3}T$  & $1.0+4.0\times10^{-5}T$   \\
     $\sigma_{\mathrm{B}} \boldsymbol{\cdot} \beta^\varepsilon$ ($\mathrm{W}\boldsymbol{\cdot}\mathrm{cm}^{-3}\boldsymbol{\cdot}\mathrm{K}^{-4}$) & $1.0\times10^{-6}$  & $1.0\times10^{-7}$  \\
    \bottomrule
\end{tabular}
\end{table}

The heat source function and initial-boundary conditions in multi-scale problem \eqref{eq:all} of this example are defined as follows:
\[
h(\bm{x},t) = 1.0\times 10^6\,\mathrm{J}/(\mathrm{cm}^2 \cdot \mathrm{s}),\
\hat T (\bm{x},t)= 300.0\,\mathrm{K},\
\widetilde T (\bm{x})= 300.0\,\mathrm{K}.
\]

The computational cost of this example is listed in \autoref{2dcomputation}, which illustrates the advantages of the proposed method by providing a tremendous saving in computing resource.
\begin{table}[!htb]
\caption{Comparison of the cost on computation resources.}
\label{2dcomputation}
\centering
\begin{tabular}{lccc}
\toprule
 & Cell eqs. & Homogenized eqs. & Multi-scale eqs. \\
\midrule
FEM elements & 8186 & 7200 & 102400 \\
FEM nodes & 4194 & 3721 & 51521 \\
\midrule
 & Off-line stage & On-line stage & FEM \\
\midrule
Computational time & 68.912s & 364.342s & 543.936s \\
\bottomrule
\end{tabular}
\end{table}
\begin{rmk}
In \autoref{2dcomputation}, "FEM" denotes the direct finite element simulation on a fine mesh for \eqref{eq:all}. The "Off-line stage" comprises the precomputation of cell problems (see \autoref{sec:offline}), and the "On-line stage" comprises the interpolation, the fully discrete scheme, and other related numerical procedures (see \autoref{sec:online}). The same notation is used in all subsequent tables (\autoref{ex3dcomputation} and \autoref{ex43dcomputation}).
\end{rmk}

\begin{figure}[!htb]
\centering
\begin{minipage}[c]{0.22\textwidth}
  \centering
  \includegraphics[width=\linewidth]{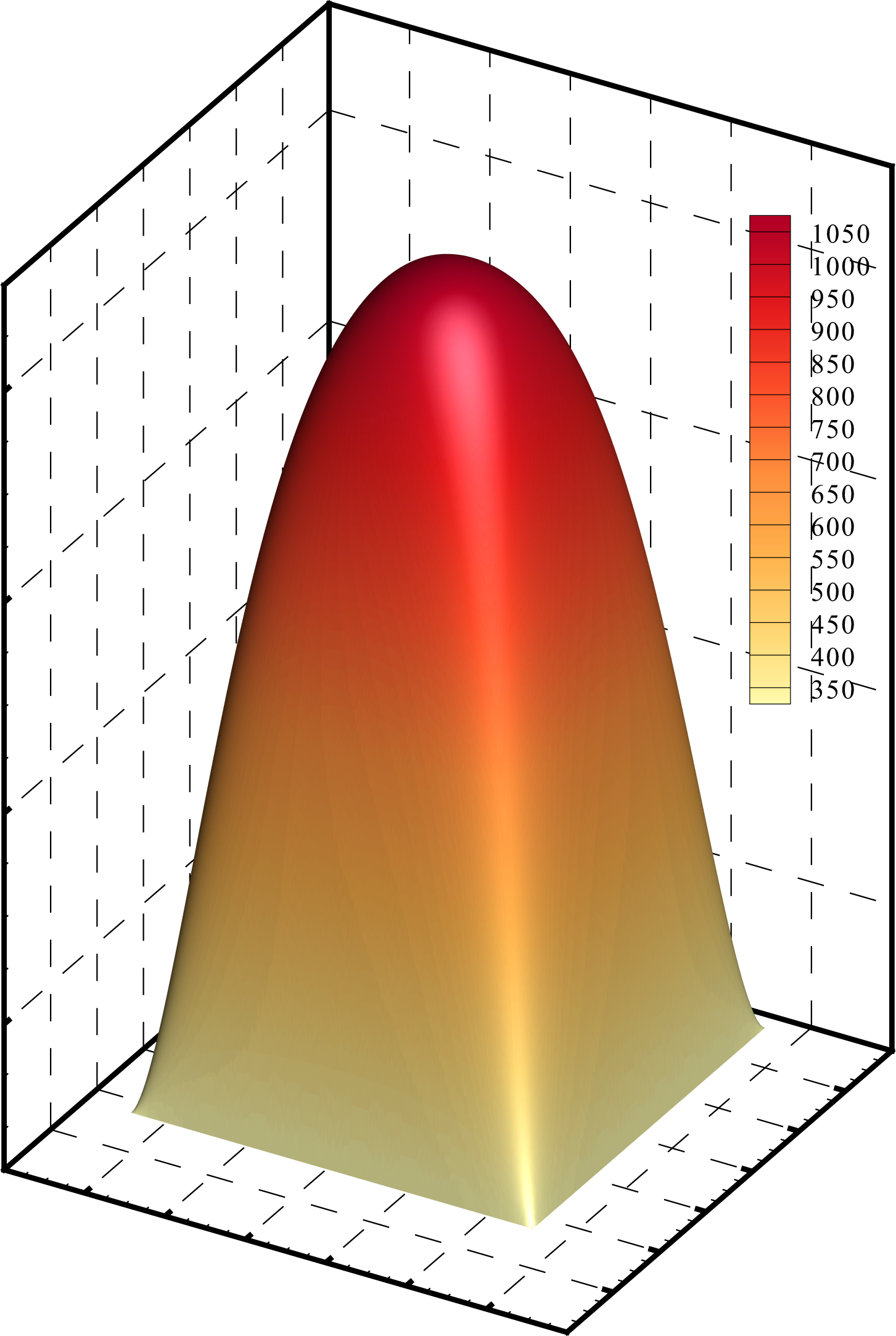}\\(a)
\end{minipage}\hspace{0.01\textwidth}
\begin{minipage}[c]{0.22\textwidth}
  \centering
  \includegraphics[width=\linewidth]{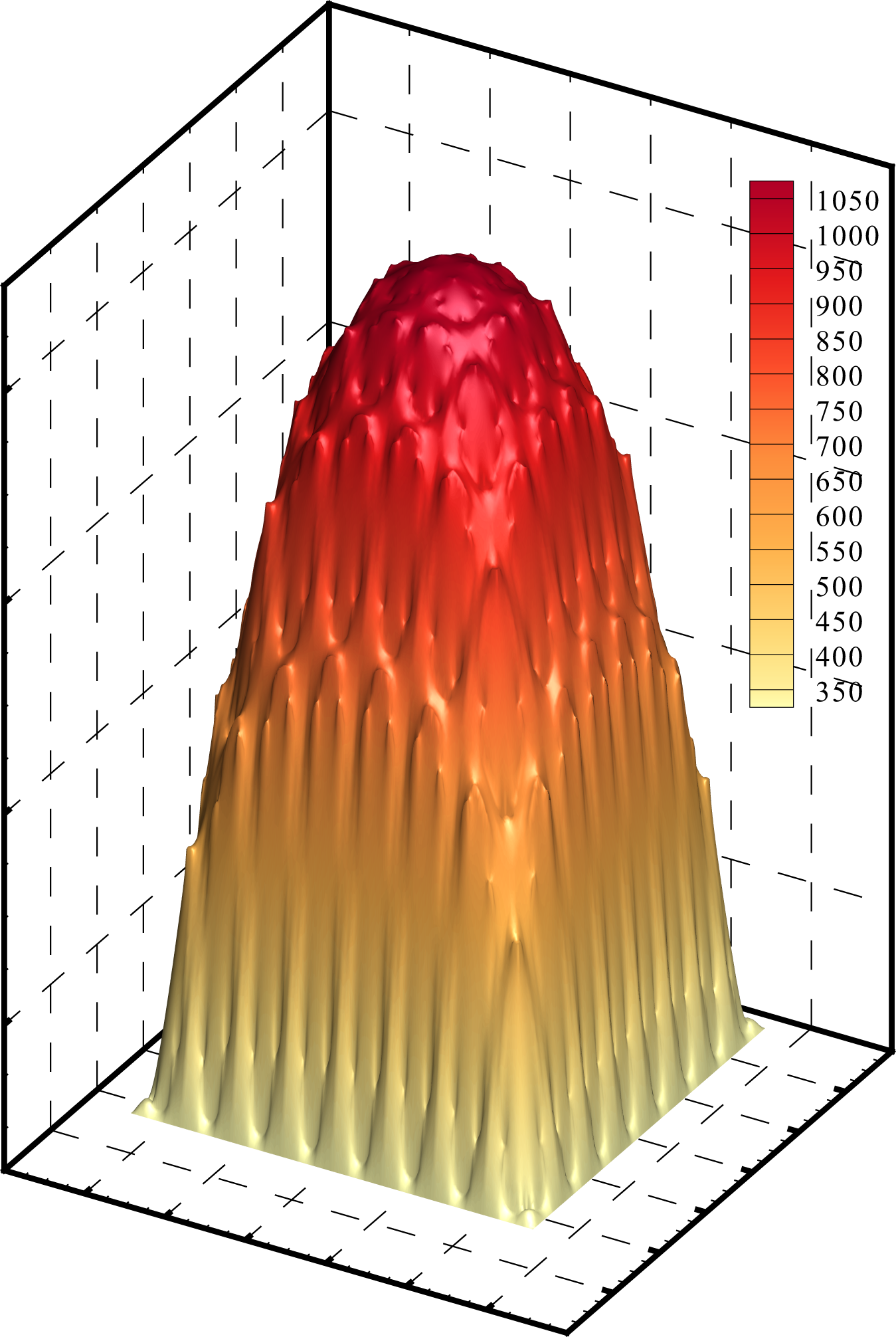}\\(b)
\end{minipage}\hspace{0.01\textwidth}
\begin{minipage}[c]{0.22\textwidth}
  \centering
  \includegraphics[width=\linewidth]{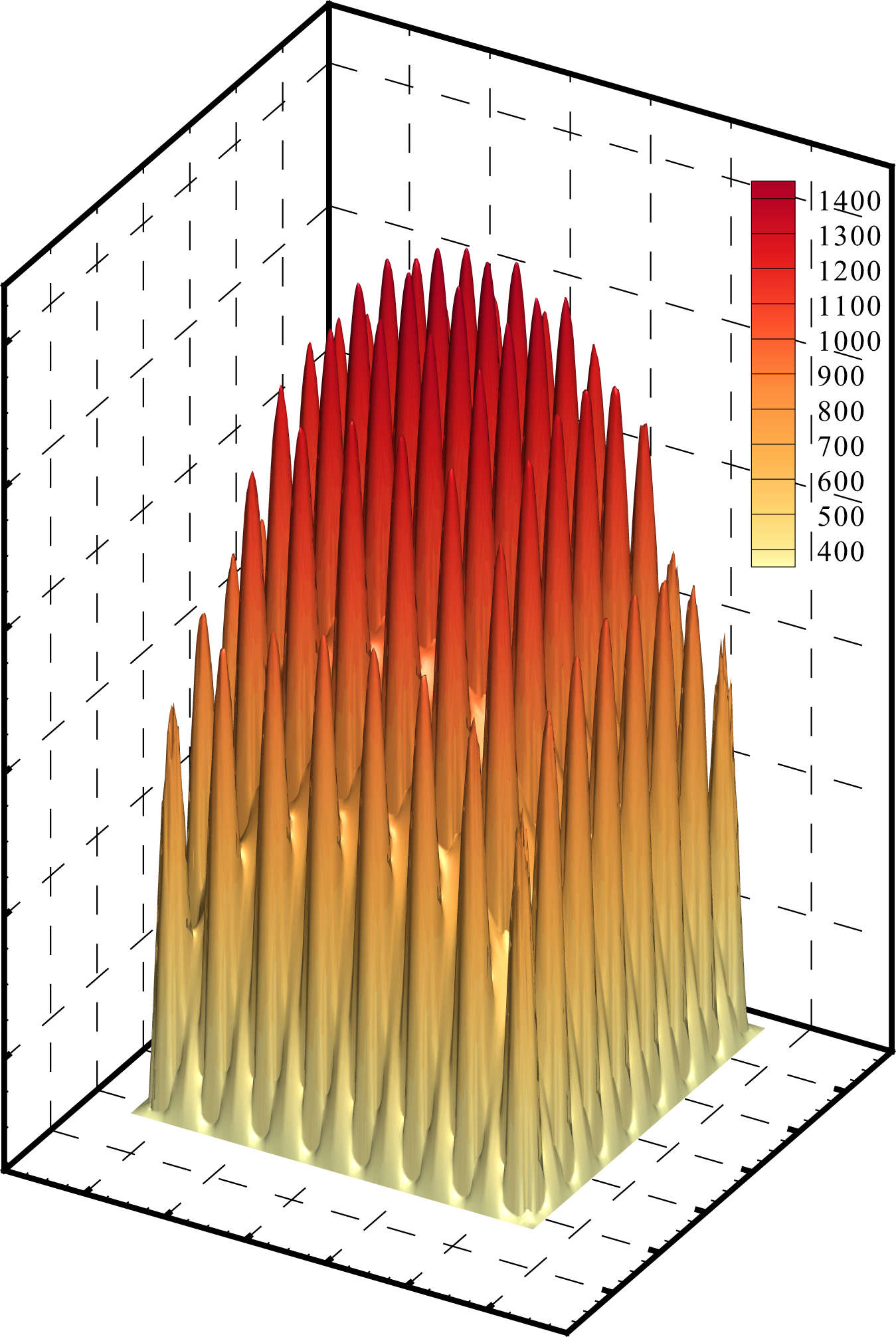}\\(c)
\end{minipage}\hspace{0.01\textwidth}
\begin{minipage}[c]{0.22\textwidth}
  \centering
  \includegraphics[width=\linewidth]{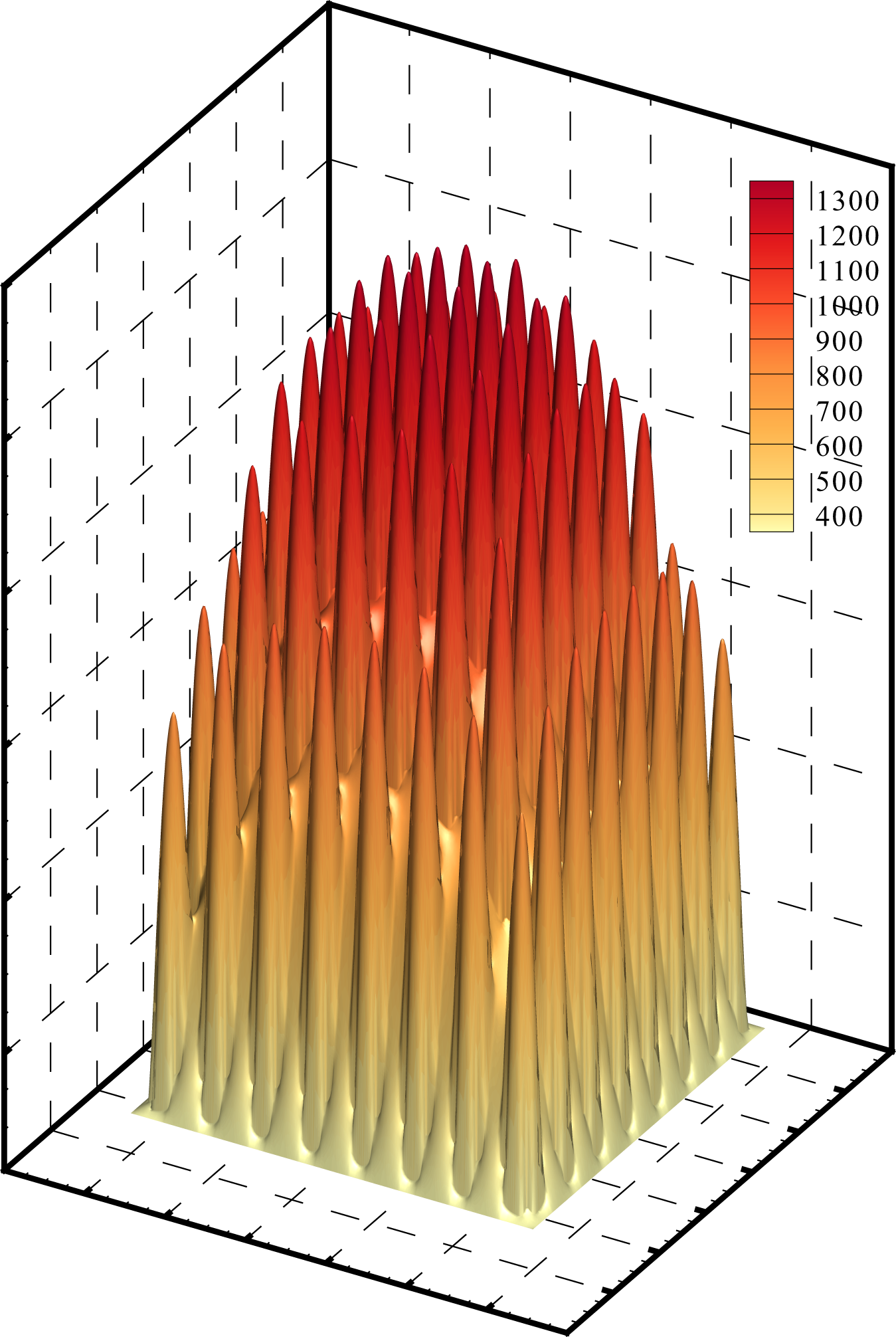}\\(d)
\end{minipage}
\caption{The temperature field at $t=1.0\,\mathrm{s}$: (a) $T_{0}$; (b) $T^{(1{\varepsilon})}$; (c) $T^{(2{\varepsilon})}$; (d) $T^{{\varepsilon}}$.}\label{2dtemperature}
\end{figure}

Ending multi-scale simulation, the resulting temperature fields $T_{0}$, $T^{(1\bm{\varepsilon})}$,
$T^{(2\bm{\varepsilon})}$ and $T^{\bm{\varepsilon}}$ are shown in \autoref{2dtemperature}.

\begin{figure}[!htb]
\centering
\begin{subfigure}[c]{0.45\textwidth}
  \centering
  \includegraphics[width=\linewidth]{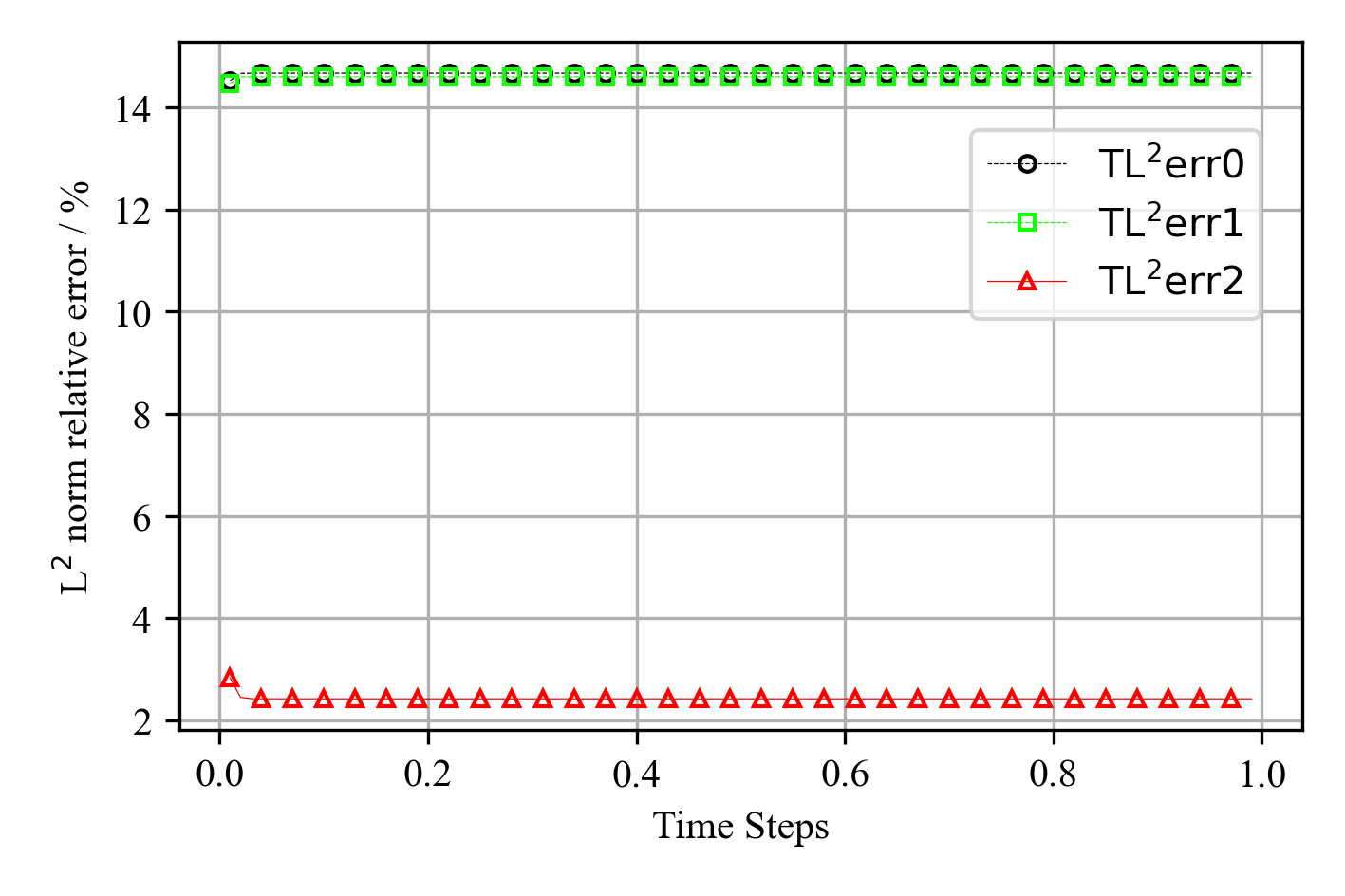}
  \caption{$L^2$ error}
  \label{fig:2dl2-error}
\end{subfigure}
\hfill
\begin{subfigure}[c]{0.45\textwidth}
  \centering
  \includegraphics[width=\linewidth]{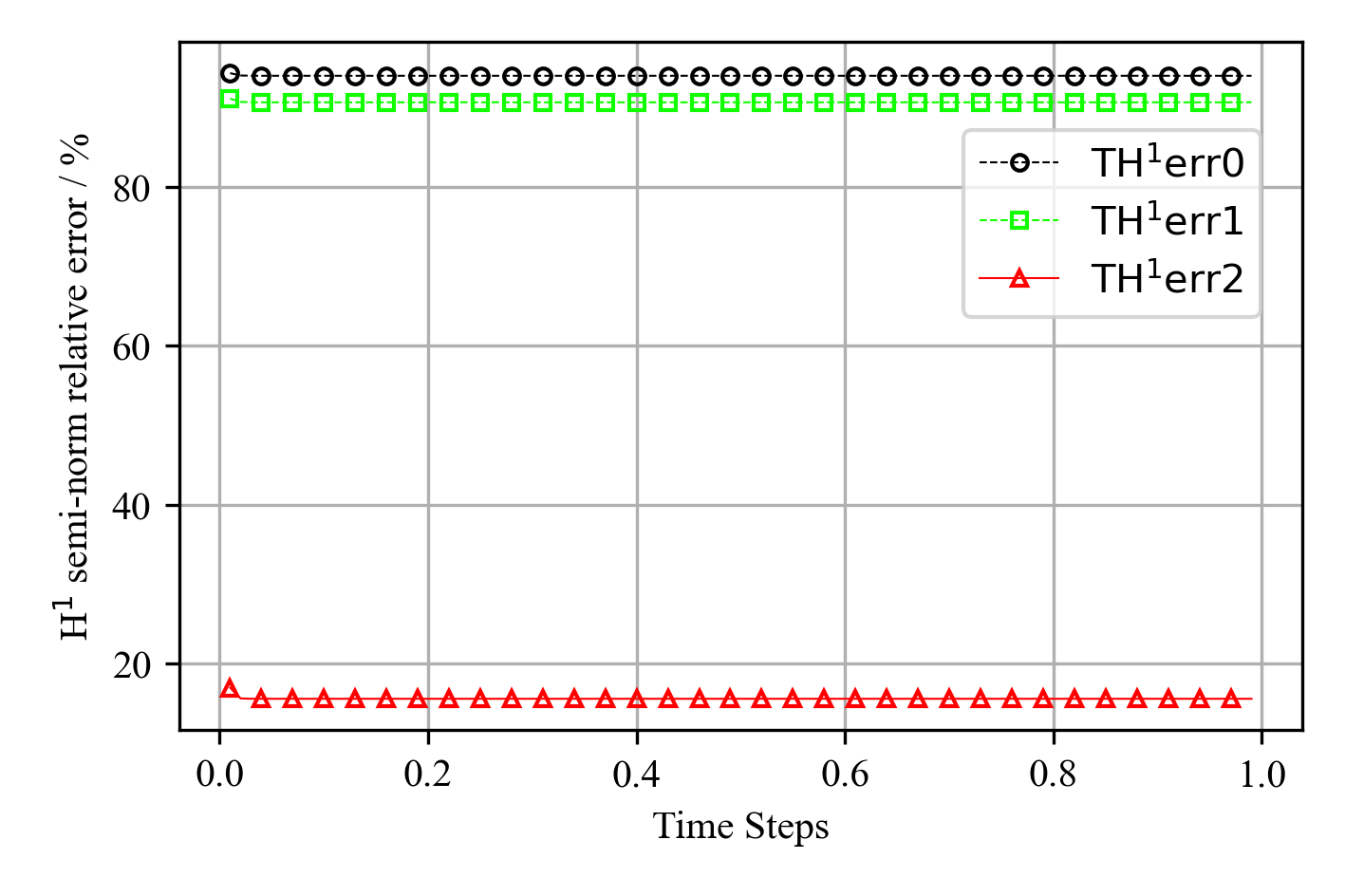}
  \caption{$H^1$ semi-norm error}
  \label{fig:2dh1-error}
\end{subfigure}
\caption{The evolution of relative errors in the $L^2$ norm and $H^1$ semi-norm.}
\label{fig:2derror-comparison}
\end{figure}

Furthermore, the detailed relative errors for the homogenization method, the low-order and high-order multi-scale approaches are presented in \autoref{fig:2derror-comparison}.

The HOMS method achieves $L^2$ norm and $H^1$ semi-norm relative errors of $2.5\%$ and $10\%$, respectively, in contrast to $15\%$ ($L^2$ norm) and $90\%$ ($H^1$ semi-norm) for both the first-order and homogenization methods. This confirms that only the HOMS solutions can accurately resolve the highly microscopic oscillating behavior under temperature-dependent coefficients. The accuracy of homogenized and LOMS solutions is far from enough especially in the $H^1$ semi-norm sense. Moreover, the presented HOMS algorithm is stable and effective after long-time numerical calculation.

\subsection{Example 3: Nonlinear radiative heat transfer problem of 3D composite structure with spherical inclusion microstructure}
\begin{figure}[!htb]
\centering
\begin{minipage}[c]{0.3\textwidth}
  \centering
  \includegraphics[width=0.9\linewidth]{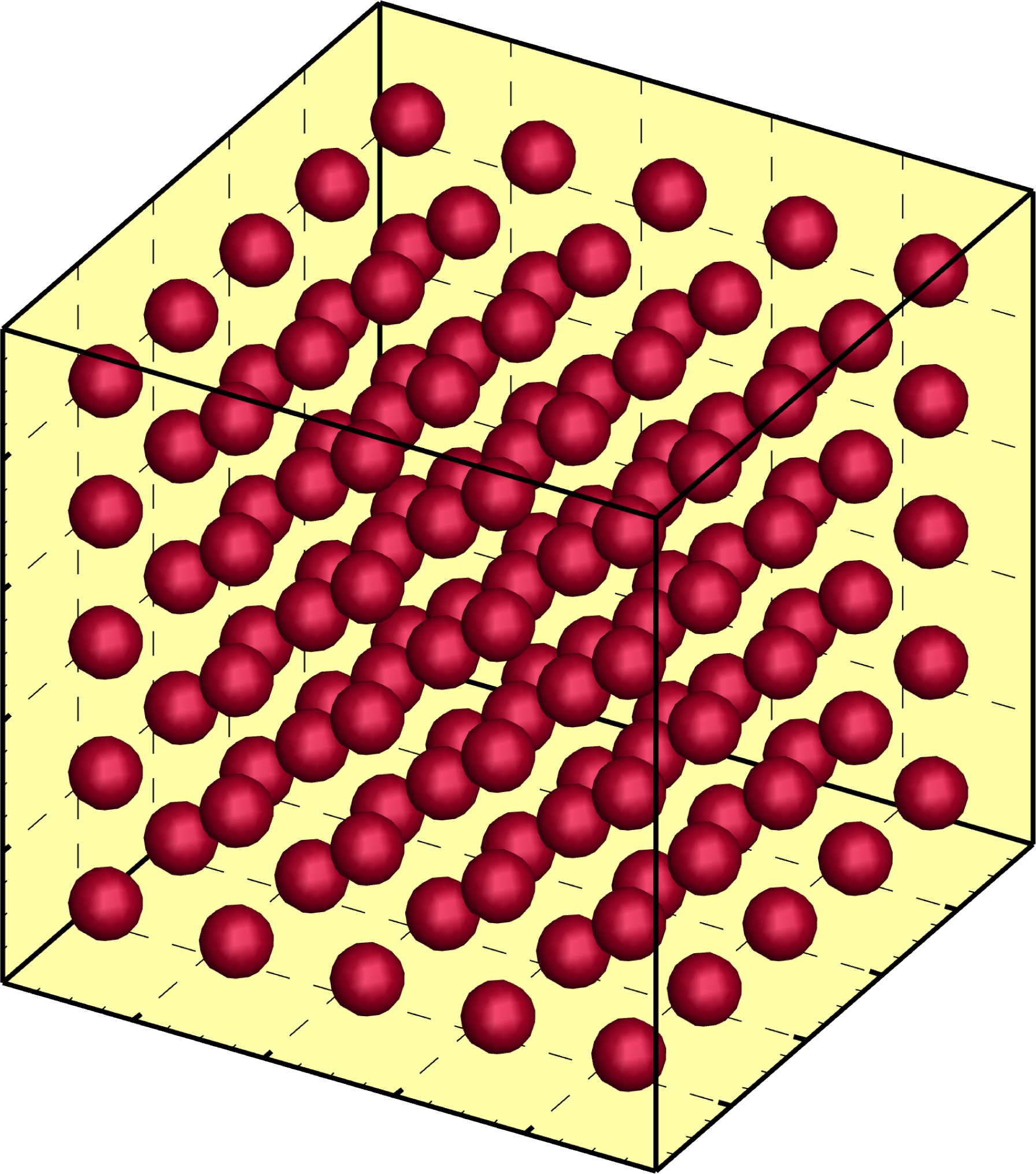}\\
  (a)
\end{minipage}
\begin{minipage}[c]{0.3\textwidth}
  \centering
  \includegraphics[width=0.9\linewidth]{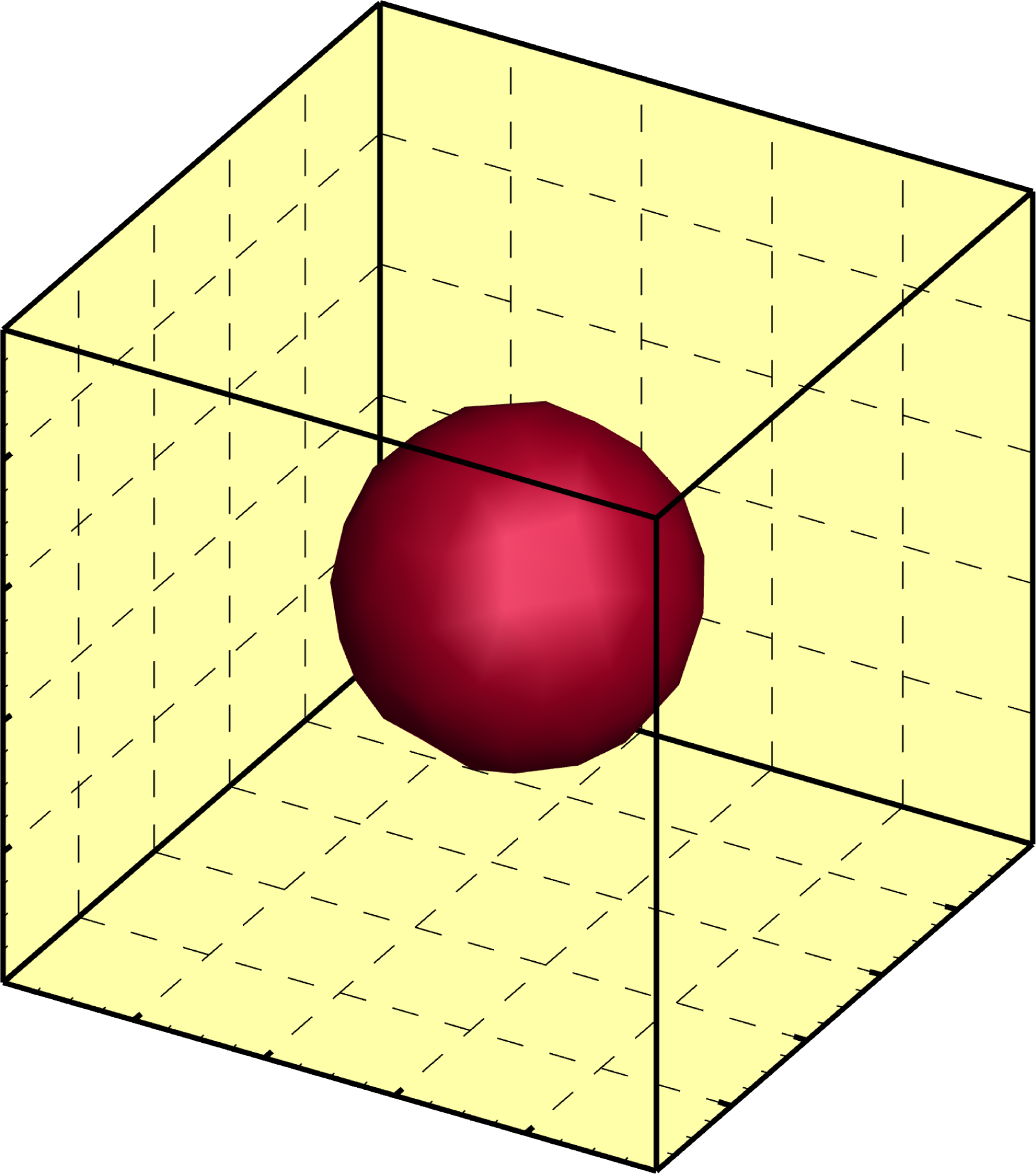}\\
  (b)
\end{minipage}
\begin{minipage}[c]{0.3\textwidth}
  \centering
  \includegraphics[width=0.9\linewidth]{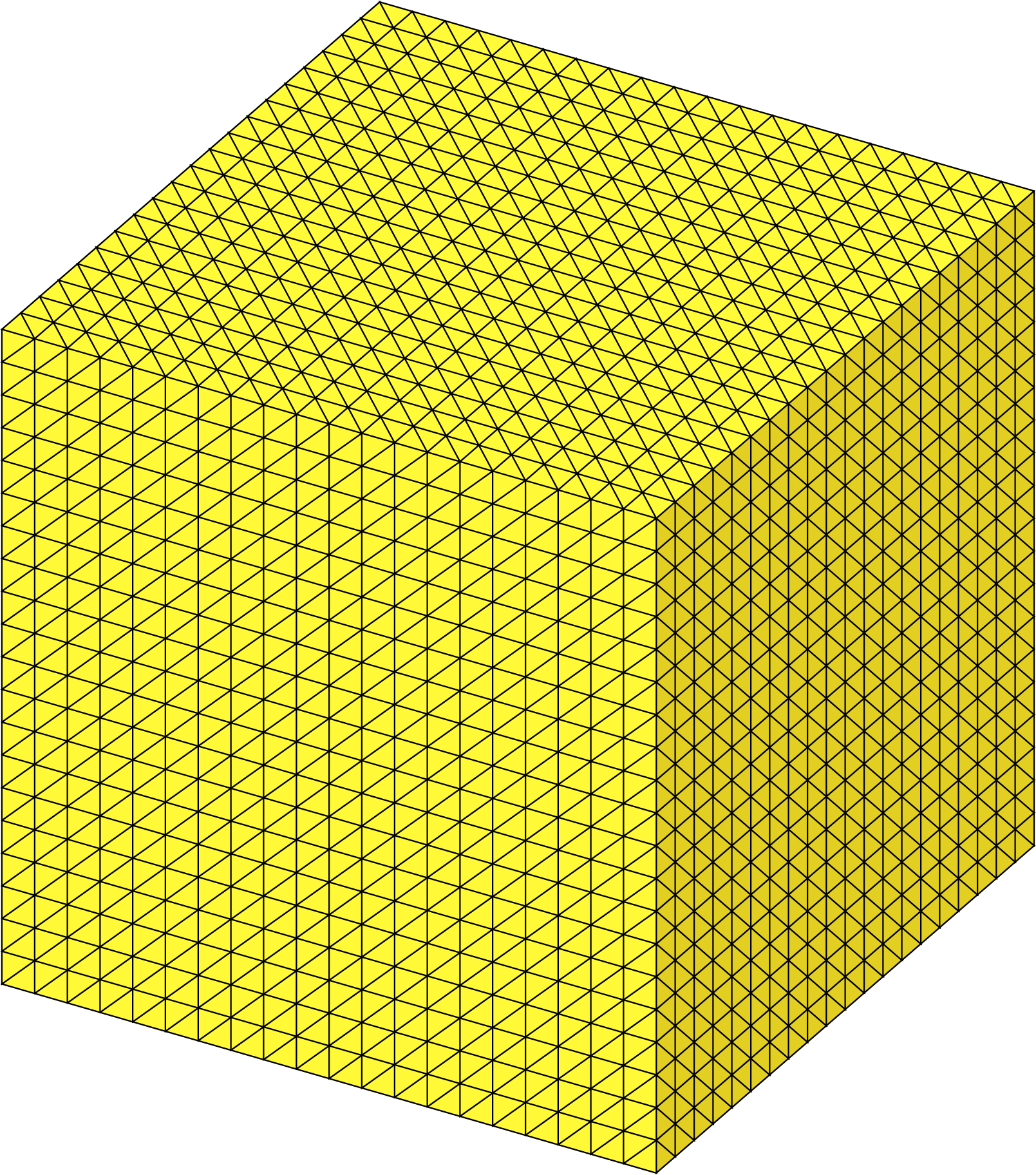}\\
  (c)
\end{minipage}
\caption{(a) The 3D multi-scale structure $\Omega$; (b) the microscopic unit cell $\Omega_{\bm{y}}$; (c) the macroscopic homogenized structure.}\label{ex3dstructure}
\end{figure}

This example studies the nonlinear radiative heat transfer problem of a 3D composite structure with temperature-dependent properties, where spherical inclusions are distributed inside the matrix. The multi-scale structure $\Omega$, macroscopic homogenization structure and microscopic unit cell $\Omega_{\bm{y}}$ are shown in \autoref{ex3dstructure}, where $\Omega = (x_1, x_2, x_3) = [0,1] \times [0,1] \times [0,1] \, \mathrm{cm}^3$, $\Omega_{\bm{y}} = (y_1, y_2, y_3) = [0,1] \times [0,1] \times [0,1]$, and $\varepsilon=1/5$. The material parameters of this example are the same as those in \autoref{tab2dcanshu}.

In addition, the heat source function and initial-boundary conditions in multi-scale problem \eqref{eq:all} of this example are defined as follows:
\[
h(\bm{x},t) = 1.0\times 10^6\,\mathrm{J}/(\mathrm{cm}^3 \cdot \mathrm{s}),\
\hat T (\bm{x},t)= 300.0\,\mathrm{K},\
\widetilde T (\bm{x})= 300.0\,\mathrm{K}.
\]

Now, the tetrahedral mesh generation is implemented to multi-scale problem \eqref{eq:all},
auxiliary cell problems \eqref{eq:Malpha1}, and \eqref{eq:Malpha1alpha2}--\eqref{eq:Nalpha1alpha2}, and
associated homogenized problem \eqref{eq:homogenized}. The computational cost of FEM elements, nodes and time
is illustrated in \autoref{ex3dcomputation}, which clearly indicates that the high-order multi-scale approach
provides a tremendous saving in computing resource, in particular, for 3D composite structure.
\begin{table}[!htb]
\caption{Comparison of the cost on computation resources.}
\label{ex3dcomputation}
\centering
\begin{tabular}{lccc}
\toprule
 & Cell eqs. & Homogenized eqs. & Multi-scale eqs. \\
\midrule
FEM elements & 72700 & 162000 & 1692398 \\
FEM nodes & 12635 & 29791 & 273008 \\
\midrule
 & Off-line stage & On-line stage & FEM \\
\midrule
Computational time & 295.182s & 1440.229s & 5128.691s \\
\bottomrule
\end{tabular}
\end{table}

\begin{figure}[!htb]
\centering
\begin{minipage}[c]{0.22\textwidth}
  \centering
  \includegraphics[width=\linewidth]{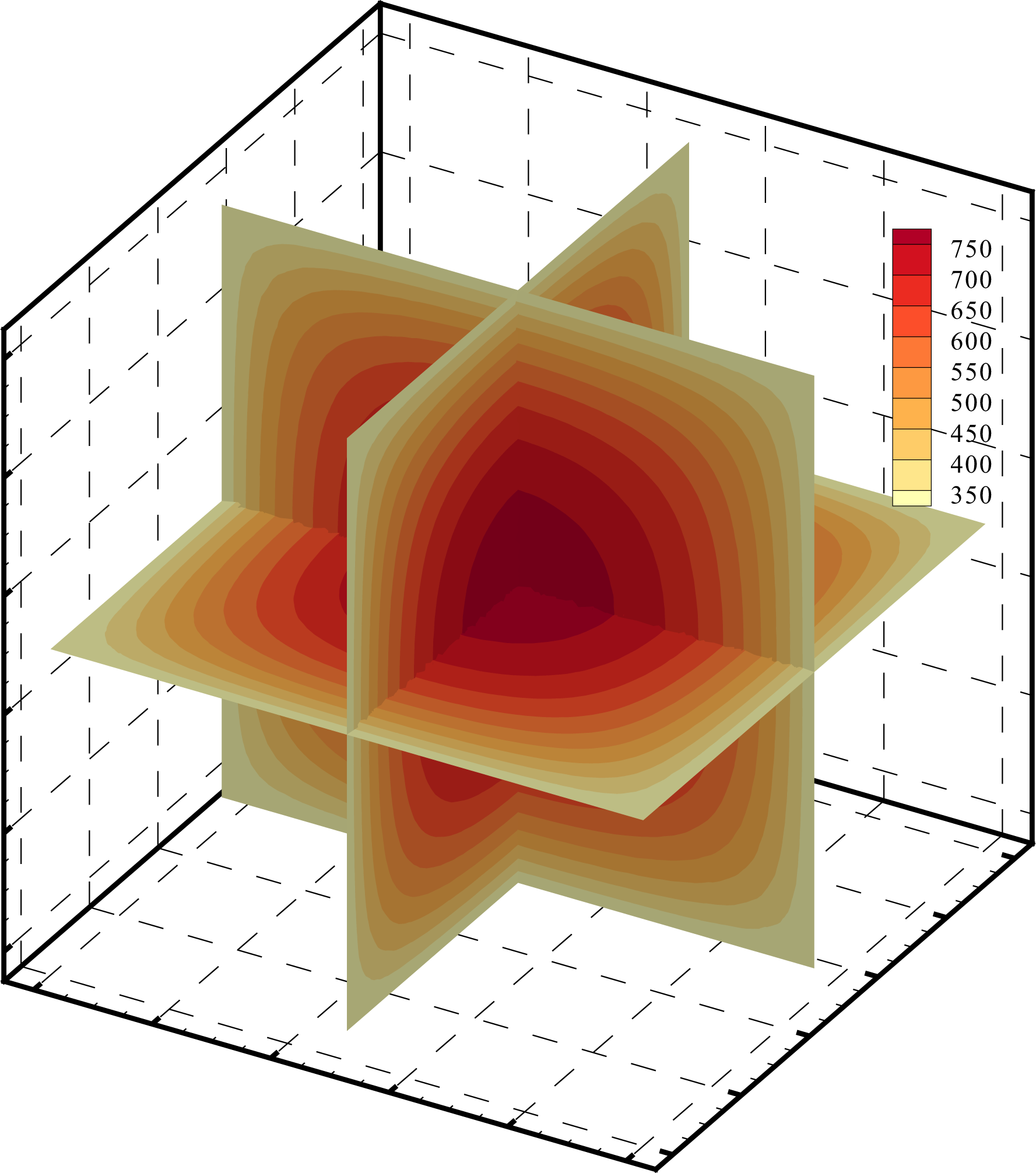}\\(a)
\end{minipage}\hspace{0.01\textwidth}
\begin{minipage}[c]{0.22\textwidth}
  \centering
  \includegraphics[width=\linewidth]{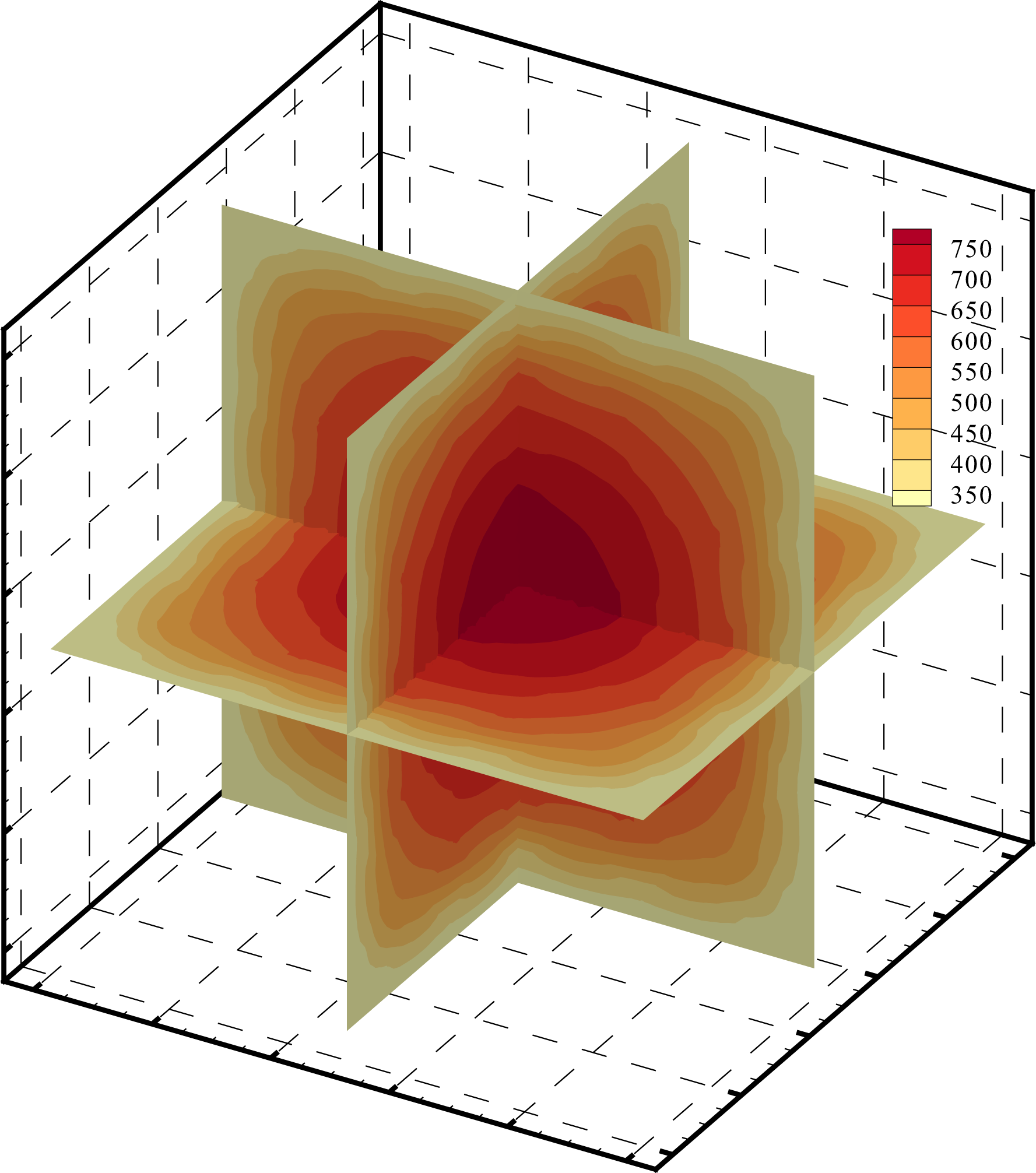}\\(b)
\end{minipage}\hspace{0.01\textwidth}
\begin{minipage}[c]{0.22\textwidth}
  \centering
  \includegraphics[width=\linewidth]{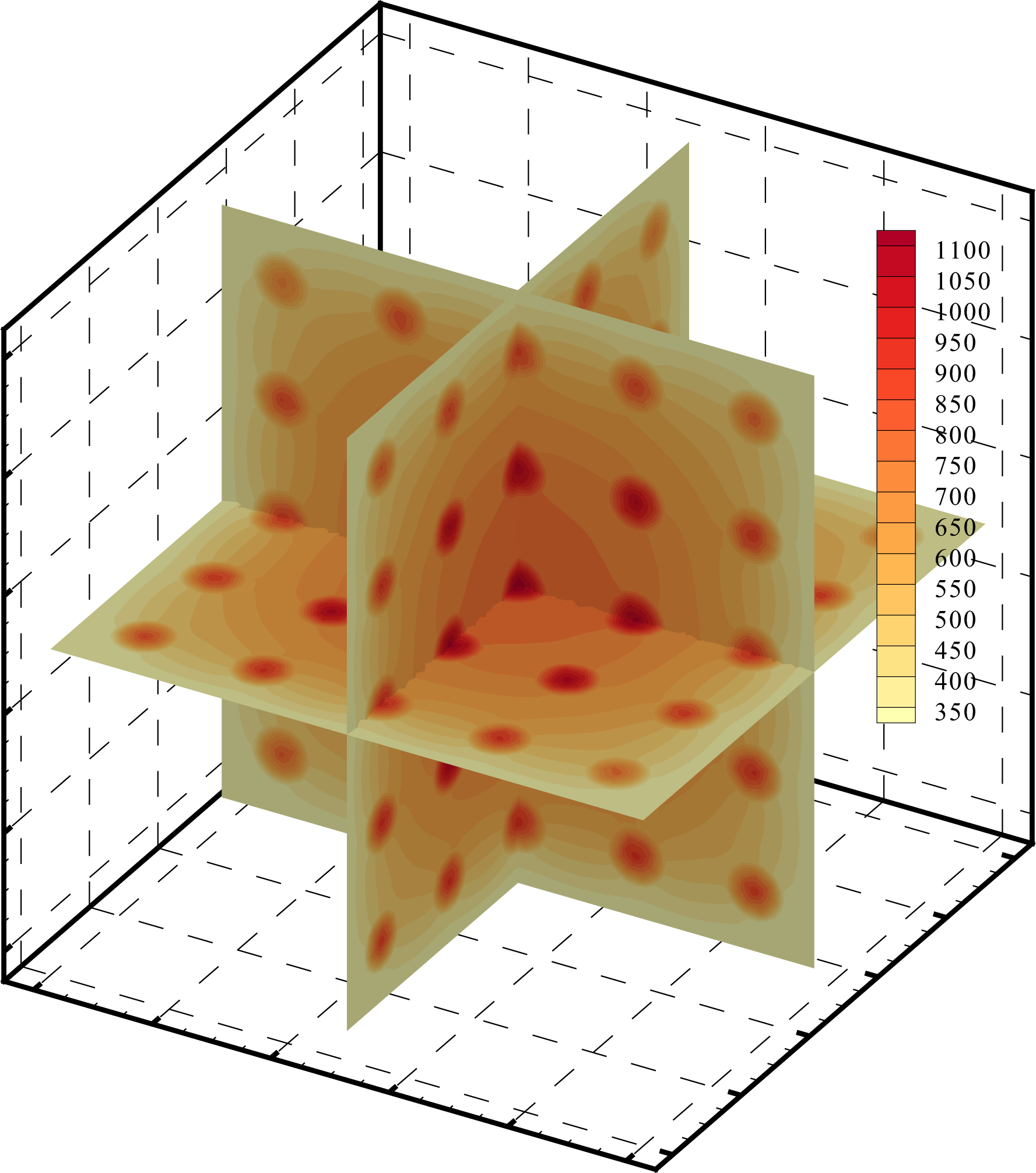}\\(c)
\end{minipage}\hspace{0.01\textwidth}
\begin{minipage}[c]{0.22\textwidth}
  \centering
  \includegraphics[width=\linewidth]{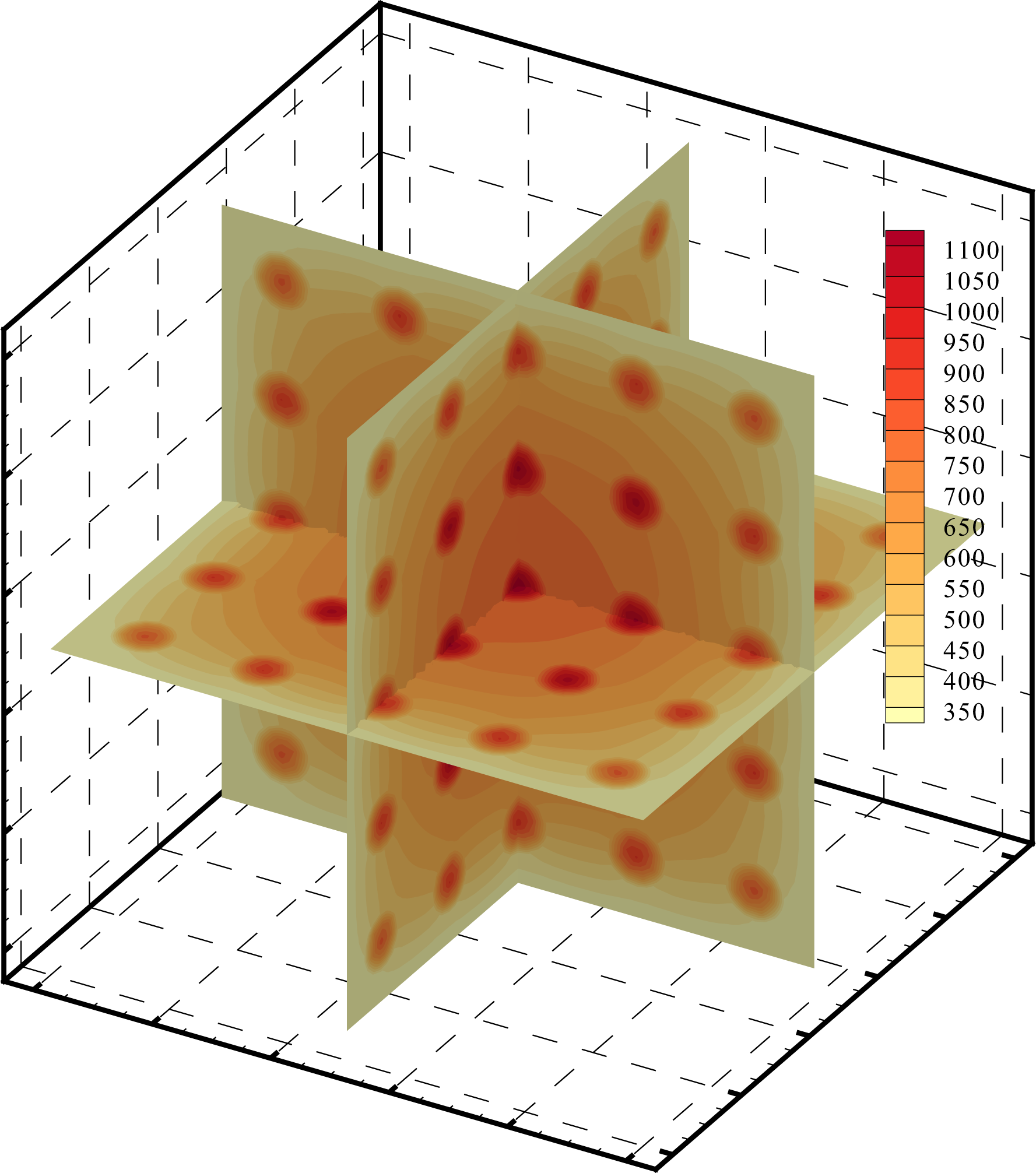}\\(d)
\end{minipage}
\caption{The temperature field slices at $t=1.0\,\mathrm{s}$: (a) $T_{0}$; (b) $T^{(1{\varepsilon})}$;
(c) $T^{(2{\varepsilon})}$; (d) $T^{{\varepsilon}}$.}\label{ex33dtemperatureSlice}
\end{figure}

\begin{figure}[!htb]
\centering
\begin{minipage}[c]{0.22\textwidth}
  \centering
  \includegraphics[width=\linewidth]{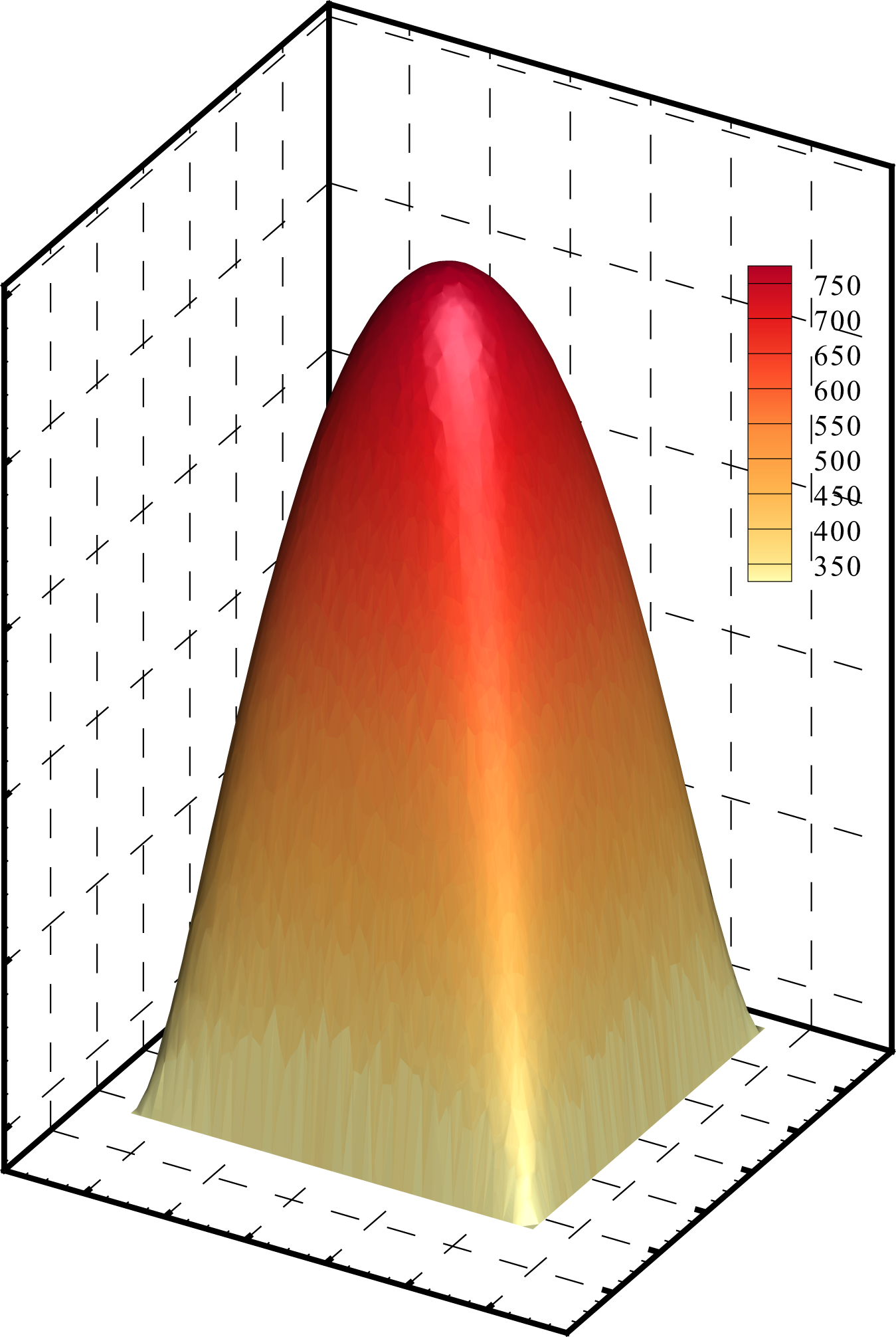}\\(a)
\end{minipage}\hspace{0.01\textwidth}
\begin{minipage}[c]{0.22\textwidth}
  \centering
  \includegraphics[width=\linewidth]{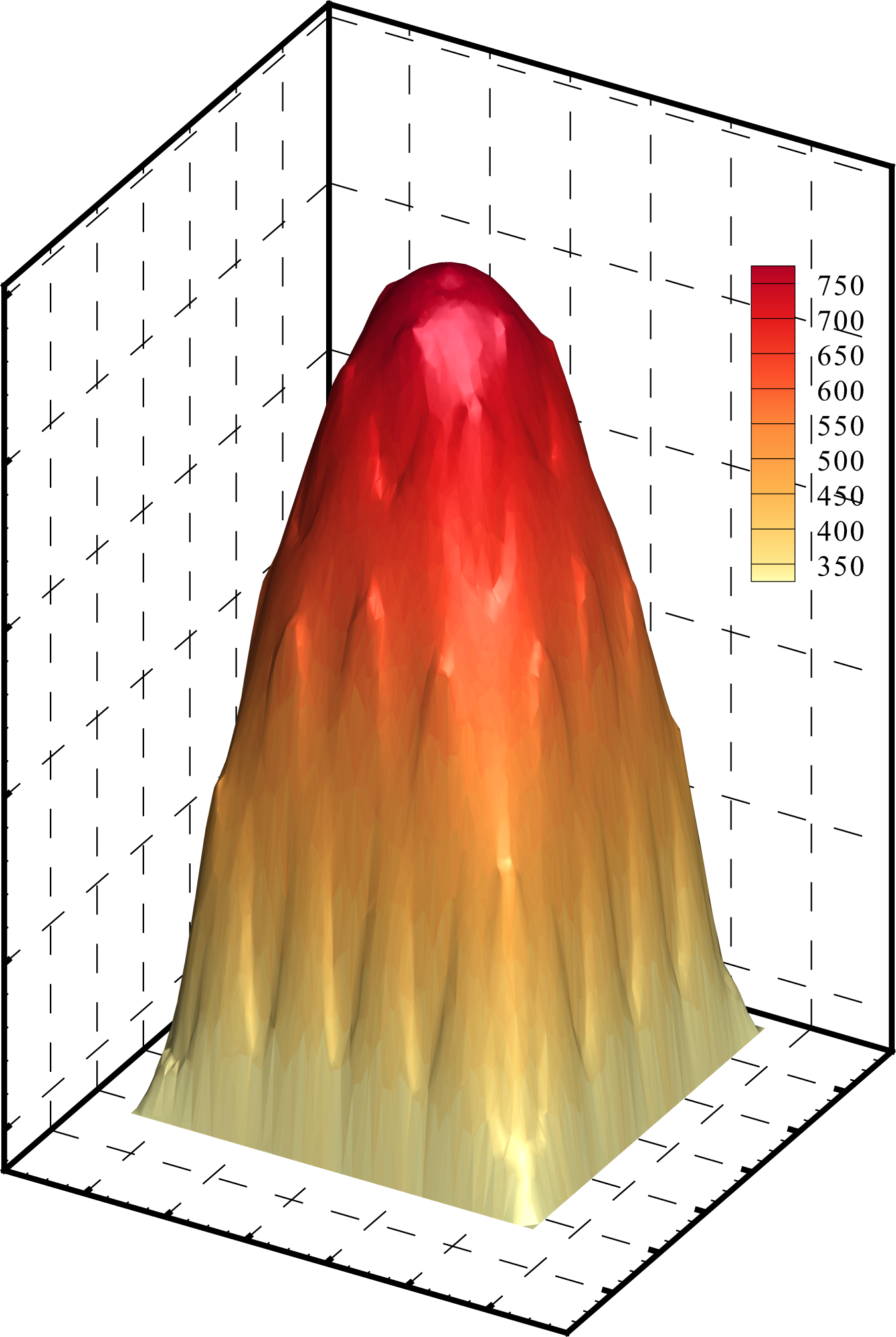}\\(b)
\end{minipage}\hspace{0.01\textwidth}
\begin{minipage}[c]{0.22\textwidth}
  \centering
  \includegraphics[width=\linewidth]{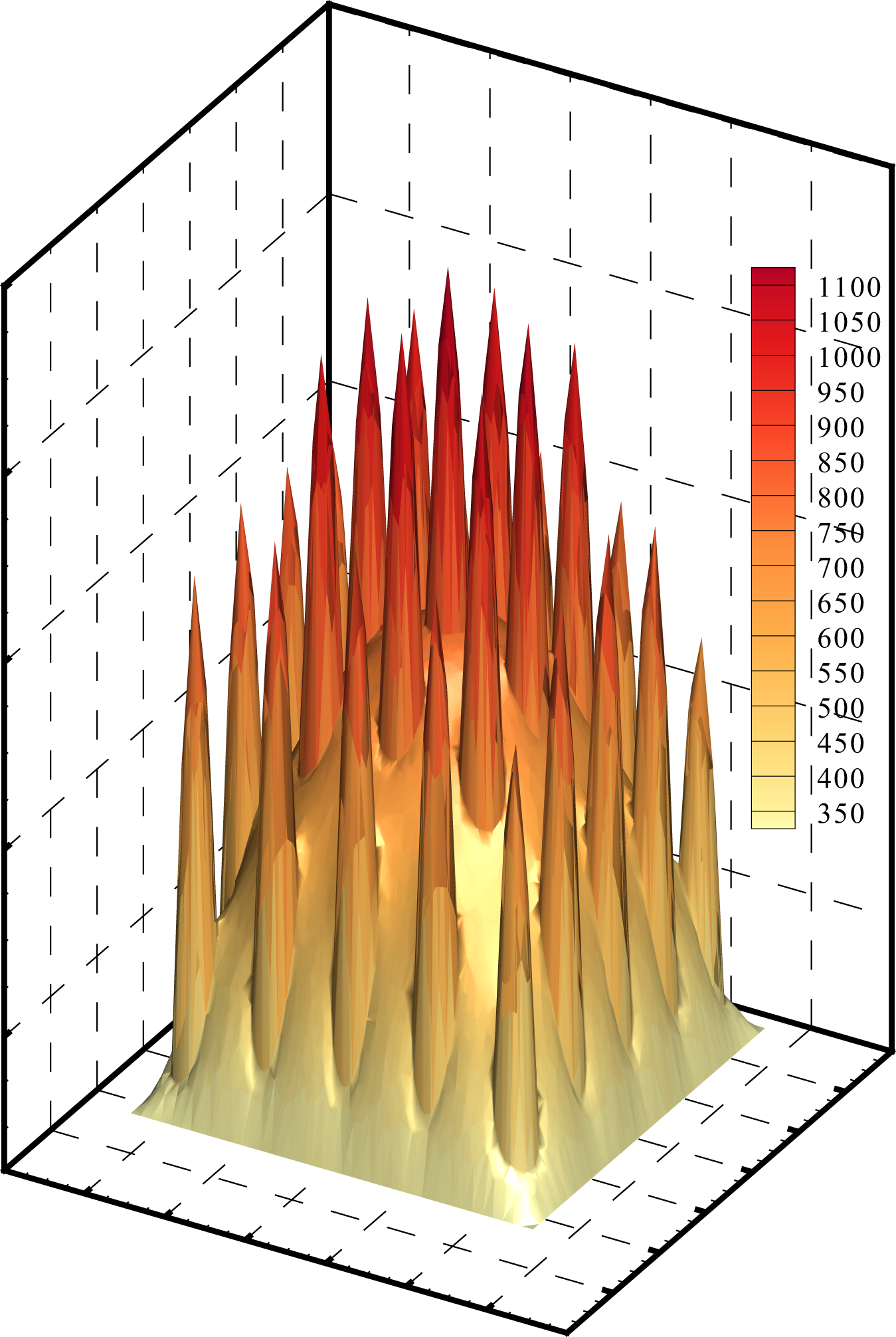}\\(c)
\end{minipage}\hspace{0.01\textwidth}
\begin{minipage}[c]{0.22\textwidth}
  \centering
  \includegraphics[width=\linewidth]{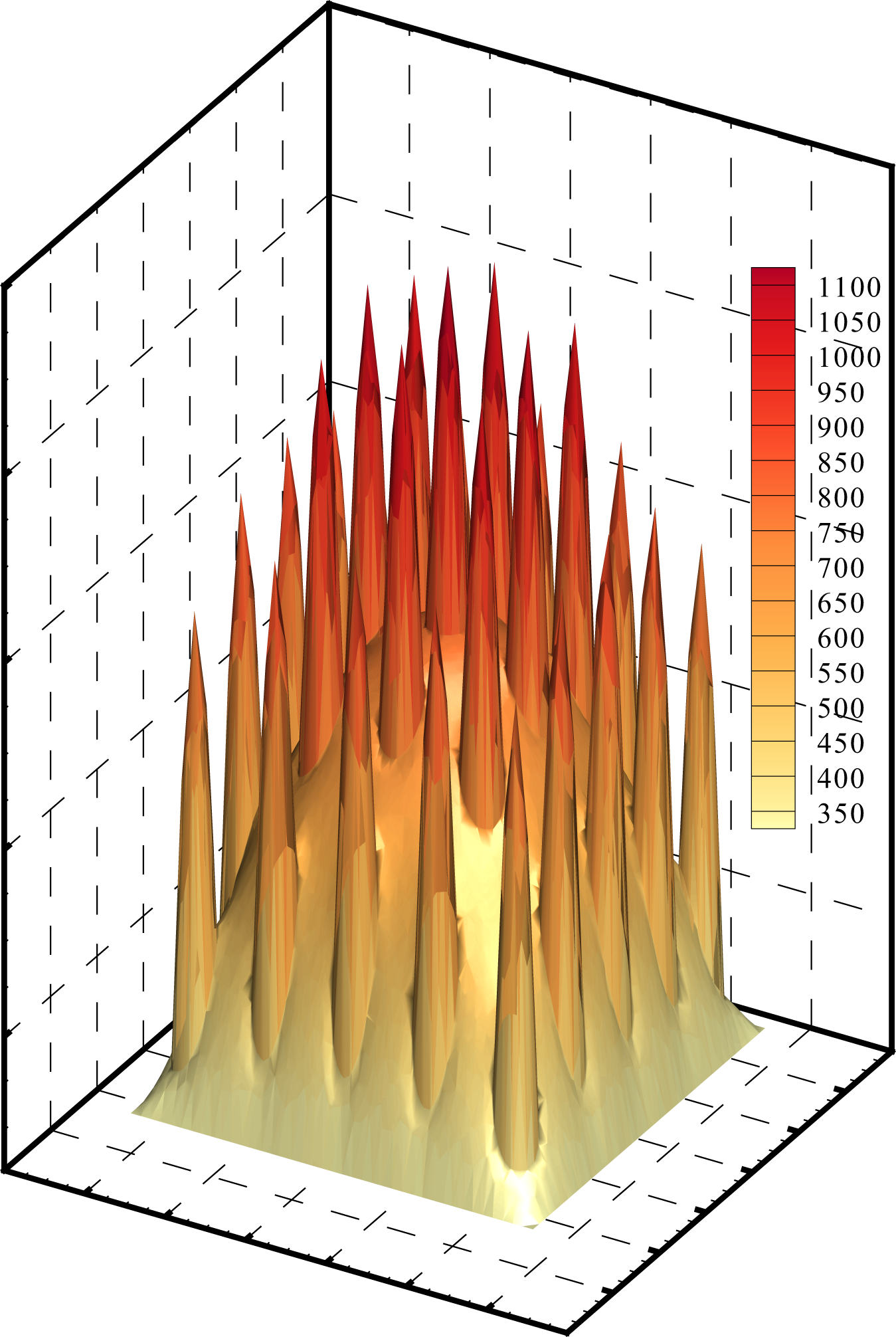}\\(d)
\end{minipage}
\caption{The temperature field on the cross-section $x_3=0.5$\,cm at $t=1.0\,\mathrm{s}$: (a) $T_{0}$; (b) $T^{(1{\varepsilon})}$; (c) $T^{(2{\varepsilon})}$; (d) $T^{{\varepsilon}}$.}\label{ex33dtemperaturez05}
\end{figure}

After multi-scale simulation, \autoref{ex33dtemperatureSlice} displays the temperature distributions on
the slices $x_1=0.5$, $x_2=0.5$, and $x_3=0.5$\,cm at $t=1.0\,\mathrm{s}$, and \autoref{ex33dtemperaturez05} shows the temperature field on the cross-section $x_3=0.5$\,cm, at $t=1.0\,\mathrm{s}$.

\begin{figure}[!htb]
\centering
\begin{subfigure}[c]{0.45\textwidth}
  \centering
  \includegraphics[width=\linewidth]{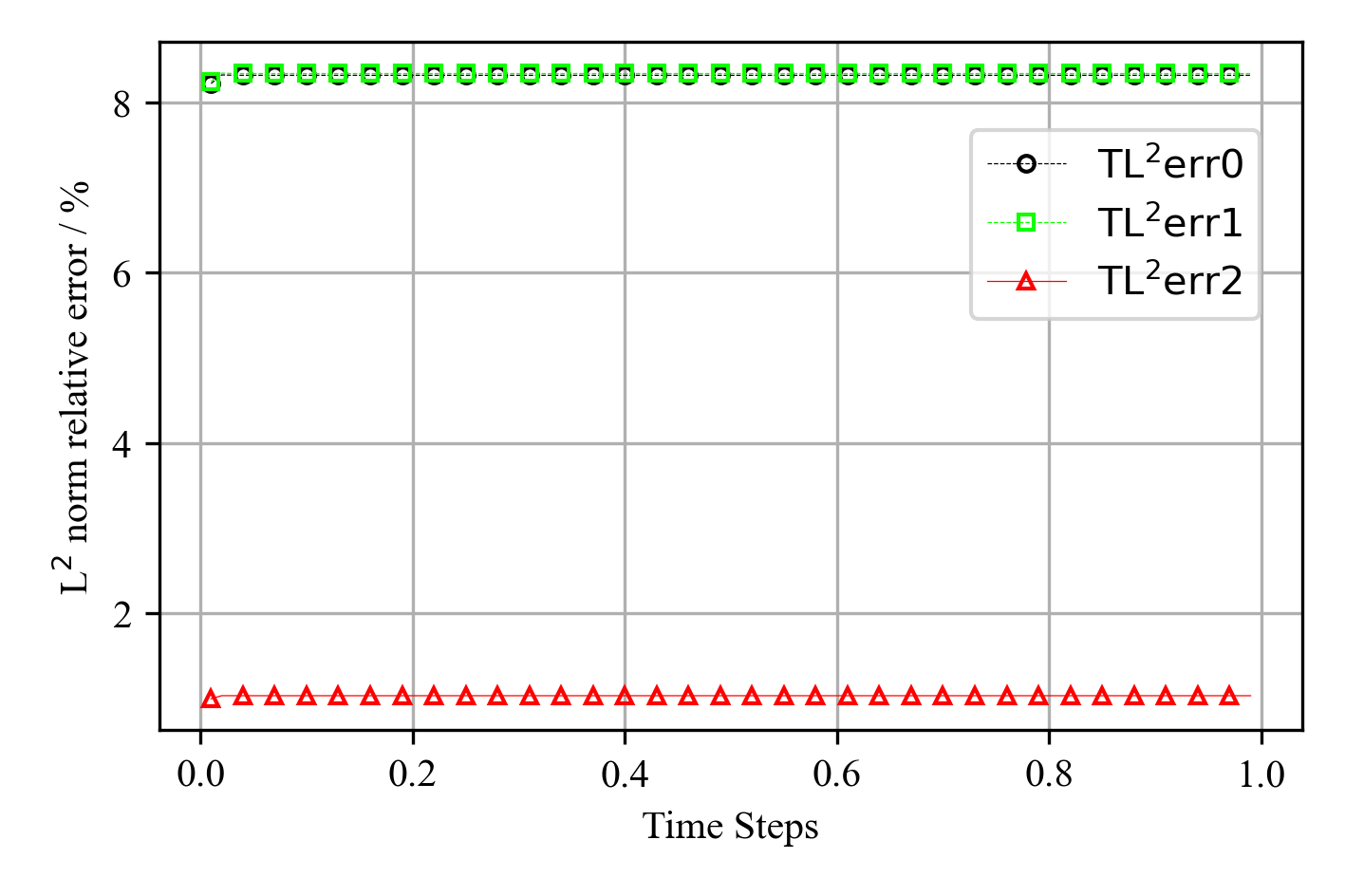}
  \caption{$L^2$ error}
  \label{fig:l2-error}
\end{subfigure}
\hfill
\begin{subfigure}[c]{0.45\textwidth}
  \centering
  \includegraphics[width=\linewidth]{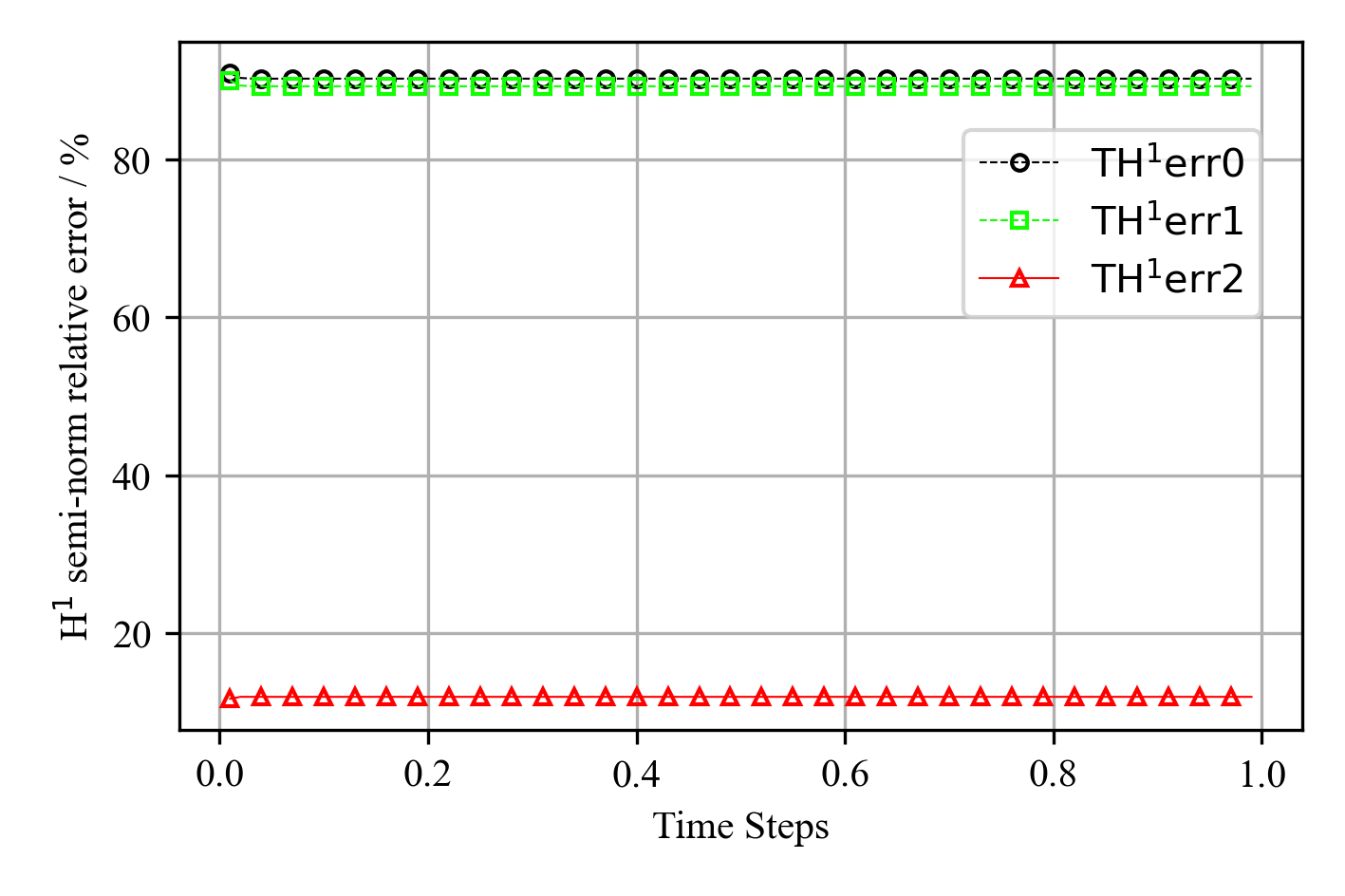}
  \caption{$H^1$ semi-norm error}
  \label{fig:h1-error}
\end{subfigure}
\caption{The evolution of relative errors in the $L^2$ norm and $H^1$ semi-norm.}
\label{fig:ex33derror-comparison}
\end{figure}

Additionally, the specific relative numerical errors for the homogenization method, the low-order and high-order multi-scale approaches are presented in \autoref{fig:ex33derror-comparison}.

According to the numerical results in \autoref{ex33dtemperatureSlice} and \autoref{ex33dtemperaturez05}, it can be clearly seen that only the HOMS solution captures the high-frequency oscillations of the 3D composite structure, while the low-order and homogenized solutions remain excessively smooth. The error comparison in \autoref{fig:ex33derror-comparison} shows that the HOMS method achieves $L^2$ norm and $H^1$ semi-norm relative errors of $2\%$ and $17\%$, respectively, whereas both the low-order and homogenization methods yield more than $10\%$ ($L^2$ norm) and $90\%$ ($H^1$ semi-norm). Thus, the HOMS methodology is suitable for simulating the highly microscopic oscillating responses of this 3D composite structure with a vast number of unit cells.

\subsection{Example 4: Nonlinear radiative heat transfer problem of 3D composite structure with complex inclusion microstructure}
\begin{figure}[!htb]
\centering
\begin{minipage}[c]{0.3\textwidth}
  \centering
  \includegraphics[width=0.9\linewidth]{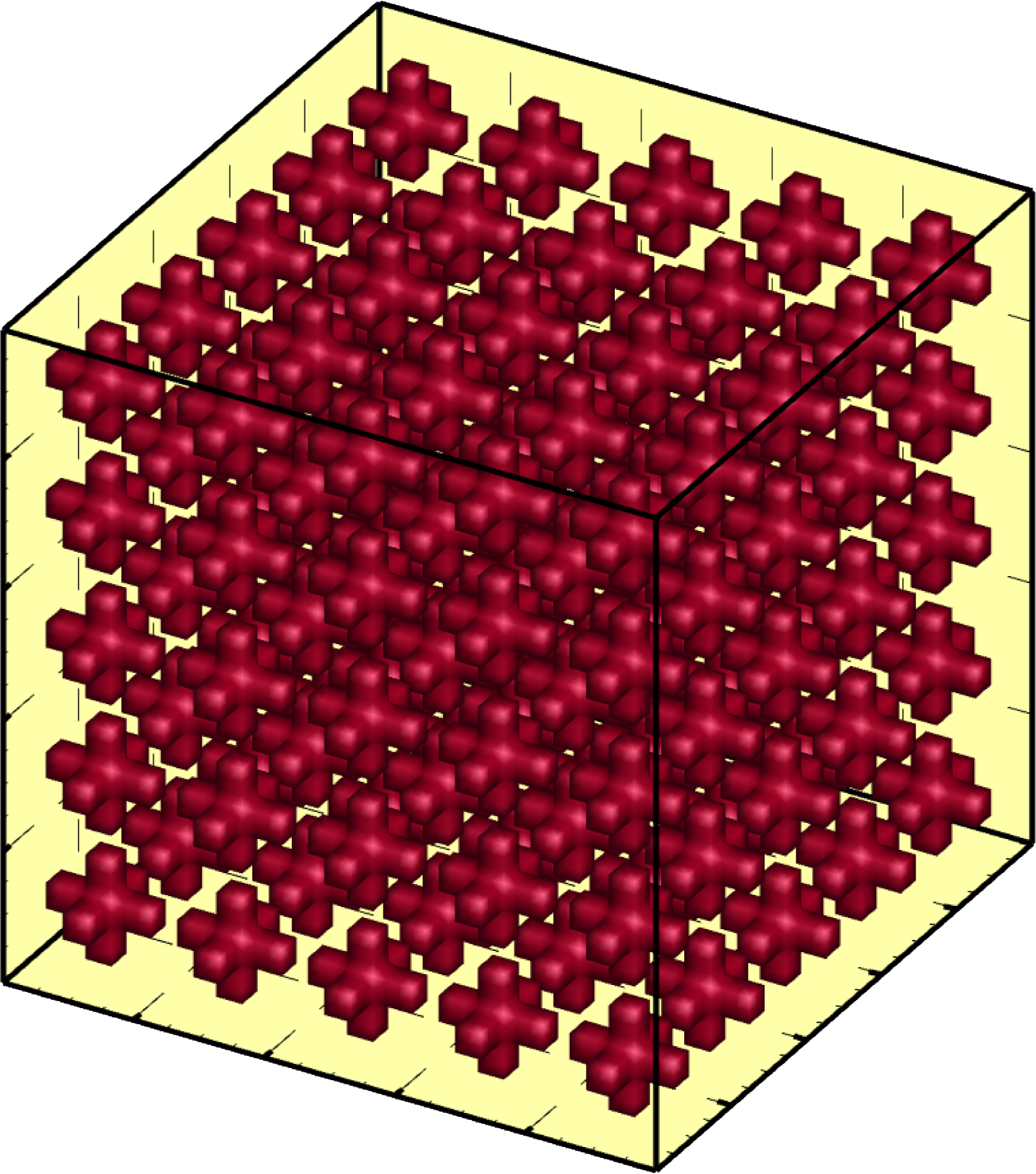}\\(a)
\end{minipage}
\begin{minipage}[c]{0.3\textwidth}
  \centering
  \includegraphics[width=0.9\linewidth]{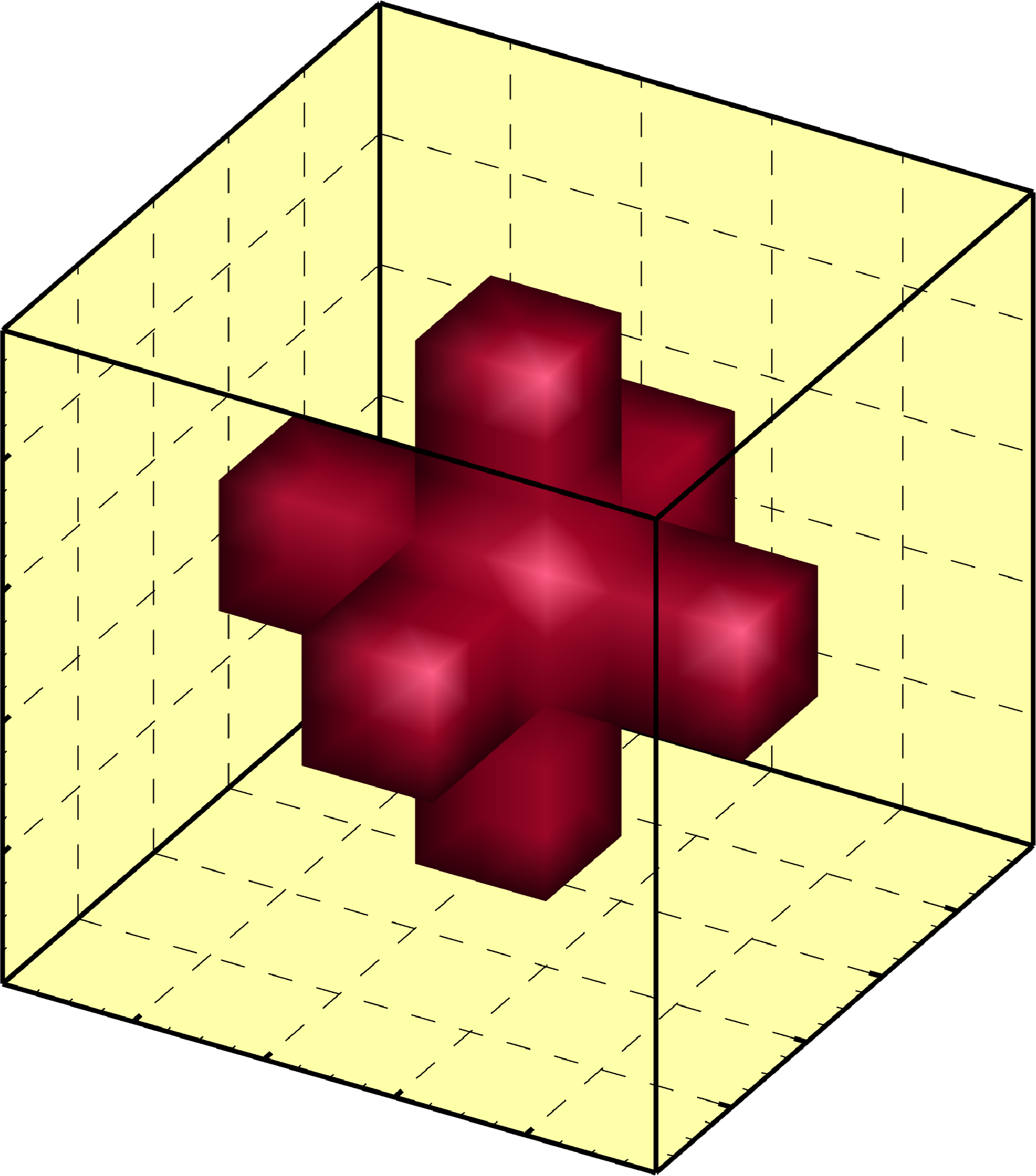}\\(b)
\end{minipage}
\begin{minipage}[c]{0.3\textwidth}
  \centering
  \includegraphics[width=0.9\linewidth]{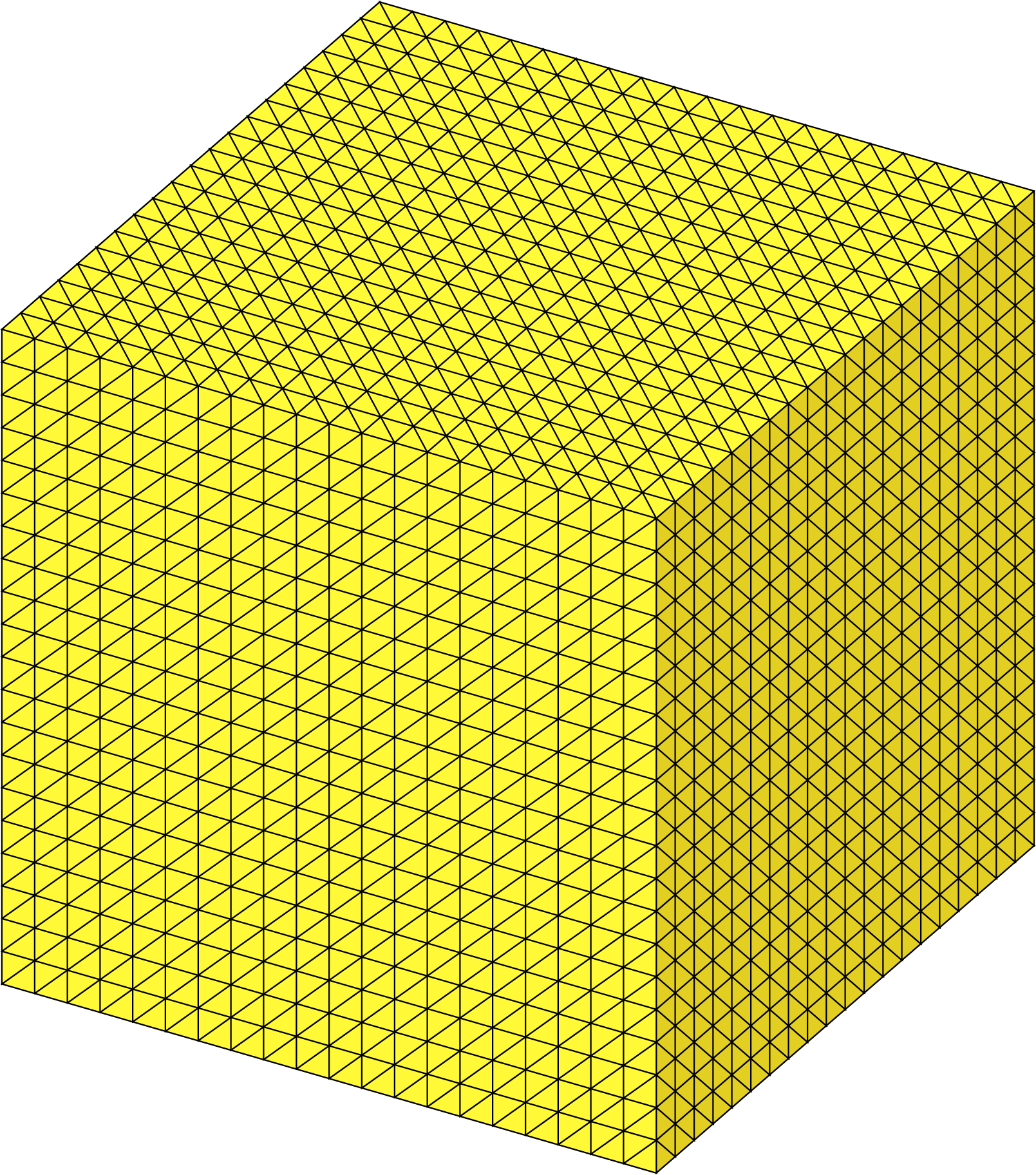}\\(c)
\end{minipage}
\caption{(a) The 3D multi-scale structure $\Omega$; (b) the microscopic unit cell $\Omega_{\bm{y}}$; (c) the macroscopic homogenized structure.}
\label{ex43dstructure}
\end{figure}

This example considers an alternative 3D composite structure with a complex inclusion microstructure. The multi-scale structure $\Omega$, macroscopic homogenization structure and microscopic unit cell $\Omega_{\bm{y}}$ are shown
in \autoref{ex43dstructure}, where $\Omega = (x_1, x_2, x_3) = [0,1] \times [0,1] \times [0,1] \, \mathrm{cm}^3$,
$\Omega_{\bm{y}} = (y_1, y_2, y_3) = [0,1] \times [0,1] \times [0,1]$, and $\varepsilon=1/5$. This example shares the same material parameters as \autoref{tab2dcanshu} and the same heat source function and initial-boundary conditions as the previous example.

After that, the tetrahedral mesh generation is implemented to multi-scale problem, auxiliary cell problems, and
corresponding homogenized problem. The computational cost is illustrated in \autoref{ex43dcomputation}, which obviously demonstrates that the HOMS approach introduced herein possesses a substantial reduction in computational resource utilization, particularly in terms of CPU memory and time, when contrasted with high-resolution DNS.
\begin{table}[!htb]
\caption{Comparison of the cost on computation resources.}
\label{ex43dcomputation}
\centering
\begin{tabular}{lccc}
\toprule
 & Cell eqs. & Homogenized eqs. & Multi-scale eqs. \\
\midrule
 FEM elements & 6000 & 162000 & 750000 \\
 FEM nodes & 1331 & 29791 & 132651 \\
\midrule
 & Off-line stage & On-line stage & FEM \\
\midrule
Computational time & 141.623s & 2091.967s & 3378.192s \\
\bottomrule
\end{tabular}
\end{table}

\begin{figure}[!htb]
\centering
\begin{minipage}[c]{0.24\textwidth}
  \centering
  \includegraphics[width=\linewidth]{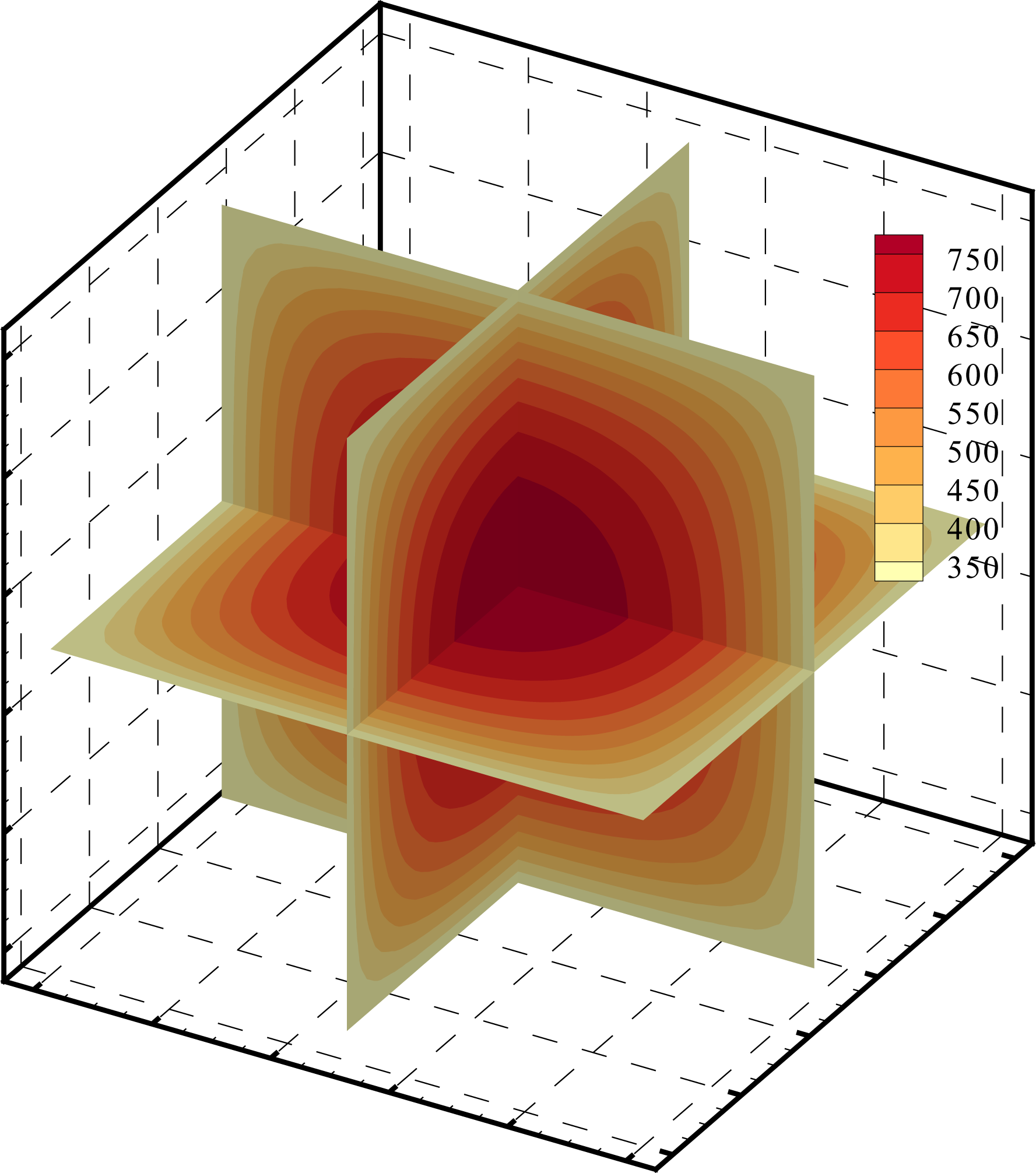}\\(a)
\end{minipage}
\begin{minipage}[c]{0.24\textwidth}
  \centering
  \includegraphics[width=\linewidth]{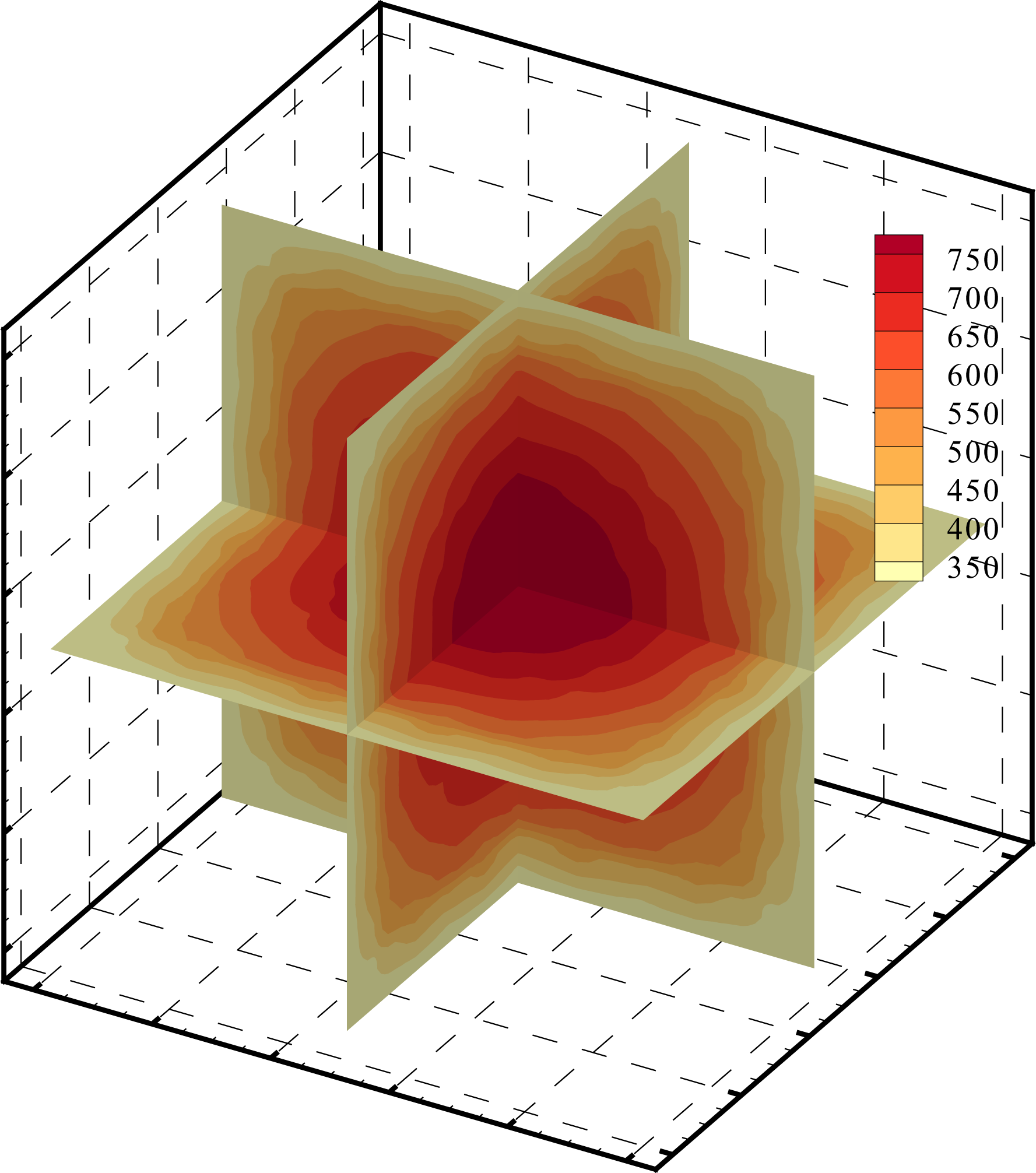}\\(b)
\end{minipage}
\begin{minipage}[c]{0.24\textwidth}
  \centering
  \includegraphics[width=\linewidth]{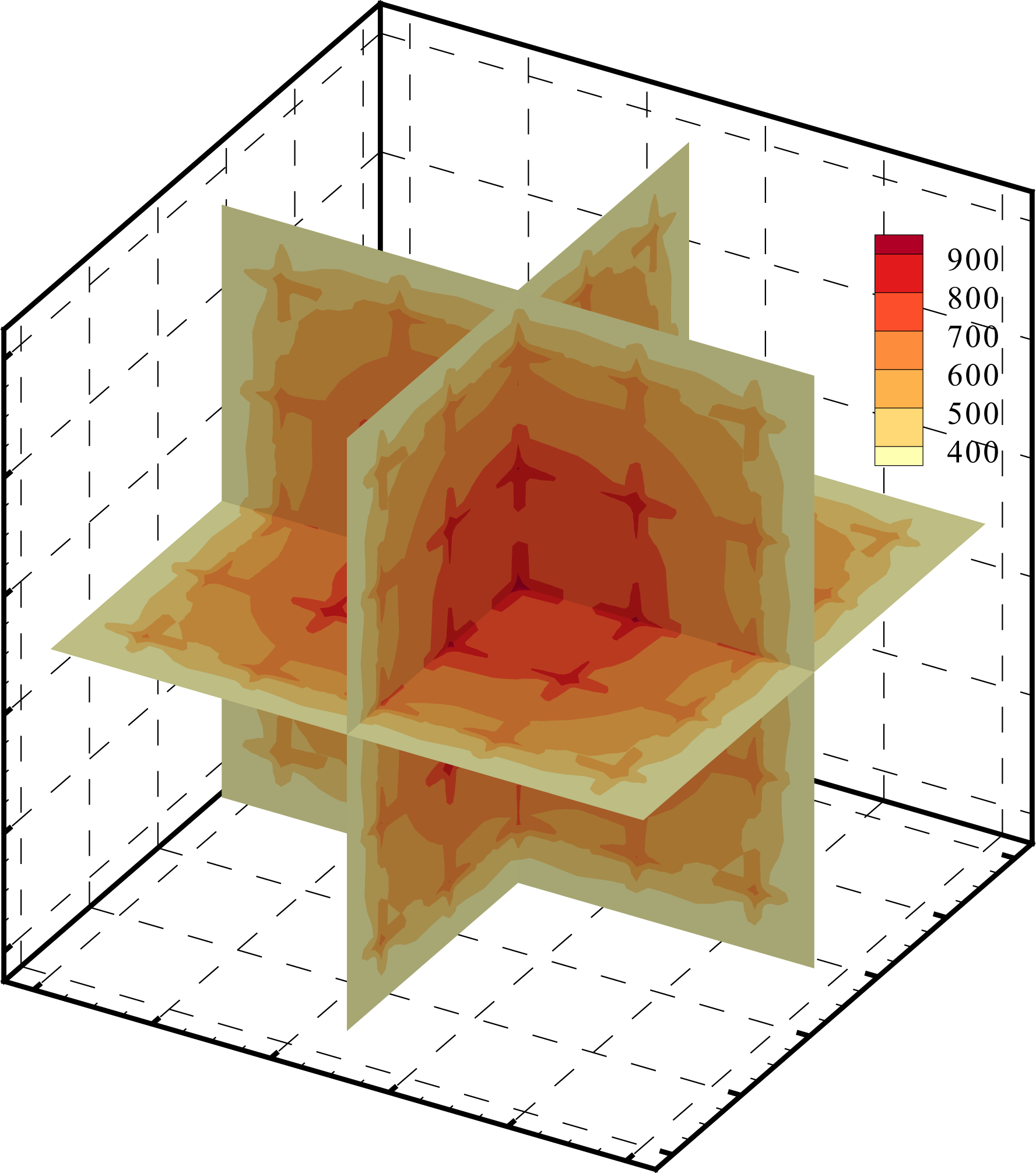}\\(c)
\end{minipage}
\begin{minipage}[c]{0.24\textwidth}
  \centering
  \includegraphics[width=\linewidth]{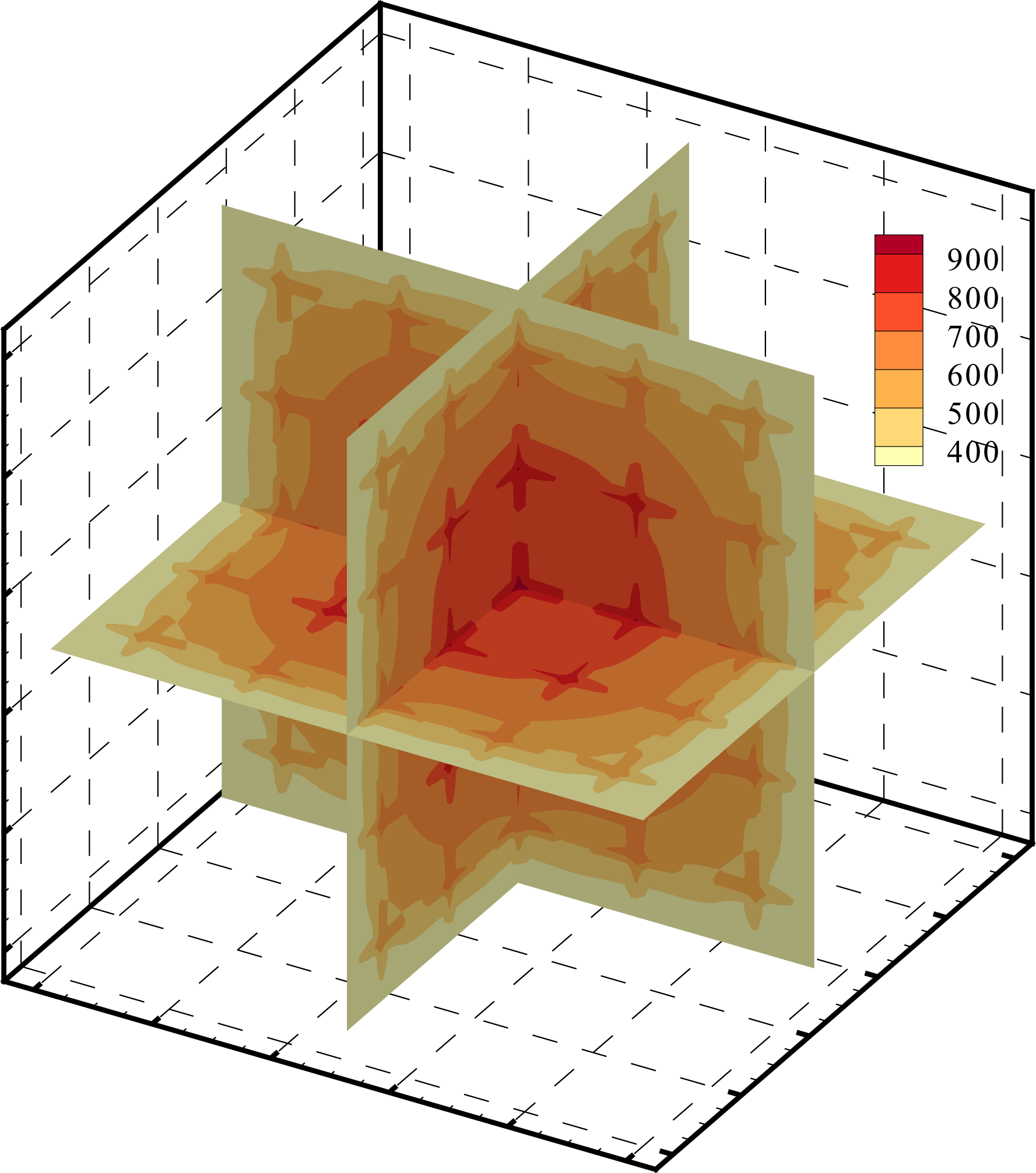}\\(d)
\end{minipage}
\caption{The temperature field slices at $t=1.0\,\mathrm{s}$: (a) $T_{0}$; (b) $T^{(1{\varepsilon})}$;
(c) $T^{(2{\varepsilon})}$; (d) $T^{{\varepsilon}}$.}\label{ex43dtemperatureSlice}
\end{figure}

\begin{figure}[!htb]
\centering
\begin{minipage}[c]{0.22\textwidth}
  \centering
  \includegraphics[width=\linewidth]{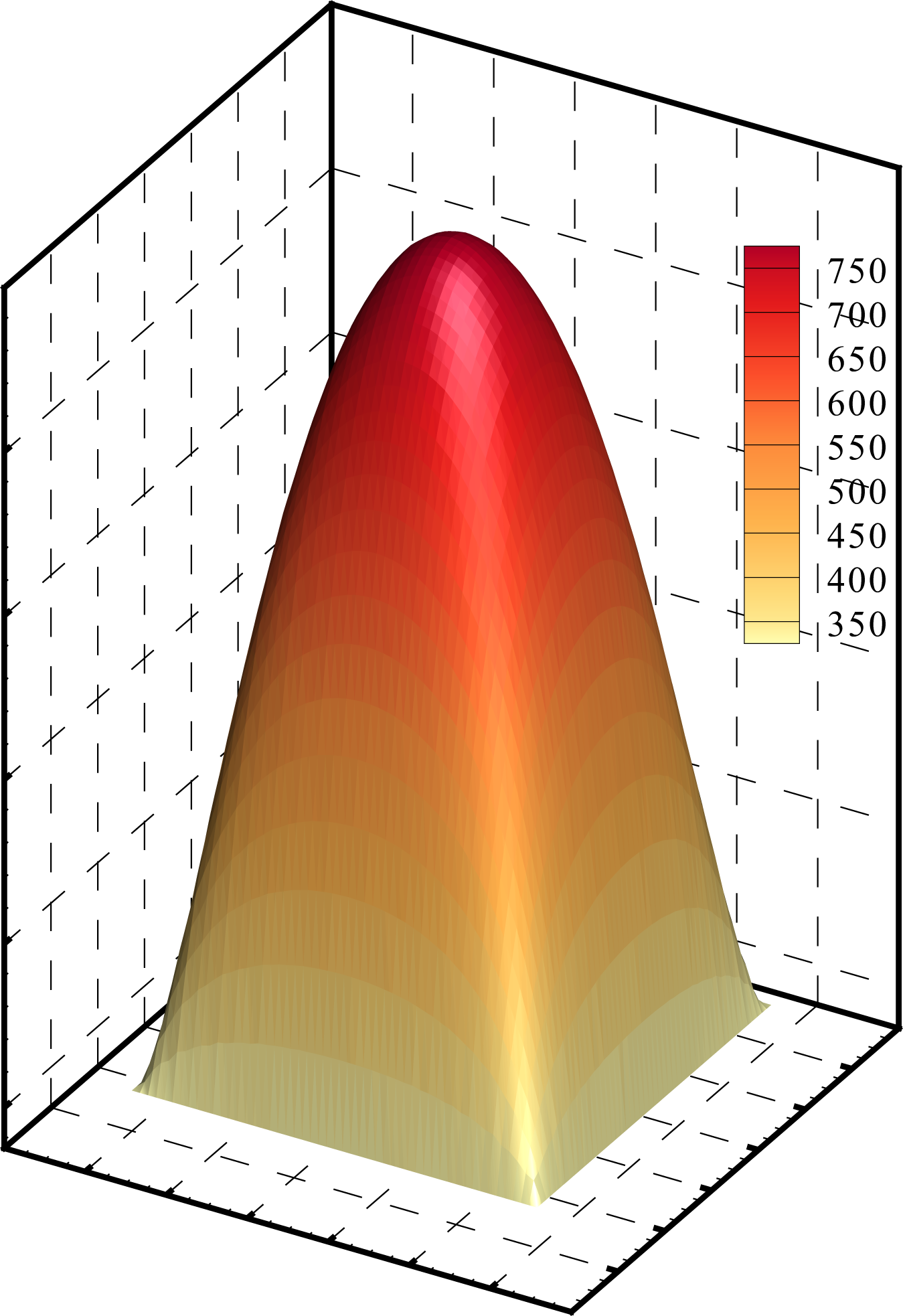}\\
  (a)
\end{minipage}
\begin{minipage}[c]{0.22\textwidth}
  \centering
  \includegraphics[width=\linewidth]{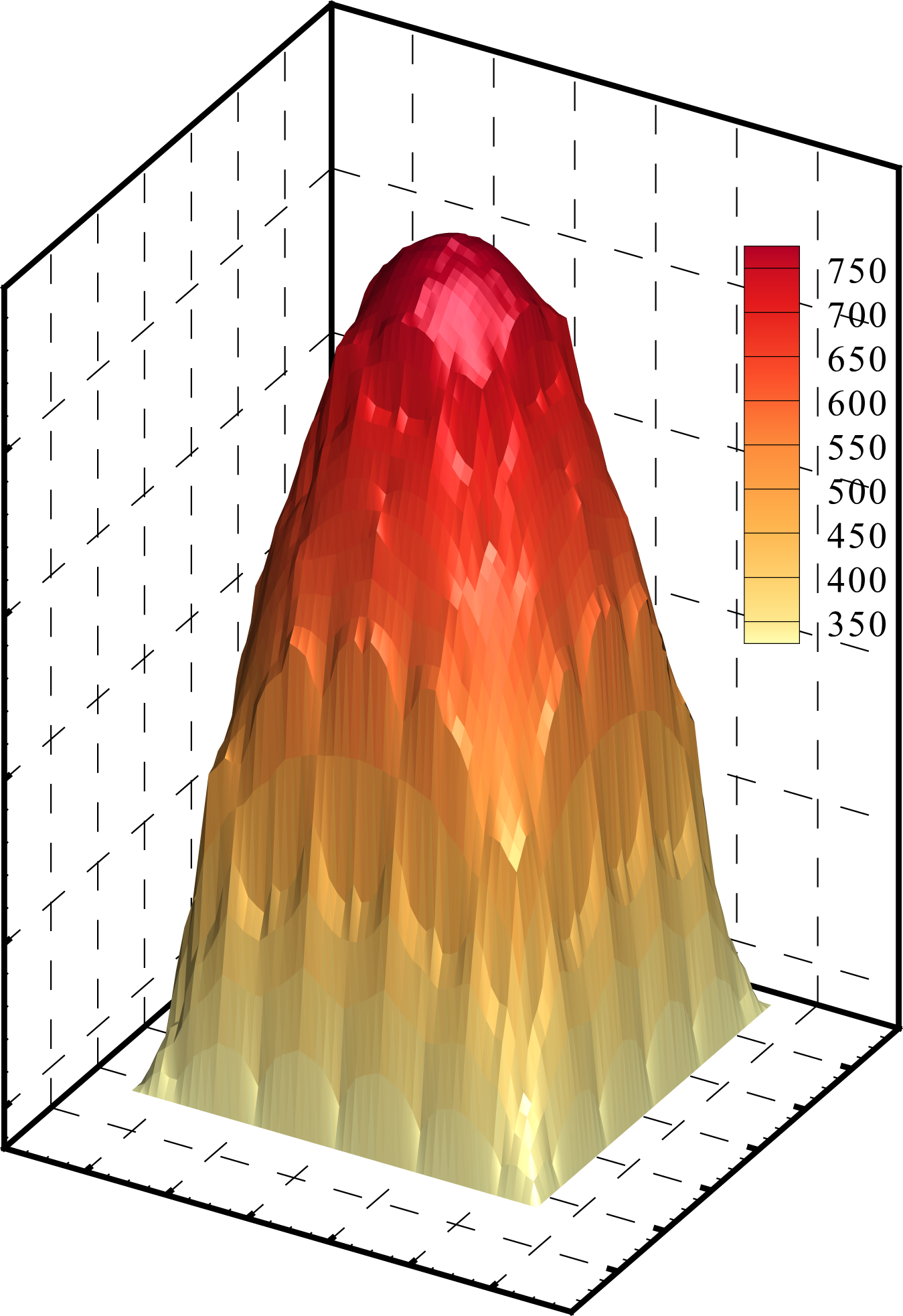}\\
  (b)
\end{minipage}
\begin{minipage}[c]{0.22\textwidth}
  \centering
  \includegraphics[width=\linewidth]{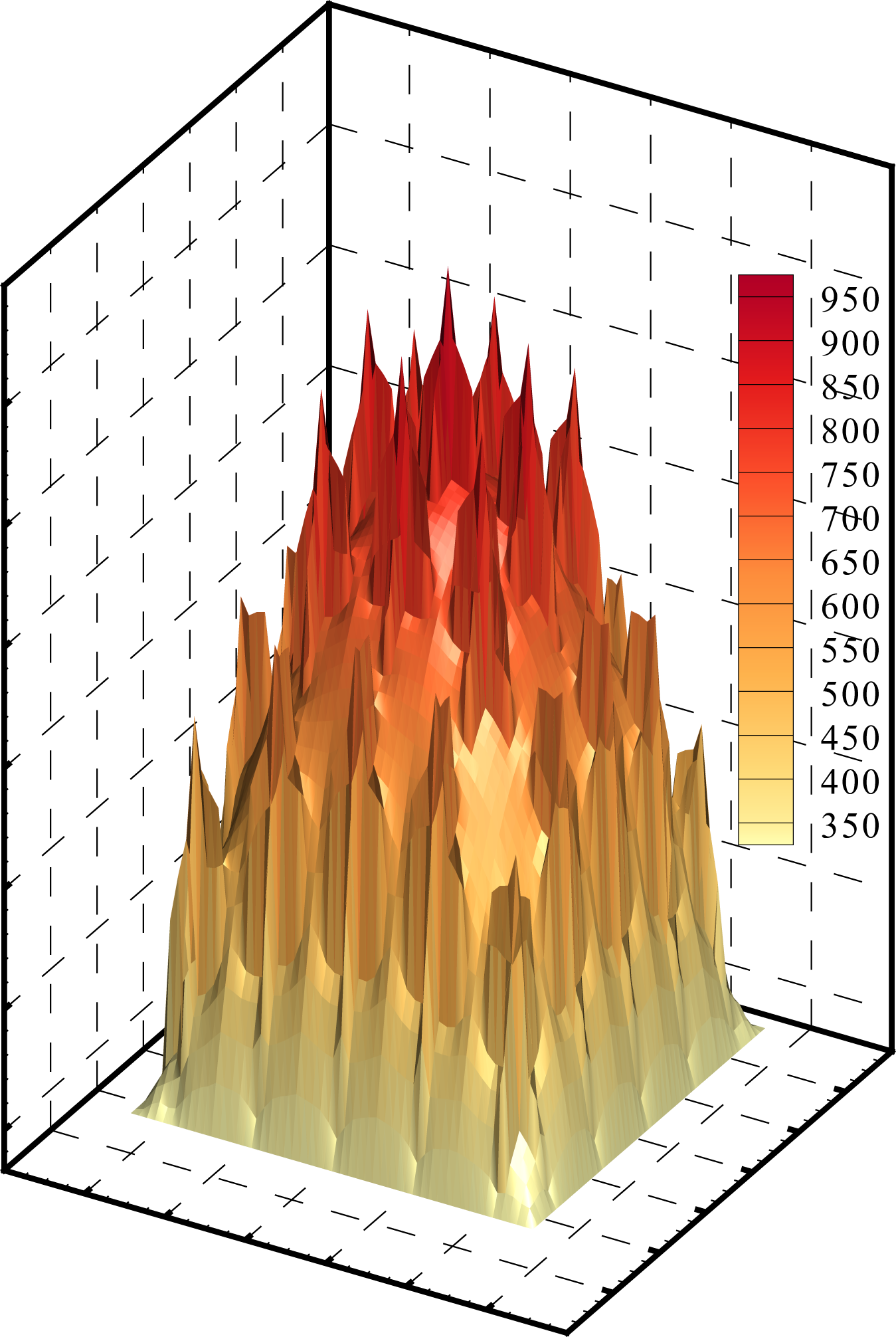}\\
  (c)
\end{minipage}
\begin{minipage}[c]{0.22\textwidth}
  \centering
  \includegraphics[width=\linewidth]{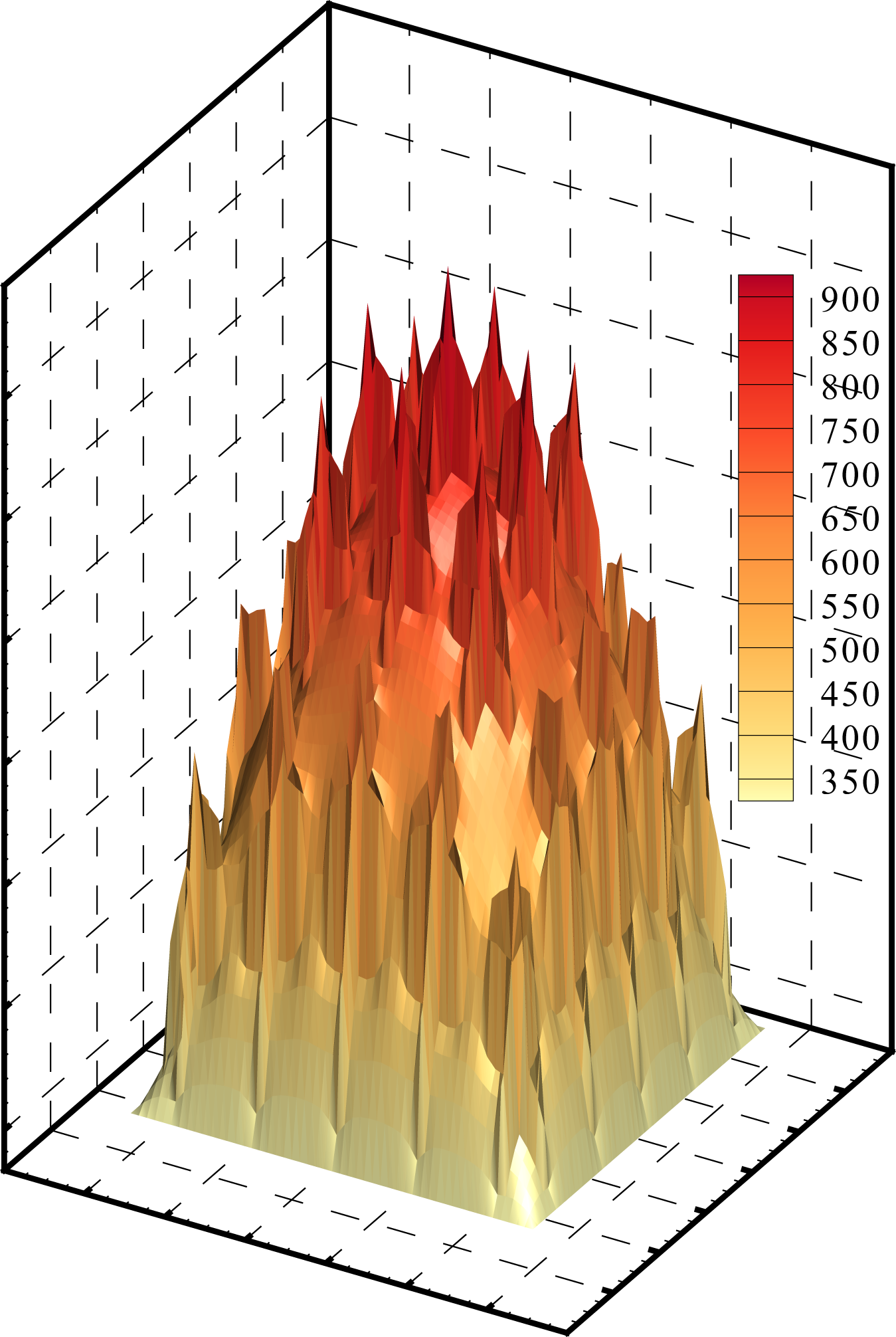}\\
  (d)
\end{minipage}

\caption{The temperature field on the cross-section $x_3=0.5$\,cm at $t=1.0\,\mathrm{s}$: (a) $T_{0}$; (b) $T^{(1{\varepsilon})}$; (c) $T^{(2{\varepsilon})}$; (d) $T^{{\varepsilon}}$.}
\label{ex43dtemperaturez05}
\end{figure}

After multi-scale simulation, \autoref{ex43dtemperatureSlice} displays the temperature distributions on
the slices $x_1=0.5$, $x_2=0.5$, and $x_3=0.5$\,cm at $t=1.0\,\mathrm{s}$, and \autoref{ex43dtemperaturez05} shows the temperature field on the cross-section $x_3=0.5$\,cm, at $t=1.0\,\mathrm{s}$.

\begin{figure}[!htb]
\centering
\begin{subfigure}[c]{0.45\textwidth}
  \centering
  \includegraphics[width=\linewidth]{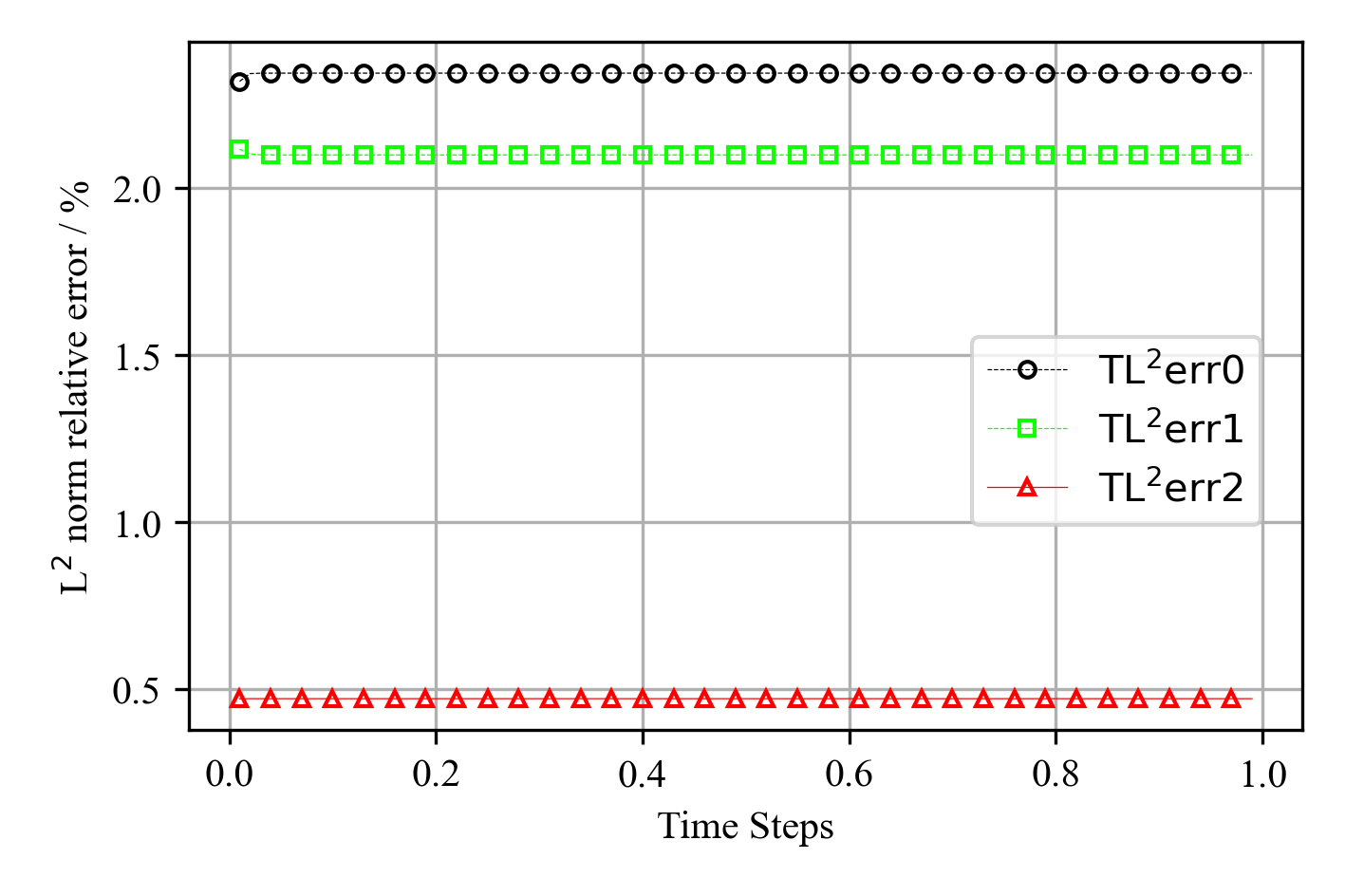}
  \caption{$L^2$ error}
\end{subfigure}
\hfill
\begin{subfigure}[c]{0.45\textwidth}
  \centering
  \includegraphics[width=\linewidth]{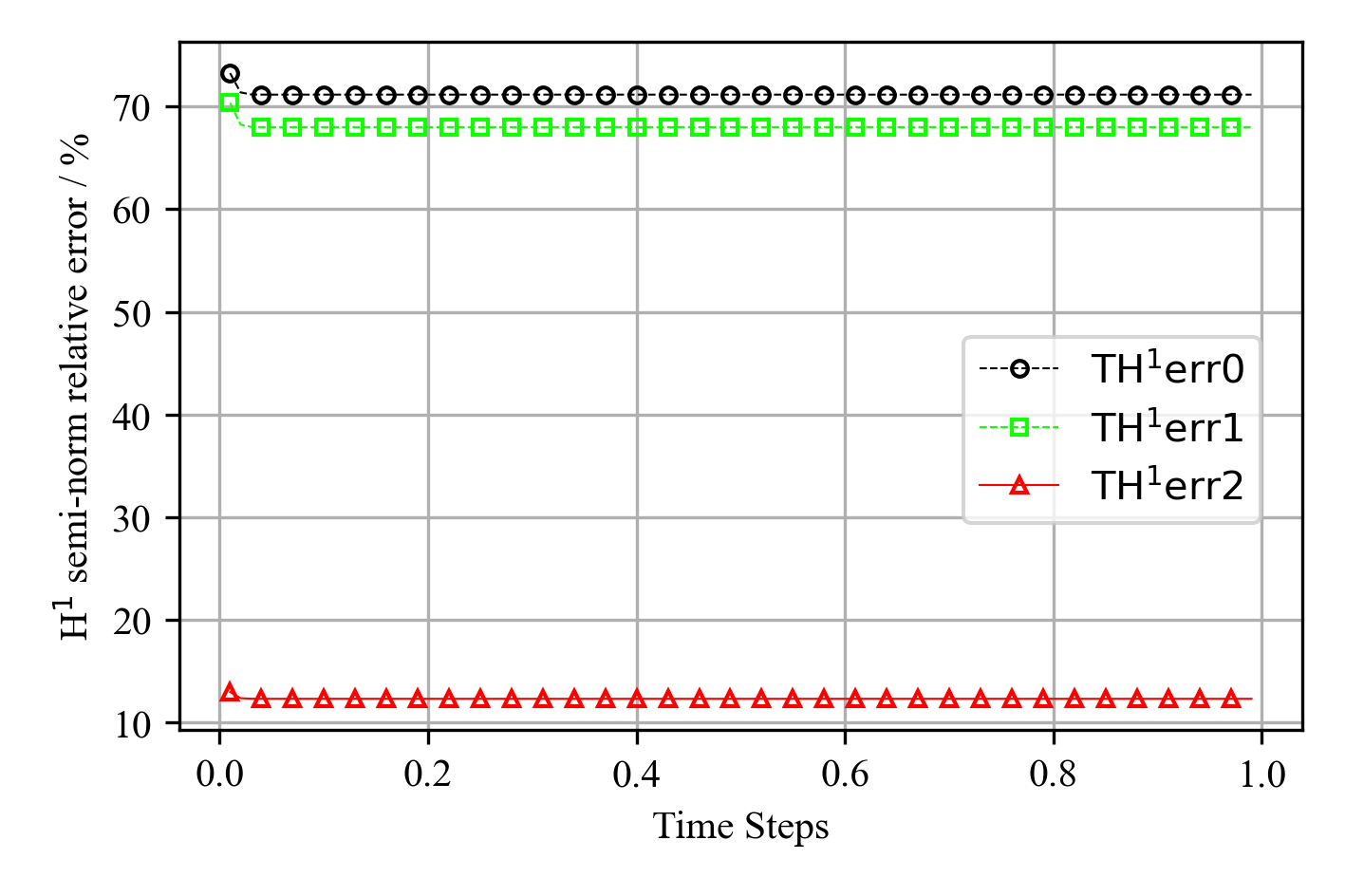}
  \caption{$H^1$ semi-norm error}
\end{subfigure}
\caption{The evolution of relative errors in the $L^2$ norm and $H^1$ semi-norm.}
\label{fig:ex43derror-comparison}
\end{figure}

Additionally, the specific relative numerical errors for the homogenization method, the low-order and high-order multi-scale approaches are presented in \autoref{fig:ex43derror-comparison}.

As illustrated by the temperature field slices (\autoref{ex43dtemperatureSlice}) and the cross-section at \(x_3=0.5\)cm (\autoref{ex43dtemperaturez05}), a consistent conclusion can be drawn: only the HOMS solution resolves the highly micro-scale oscillations. The error comparison (\autoref{fig:ex43derror-comparison}) shows that the HOMS method achieves $L^2$ norm and $H^1$ semi-norm relative errors of $0.5\%$ and $10\%$, respectively, versus $20\%$ ($L^2$ norm) and $70\%$ ($H^1$ semi-norm) for both the first-order and homogenization methods. This further confirms the effectiveness of the HOMS approach for 3D problems. In real-world engineering applications, it is impractical to obtain reference FEM solutions for large-scale composite structures. However, the HOMS method proposed herein can effectively simulate nonlinear thermal radiation problems of large-scale composite structures with minimal computational resource consumption.
\section{Conclusions}
In real-world applications, novel composite materials and structures always served in extreme high-temperature environment, which exhibit complex radiative heat transfer phenomena. These materials also possess significant temperature-dependent material properties and feature highly heterogeneous properties that vary sharply at the micro-scale. The simulation and analysis of radiative heat transfer behaviors of composite structures with temperature-dependent effect should not only take into account the overall macroscopic performances, but also consider the microscopic responses inside composite structures. Consequently, the development of robust and efficient multi-scale computational method for nonlinear multi-scale simulation of the composites is of highly practical and theoretical values.

This paper developed the HOMS computational method for nonlinear radiative heat transfer problems of the composites with heterogeneous periodic configurations. The new contributions of this study are threefold: the establishment of high-accuracy macro-micro correlated HOMS model, rigorous error analysis in both the local pointwise and global integral senses for the HOMS solutions, and corresponding multi-scale numerical algorithm with convergence analysis. Numerical experiments demonstrate that the proposed HOMS approach is efficient and robust for temperature-dependent radiative heat transfer problem \eqref{eq:all} of the composites. Furthermore, numerical results demonstrate markedly that only HOMS solutions can accurately capture the microscopic oscillating information and offer high-accuracy numerical solutions for large-scale engineering problems, which validate the theoretical results of this study.

Future research will focus on two key aspects: First, the computational efficiency will be further improved in the off-line stage via parallel computing and on-line stage by model reduction techniques. Second, the proposed multi-scale framework shall be extended to other complicated multi-physics problems with radiative thermal effect of composite structures with multiple periodic configurations.

\section*{Acknowledgments}
This work was supported by the National Natural Science Foundation of China (No.\hspace{1mm}12471387 and 12401523), Xidian University Specially Funded Project for Interdisciplinary Exploration (No.\hspace{1mm}TZJH2024008), Fundamental Research Funds for the Central Universities (No.\hspace{1mm}QTZX25082), Innovation Capability Support Program of Shaanxi Province (No.\hspace{1mm}2024RS-CXTD-88), and also supported by the Center for high performance computing of Xidian University.
\appendix
\renewcommand{\appendixname}{Appendix~\Alph{section}}

\setcounter{equation}{0}
\renewcommand{\theequation}{A.\arabic{equation}}
\section*{Appendix A. Detailed mathematical expressions of some functions}\label{appendix:A}
The detailed expressions of ${{F}_{0}}({\bm{x}},{\bm{y}},t)$, ${{F}_{1}}({\bm{x}},{\bm{y}},t)$, and ${{F}_{2}}({\bm{x}},{\bm{y}},t)$ are presented as follows.
\begin{equation}
\begin{aligned}
& F_{0}({\bm{x}},{\bm{y}},t) = \Big[ \bar S - \rho^{(0)} c^{(0)} \Big]\frac{\partial T_0}{\partial t} + \Big[ \bar\beta \sigma_{\mathrm{B}} - \beta^{(0)} \sigma_{\mathrm{B}} \Big] T_0^4 \\
&+ \Big[ k_{\alpha_1\alpha_2}^{(0)} + \frac{\partial}{\partial y_i}\big( k_{i\alpha_1}^{(0)} M_{\alpha_2} \big) + k_{\alpha_1 j}^{(0)} \frac{\partial M_{\alpha_2}}{\partial y_j} - \bar k_{\alpha_1\alpha_2} \Big]\frac{\partial^2 T_0}{\partial x_{\alpha_1}\partial x_{\alpha_2}} \\
&+ \Big[ \frac{\partial k_{i\alpha_1}^{(0)}}{\partial x_i} + \frac{\partial}{\partial y_i}\big( k_{ij}^{(0)} \frac{\partial M_{\alpha_1}}{\partial x_j} \big) + \frac{\partial}{\partial x_i}\big( k_{ij}^{(0)} \frac{\partial M_{\alpha_1}}{\partial y_j} \big) - \frac{\partial \bar k_{i\alpha_1}}{\partial x_i} \Big]\frac{\partial T_0}{\partial x_{\alpha_1}} \\
&+ \frac{\partial}{\partial y_i}\Big( M_{\alpha_1} \mathbf{D}^{(0,1)} k_{i\alpha_2}^{(0)} + M_{\alpha_1} \mathbf{D}^{(0,1)} k_{ij}^{(0)} \frac{\partial M_{\alpha_2}}{\partial y_j} \Big)\frac{\partial T_0}{\partial x_{\alpha_1}}\frac{\partial T_0}{\partial x_{\alpha_2}}.
\end{aligned}
\end{equation}

\begin{equation}
\begin{aligned}
&F_1({\bm{x}},{\bm{y}},t) =
\rho^{(0)} c^{(0)} \frac{\partial T_1}{\partial t}
- \frac{\partial}{\partial y_i} \left( k_{ij}^{(2)} \frac{\partial T_1}{\partial y_j} \right)
- \frac{\partial}{\partial x_i} \left( k_{ij}^{(1)} \frac{\partial T_1}{\partial y_j} \right) \\
&- \frac{\partial}{\partial y_i} \left( k_{ij}^{(1)} \frac{\partial T_1}{\partial x_j} \right)
- \frac{\partial}{\partial x_i} \left( k_{ij}^{(0)} \frac{\partial T_1}{\partial x_j} \right) \\
&- \sigma_{\mathrm{B}} \beta^{(0)} T_1 ( 2T_0+{\varepsilon}T_1+{\varepsilon}^2T_2 ) \left[ T_0^2 + (T_0+{\varepsilon}T_1+{\varepsilon}^2T_2)^2  \right]+ \mathrm{O}(\varepsilon).
\end{aligned}
\end{equation}

\begin{equation}
\begin{aligned}
&F_2({\bm{x}},{\bm{y}},t) =
\rho^{(0)} c^{(0)} \frac{\partial T_1}{\partial t}
- \frac{\partial}{\partial y_i} \left( k_{ij}^{(2)} \frac{\partial T_1}{\partial y_j} \right)
- \frac{\partial}{\partial x_i} \left( k_{ij}^{(1)} \frac{\partial T_1}{\partial y_j} \right)\\
&- \frac{\partial}{\partial y_i} \left[ k_{ij}^{(1)} \left( \frac{\partial T_1}{\partial x_j} + \frac{\partial T_2}{\partial y_j} \right) \right]
- \frac{\partial}{\partial y_i} \left( k_{ij}^{(0)} \frac{\partial T_2}{\partial x_j} \right)\\
&- \frac{\partial}{\partial x_i} \left[ k_{ij}^{(0)} \left( \frac{\partial T_1}{\partial x_j} + \frac{\partial T_2}{\partial y_j} \right) \right]
- \varepsilon \frac{\partial}{\partial y_i} \left( k_{ij}^{(0)} \frac{\partial T_2}{\partial y_j} \right) \\
&- \sigma_{\mathrm{B}}\beta^{(0)} (T_1 \!+\! \varepsilon T_2) ( 2T_0 \!+\! {\varepsilon}T_1 \!+\! {\varepsilon}^2T_2 ) \!\left[ T_0^2 \!+\! (T_0 \!+\! {\varepsilon}T_1 \!+\! {\varepsilon}^2T_2)^2  \right]\!+\! \mathrm{O}(\varepsilon).
\end{aligned}
\end{equation}

\bibliographystyle{model1-num-names}
\bibliography{radiation}

\begin{thebibliography}{41}
\expandafter\ifx\csname natexlab\endcsname\relax\def\natexlab#1{#1}\fi
\providecommand{\url}[1]{\texttt{#1}}
\providecommand{\href}[2]{#2}
\providecommand{\path}[1]{#1}
\providecommand{\DOIprefix}{doi:}
\providecommand{\ArXivprefix}{arXiv:}
\providecommand{\URLprefix}{URL: }
\providecommand{\Pubmedprefix}{pmid:}
\providecommand{\doi}[1]{\href{http://dx.doi.org/#1}{\path{#1}}}
\providecommand{\Pubmed}[1]{\href{pmid:#1}{\path{#1}}}
\providecommand{\bibinfo}[2]{#2}
\ifx\xfnm\relax \def\xfnm[#1]{\unskip,\space#1}\fi
\bibitem[{Tanigawa et~al.(1996)Tanigawa, Akai, Kawamura, and Oka}]{R1}
\bibinfo{author}{Y.~Tanigawa}, \bibinfo{author}{T.~Akai},
  \bibinfo{author}{R.~Kawamura}, \bibinfo{author}{N.~Oka},
\newblock \bibinfo{title}{Transient heat conduction and thermal stress problems
  of a nonhomogeneous plate with temperature-dependent material properties},
\newblock \bibinfo{journal}{Journal of Thermal Stresses} \bibinfo{volume}{19}
  (\bibinfo{year}{1996}) \bibinfo{pages}{77--102}.
\bibitem[{Saad et~al.(2023)Saad, Martinez, and Trice}]{R2}
\bibinfo{author}{A.~A. Saad}, \bibinfo{author}{C.~Martinez},
  \bibinfo{author}{R.~W. Trice},
\newblock \bibinfo{title}{Radiation heat transfer during hypersonic flight: A
  review of emissivity measurement and enhancement approaches of ultra-high
  temperature ceramics},
\newblock \bibinfo{journal}{International Journal of Ceramic Engineering \&
  Science} \bibinfo{volume}{5} (\bibinfo{year}{2023}) \bibinfo{pages}{e10171}.
\bibitem[{Modest and Mazumder(2021)}]{R3}
\bibinfo{author}{M.~F. Modest}, \bibinfo{author}{S.~Mazumder},
  \bibinfo{title}{Radiative Heat Transfer}, \bibinfo{publisher}{Academic
  Press}, \bibinfo{year}{2021}.
\bibitem[{Cioranescu and Donato(1999)}]{R4}
\bibinfo{author}{D.~Cioranescu}, \bibinfo{author}{P.~Donato},
  \bibinfo{title}{An introduction to homogenization},
  \bibinfo{publisher}{Oxford university press}, \bibinfo{year}{1999}.
\bibitem[{Babu{\v s}ka(1976)}]{R5}
\bibinfo{author}{I.~Babu{\v s}ka},
\newblock \bibinfo{title}{Solution of interface problems by homogenization.
  {I}},
\newblock \bibinfo{journal}{SIAM Journal on Mathematical Analysis}
  \bibinfo{volume}{7} (\bibinfo{year}{1976}) \bibinfo{pages}{603--634}.
\bibitem[{De~Giorgi(1977)}]{R6}
\bibinfo{author}{E.~De~Giorgi},
\newblock \bibinfo{title}{$\gamma$-convergenza e g-convergenza},
\newblock \bibinfo{journal}{Ennio De Giorgi}  (\bibinfo{year}{1977})
  \bibinfo{pages}{437}.
\bibitem[{Lewinski and Telega(2000)}]{R7}
\bibinfo{author}{T.~Lewinski}, \bibinfo{author}{J.~J. Telega},
  \bibinfo{title}{Plates, laminates and shells: asymptotic analysis and
  homogenization}, volume~\bibinfo{volume}{52}, \bibinfo{publisher}{World
  Scientific}, \bibinfo{year}{2000}.
\bibitem[{Bakhvalov and Panasenko(1989)}]{R8}
\bibinfo{author}{N.~S. Bakhvalov}, \bibinfo{author}{G.~P. Panasenko},
  \bibinfo{title}{Homogenisation: Averaging Processes in Periodic Media},
  \bibinfo{publisher}{Kluwer Academic Publishers}, \bibinfo{year}{1989}.
\bibitem[{Oleinik et~al.(1992)Oleinik, Shamaev, and Yosifian}]{R9}
\bibinfo{author}{O.~A. Oleinik}, \bibinfo{author}{A.~S. Shamaev},
  \bibinfo{author}{G.~A. Yosifian}, \bibinfo{title}{Mathematical Problems in
  Elasticity and Homogenization}, \bibinfo{publisher}{North-Holland},
  \bibinfo{year}{1992}.
\bibitem[{Bensoussan et~al.(2011)Bensoussan, Lions, and Papanicolaou}]{R10}
\bibinfo{author}{A.~Bensoussan}, \bibinfo{author}{J.-L. Lions},
  \bibinfo{author}{G.~Papanicolaou}, \bibinfo{title}{Asymptotic analysis for
  periodic structures}, volume \bibinfo{volume}{374},
  \bibinfo{publisher}{American Mathematical Soc.}, \bibinfo{address}{Rhode
  Island}, \bibinfo{year}{2011}.
\bibitem[{Hou et~al.(1999)Hou, Wu, and Cai}]{R11}
\bibinfo{author}{T.~Y. Hou}, \bibinfo{author}{X.-H. Wu},
  \bibinfo{author}{Z.~Cai},
\newblock \bibinfo{title}{Convergence of a multiscale finite element method for
  elliptic problems with rapidly oscillating coefficients},
\newblock \bibinfo{journal}{Mathematics of Computation} \bibinfo{volume}{68}
  (\bibinfo{year}{1999}) \bibinfo{pages}{913--943}.
\bibitem[{E and Engquist(2003)}]{R12}
\bibinfo{author}{W.~E}, \bibinfo{author}{B.~Engquist},
\newblock \bibinfo{title}{The heterogeneous multiscale methods},
\newblock \bibinfo{journal}{Communications in Mathematical Sciences}
  \bibinfo{volume}{1} (\bibinfo{year}{2003}) \bibinfo{pages}{87--132}.
\bibitem[{Ming and Zhang(2005)}]{R13}
\bibinfo{author}{P.~Ming}, \bibinfo{author}{P.~Zhang},
\newblock \bibinfo{title}{Analysis of the heterogeneous multiscale method for
  elliptic homogenization problems},
\newblock \bibinfo{journal}{Journal of the American Mathematical Society}
  \bibinfo{volume}{18} (\bibinfo{year}{2005}) \bibinfo{pages}{121--156}.
\bibitem[{Hughes et~al.(1998)Hughes, Feij{\'o}o, Mazzei, and Quincy}]{R14}
\bibinfo{author}{T.~J.~R. Hughes}, \bibinfo{author}{G.~R. Feij{\'o}o},
  \bibinfo{author}{L.~Mazzei}, \bibinfo{author}{J.-B. Quincy},
\newblock \bibinfo{title}{The variational multiscale method---a paradigm for
  computational mechanics},
\newblock \bibinfo{journal}{Computer Methods in Applied Mechanics and
  Engineering} \bibinfo{volume}{166} (\bibinfo{year}{1998})
  \bibinfo{pages}{3--24}.
\bibitem[{Xing et~al.(2010)Xing, Yang, and Wang}]{R15}
\bibinfo{author}{Y.~F. Xing}, \bibinfo{author}{Y.~Yang}, \bibinfo{author}{X.~M.
  Wang},
\newblock \bibinfo{title}{A multiscale eigenelement method and its application
  to periodical composite structures},
\newblock \bibinfo{journal}{Composite Structures} \bibinfo{volume}{92}
  (\bibinfo{year}{2010}) \bibinfo{pages}{2265--2275}.
\bibitem[{Henning and M{\aa}lqvist(2014)}]{R16}
\bibinfo{author}{P.~Henning}, \bibinfo{author}{A.~M{\aa}lqvist},
\newblock \bibinfo{title}{Localized orthogonal decomposition techniques for
  boundary value problems},
\newblock \bibinfo{journal}{SIAM Journal on Scientific Computing}
  \bibinfo{volume}{36} (\bibinfo{year}{2014}) \bibinfo{pages}{A1609--A1634}.
\bibitem[{He and Pindera(2021)}]{R17}
\bibinfo{author}{Z.~He}, \bibinfo{author}{M.-J. Pindera},
\newblock \bibinfo{title}{Finite volume based asymptotic homogenization theory
  for periodic materials under anti-plane shear},
\newblock \bibinfo{journal}{European Journal of Mechanics - A/Solids}
  \bibinfo{volume}{85} (\bibinfo{year}{2021}) \bibinfo{pages}{104122}.
\bibitem[{Ma et~al.(2023)Ma, Alber, and Scheichl}]{R18}
\bibinfo{author}{C.~Ma}, \bibinfo{author}{C.~Alber},
  \bibinfo{author}{R.~Scheichl},
\newblock \bibinfo{title}{Wavenumber explicit convergence of a multiscale
  generalized finite element method for heterogeneous helmholtz problems},
\newblock \bibinfo{journal}{SIAM Journal on Numerical Analysis}
  \bibinfo{volume}{61} (\bibinfo{year}{2023}) \bibinfo{pages}{1546--1584}.
\bibitem[{Wu et~al.(2014)Wu, Nie, and Yang}]{R19}
\bibinfo{author}{Y.~T. Wu}, \bibinfo{author}{Y.~F. Nie}, \bibinfo{author}{Z.~H.
  Yang},
\newblock \bibinfo{title}{Comparison of four multiscale methods for elliptic
  problems},
\newblock \bibinfo{journal}{CMES-Computer Modeling in Engineering \& Sciences}
  \bibinfo{volume}{99} (\bibinfo{year}{2014}) \bibinfo{pages}{297--325}.
\bibitem[{Gao et~al.(2020)Gao, Xing, Huang, Li, and Yang}]{R20}
\bibinfo{author}{Y.~Gao}, \bibinfo{author}{Y.~Xing},
  \bibinfo{author}{Z.~Huang}, \bibinfo{author}{M.~Li},
  \bibinfo{author}{Y.~Yang},
\newblock \bibinfo{title}{An assessment of multiscale asymptotic expansion
  method for linear static problems of periodic composite structures},
\newblock \bibinfo{journal}{European Journal of Mechanics - A/Solids}
  \bibinfo{volume}{81} (\bibinfo{year}{2020}) \bibinfo{pages}{103951}.
\bibitem[{Cao et~al.(2002)Cao, Cui, and Zhu}]{R21}
\bibinfo{author}{L.-Q. Cao}, \bibinfo{author}{J.-Z. Cui},
  \bibinfo{author}{D.-C. Zhu},
\newblock \bibinfo{title}{Multiscale asymptotic analysis and numerical
  simulation for the second order {Helmholtz} equations with rapidly
  oscillating coefficients over general convex domains},
\newblock \bibinfo{journal}{SIAM Journal on Numerical Analysis}
  \bibinfo{volume}{40} (\bibinfo{year}{2002}) \bibinfo{pages}{543--577}.
\bibitem[{Dong et~al.(2023)Dong, Cui, Nie, Ma, Jin, and Huang}]{R22}
\bibinfo{author}{H.~Dong}, \bibinfo{author}{J.~Cui}, \bibinfo{author}{Y.~Nie},
  \bibinfo{author}{R.~Ma}, \bibinfo{author}{K.~Jin},
  \bibinfo{author}{D.~Huang},
\newblock \bibinfo{title}{Multi-scale computational method for nonlinear
  dynamic thermo-mechanical problems of composite materials with
  temperature-dependent properties},
\newblock \bibinfo{journal}{Communications in Nonlinear Science and Numerical
  Simulation} \bibinfo{volume}{118} (\bibinfo{year}{2023})
  \bibinfo{pages}{107000}.
\bibitem[{Dong et~al.(2025)Dong, Guan, and Nie}]{R23}
\bibinfo{author}{H.~Dong}, \bibinfo{author}{X.~Guan}, \bibinfo{author}{Y.~Nie},
\newblock \bibinfo{title}{Multiscale method and convergence analysis for
  coupled nonlinear thermomechanical problems in heterogeneous shells},
\newblock \bibinfo{journal}{SIAM Journal on Scientific Computing}
  \bibinfo{volume}{47} (\bibinfo{year}{2025}) \bibinfo{pages}{B190--B219}.
\bibitem[{Cao and Cui(2004)}]{R24}
\bibinfo{author}{L.-Q. Cao}, \bibinfo{author}{J.-Z. Cui},
\newblock \bibinfo{title}{Asymptotic expansions and numerical algorithms of
  eigenvalues and eigenfunctions of the {Dirichlet} problem for second order
  elliptic equations in perforated domains},
\newblock \bibinfo{journal}{Numerische Mathematik} \bibinfo{volume}{96}
  (\bibinfo{year}{2004}) \bibinfo{pages}{525--581}.
\bibitem[{Feng and Cui(2004)}]{R25}
\bibinfo{author}{Y.-P. Feng}, \bibinfo{author}{J.-Z. Cui},
\newblock \bibinfo{title}{Multi-scale analysis and {FE} computation for the
  structure of composite materials with small periodic configuration under
  condition of coupled thermoelasticity},
\newblock \bibinfo{journal}{International Journal for Numerical Methods in
  Engineering} \bibinfo{volume}{60} (\bibinfo{year}{2004})
  \bibinfo{pages}{1879--1910}.
\bibitem[{Wang et~al.(2015)Wang, Cao, and Wong}]{R26}
\bibinfo{author}{X.~Wang}, \bibinfo{author}{L.-Q. Cao}, \bibinfo{author}{Y.~S.
  Wong},
\newblock \bibinfo{title}{Multiscale computation and convergence for coupled
  thermoelastic system in composite materials},
\newblock \bibinfo{journal}{Multiscale Modeling \& Simulation}
  \bibinfo{volume}{13} (\bibinfo{year}{2015}) \bibinfo{pages}{661--690}.
\bibitem[{Ma et~al.(2016)Ma, Cui, Li, and Wang}]{R27}
\bibinfo{author}{Q.~Ma}, \bibinfo{author}{J.~Cui}, \bibinfo{author}{Z.~Li},
  \bibinfo{author}{Z.~Wang},
\newblock \bibinfo{title}{Second-order asymptotic algorithm for heat conduction
  problems of periodic composite materials in curvilinear coordinates},
\newblock \bibinfo{journal}{Journal of Computational and Applied Mathematics}
  \bibinfo{volume}{306} (\bibinfo{year}{2016}) \bibinfo{pages}{87--115}.
\bibitem[{Dong et~al.(2017)Dong, Cui, Nie, and Yang}]{R28}
\bibinfo{author}{H.~Dong}, \bibinfo{author}{J.~Cui}, \bibinfo{author}{Y.~Nie},
  \bibinfo{author}{Z.~Yang},
\newblock \bibinfo{title}{Second-order two-scale computational method for
  nonlinear dynamic thermo-mechanical problems of composites with cylindrical
  periodicity},
\newblock \bibinfo{journal}{Communications in Computational Physics}
  \bibinfo{volume}{21} (\bibinfo{year}{2017}) \bibinfo{pages}{1173--1206}.
\bibitem[{Amosov(2010)}]{R29}
\bibinfo{author}{A.~A. Amosov},
\newblock \bibinfo{title}{Nonstationary radiative--conductive heat transfer
  problem in a periodic system of grey heat shields},
\newblock \bibinfo{journal}{Journal of Mathematical Sciences}
  \bibinfo{volume}{169} (\bibinfo{year}{2010}) \bibinfo{pages}{1--45}.
\bibitem[{Amosov(2011)}]{R30}
\bibinfo{author}{A.~A. Amosov},
\newblock \bibinfo{title}{Semidiscrete and asymptotic approximations for the
  nonstationary radiative--conductive heat transfer problem in a periodic
  system of grey heat shields},
\newblock \bibinfo{journal}{Journal of Mathematical Sciences}
  \bibinfo{volume}{176} (\bibinfo{year}{2011}) \bibinfo{pages}{361--408}.
\bibitem[{Bakhvalov(1981)}]{R31}
\bibinfo{author}{N.~S. Bakhvalov},
\newblock \bibinfo{title}{Averaging of the heat-transfer process in periodic
  media with radiation},
\newblock \bibinfo{journal}{Differential Equations} \bibinfo{volume}{17}
  (\bibinfo{year}{1981}) \bibinfo{pages}{1094--1100}.
\bibitem[{Allaire and El~Ganaoui(2009)}]{R32}
\bibinfo{author}{G.~Allaire}, \bibinfo{author}{K.~El~Ganaoui},
\newblock \bibinfo{title}{Homogenization of a conductive and radiative heat
  transfer problem},
\newblock \bibinfo{journal}{Multiscale Modeling \& Simulation}
  \bibinfo{volume}{7} (\bibinfo{year}{2009}) \bibinfo{pages}{1148--1170}.
\bibitem[{Allaire and Habibi(2013{\natexlab{a}})}]{R33}
\bibinfo{author}{G.~Allaire}, \bibinfo{author}{Z.~Habibi},
\newblock \bibinfo{title}{Homogenization of a conductive, convective, and
  radiative heat transfer problem in a heterogeneous domain},
\newblock \bibinfo{journal}{SIAM Journal on Mathematical Analysis}
  \bibinfo{volume}{45} (\bibinfo{year}{2013}{\natexlab{a}})
  \bibinfo{pages}{1136--1178}.
\bibitem[{Allaire and Habibi(2013{\natexlab{b}})}]{R34}
\bibinfo{author}{G.~Allaire}, \bibinfo{author}{Z.~Habibi},
\newblock \bibinfo{title}{Second order corrector in the homogenization of a
  conductive-radiative heat transfer problem},
\newblock \bibinfo{journal}{Discrete and Continuous Dynamical Systems - Series
  B} \bibinfo{volume}{18} (\bibinfo{year}{2013}{\natexlab{b}})
  \bibinfo{pages}{1--36}.
\bibitem[{Huang and Cao(2014)}]{R35}
\bibinfo{author}{J.~Huang}, \bibinfo{author}{L.~Cao},
\newblock \bibinfo{title}{Global regularity and multiscale approach for thermal
  radiation heat transfer},
\newblock \bibinfo{journal}{Multiscale Modeling \& Simulation}
  \bibinfo{volume}{12} (\bibinfo{year}{2014}) \bibinfo{pages}{694--724}.
\bibitem[{Huang et~al.(2015)Huang, Cao, and Yang}]{R36}
\bibinfo{author}{J.~Huang}, \bibinfo{author}{L.~Cao},
  \bibinfo{author}{C.~Yang},
\newblock \bibinfo{title}{A multiscale algorithm for radiative heat transfer
  equation with rapidly oscillating coefficients},
\newblock \bibinfo{journal}{Applied Mathematics and Computation}
  \bibinfo{volume}{266} (\bibinfo{year}{2015}) \bibinfo{pages}{149--168}.
\bibitem[{Han et~al.(2018)Han, Nie, and Dong}]{R37}
\bibinfo{author}{Y.~Han}, \bibinfo{author}{Y.~Nie}, \bibinfo{author}{H.~Dong},
\newblock \bibinfo{title}{A fast multipole algorithm for radiative heat
  transfer in {3D} semitransparent media},
\newblock \bibinfo{journal}{Journal of Quantitative Spectroscopy and Radiative
  Transfer} \bibinfo{volume}{221} (\bibinfo{year}{2018})
  \bibinfo{pages}{8--17}.
\bibitem[{Li et~al.(2014)Li, Gao, and Sun}]{R38}
\bibinfo{author}{B.~Li}, \bibinfo{author}{H.~Gao}, \bibinfo{author}{W.~Sun},
\newblock \bibinfo{title}{Unconditionally optimal error estimates of a
  {Crank--Nicolson} {Galerkin} method for the nonlinear thermistor equations},
\newblock \bibinfo{journal}{SIAM Journal on Numerical Analysis}
  \bibinfo{volume}{52} (\bibinfo{year}{2014}) \bibinfo{pages}{933--954}.
\bibitem[{Dong and Cao(2009)}]{R39}
\bibinfo{author}{Q.-L. Dong}, \bibinfo{author}{L.-Q. Cao},
\newblock \bibinfo{title}{Multiscale asymptotic expansions and numerical
  algorithms for the wave equations of second order with rapidly oscillating
  coefficients},
\newblock \bibinfo{journal}{Applied Numerical Mathematics} \bibinfo{volume}{59}
  (\bibinfo{year}{2009}) \bibinfo{pages}{3008--3032}.
\bibitem[{Dong and Cao(2014)}]{R40}
\bibinfo{author}{Q.-L. Dong}, \bibinfo{author}{L.-Q. Cao},
\newblock \bibinfo{title}{Multiscale asymptotic expansions methods and
  numerical algorithms for the wave equations in perforated domains},
\newblock \bibinfo{journal}{Applied Mathematics and Computation}
  \bibinfo{volume}{232} (\bibinfo{year}{2014}) \bibinfo{pages}{872--887}.
\bibitem[{Cao(2006)}]{R41}
\bibinfo{author}{L.-Q. Cao},
\newblock \bibinfo{title}{Multiscale asymptotic expansion and finite element
  methods for the mixed boundary value problems of second order elliptic
  equation in perforated domains},
\newblock \bibinfo{journal}{Numerische Mathematik} \bibinfo{volume}{103}
  (\bibinfo{year}{2006}) \bibinfo{pages}{11--45}.

\end{thebibliography}







\end{document}